\documentclass{article}

\usepackage{latexsym}
\usepackage{a4wide}
\usepackage{amscd}
\usepackage{graphics}
\usepackage{amsmath}
\usepackage{amssymb}
\usepackage{amsthm}
\input xy
\xyoption{all}

\newcommand{\N}{\mathbb{N}}                     
\newcommand{\Z}{\mathbb{Z}}                     
\newcommand{\R}{\mathbb{R}}                     
\newcommand{\C}{\mathbb{C}}                     
\newcommand{\T}{\mathbb{T}}                     
\newcommand{\set}[2]{\left\{{#1}\mid{#2}\right\}}       
\newcommand{\re}{\mathrm{Re\,}}                 
\newcommand{\dist}{\mathrm{dist\,}}             
\newcommand{\coker}{\mathrm{coker\,}}           
\newcommand{\Span}{\mathrm{span\,}}             
\newcommand{\ind}{\mathrm{ind\,}}               
\newcommand{\codim}{\mathrm{codim}}           
\newcommand{\graf}{\mathrm{graph\,}}       
\newcommand{\rank}{\mathrm{rank\,}}     
\newcommand{\ran}{\mathrm{ran\,}}   
\newcommand{\dom}{\mathrm{dom}\,}       

\newcommand{\crit}{\mathrm{crit}}

\newcommand{\grad}{\mathrm{grad\,}}

\newcommand{\hess}{\mathrm{Hess\,}}

\newcommand{\co}{\mathrm{co}\,}             
\newcommand{\Det}{\mathrm{Det}}
\newcommand{\rest}{\mathrm{rest}\,}

\swapnumbers
\newtheorem{thm}{\sc Theorem}[section]      

\newtheorem{cor}[thm]{\sc Corollary}        
\newtheorem{lem}[thm]{\sc Lemma}            
\newtheorem{prop}[thm]{\sc Proposition}     

\newtheorem{defn}[thm]{\sc Definition}      
\newtheorem{rem}[thm]{\sc Remark}
\newtheorem{ex}[thm]{\sc Example}

\title{A Morse complex for infinite dimensional manifolds\\ Part I}
\author{Alberto Abbondandolo\footnote{Scuola Normale Superiore di
    Pisa, Piazza dei Cavalieri 7,
  56126 Pisa, Italy, e-mail: abbo@sns.it.}\hspace{5pt} and Pietro 
Majer\footnote{
Dipartimento di Matematica ``L. Tonelli'', via
Buonarroti 2, 56127 Pisa, Italy, e-mail: majer@dm.unipi.it}} 
\date{December 21, 2004}

\begin{document}

\renewcommand{\theenumi}{\roman{enumi}}
\renewcommand{\labelenumi}{(\theenumi)}

\maketitle

\begin{abstract}
In this paper and in the forthcoming Part II we introduce a Morse
complex for a class of functions $f$ defined on an infinite
dimensional Hilbert manifold $M$, possibly having critical points
of infinite Morse index and co-index. The idea is to consider an
infinite dimensional subbundle - or more generally an essential
subbundle - of the tangent bundle of $M$, suitably related with
the gradient flow of $f$. This Part I deals with the following
questions about the intersection $W$ of the unstable manifold of a
critical point $x$ and the stable manifold of another critical
point $y$: finite dimensionality of $W$, possibility that
different components of $W$ have different dimension,
orientability of $W$ and coherence in the choice of an
orientation, compactness of the closure of $W$, classification, up
to topological conjugacy, of the gradient flow on the closure of
$W$, in the case $\dim W=2$.
\end{abstract}

\section*{Introduction}

Morse theory \cite{mor25} relates the topology of a compact
differentiable manifold $M$ to the combinatorics of the critical
points of a smooth Morse function $f:M\rightarrow \R$: if
$\beta_q(M) = \rank H_q(M)$ denotes the $q$-th Betti number of
$M$, and $c_q(f)$ is the number of critical points $x$ of $f$
with Morse index $m(x)=q$, then the identity
\begin{equation}
\label{mr} \sum_{q=0}^{\dim M} c_q(f) t^q = \sum_{q=0}^{\dim M}
\beta_q(M) t^q + (1+t) Q(t),
\end{equation}
holds, with $Q$ a polynomial with positive integer coefficients.
Denoting by $C_q(f)$ the free Abelian group generated by the
critical points of $f$ of index $q$, $q=0,1,\dots,\dim M$, it is
readily seen that (\ref{mr}) is implied\footnote{The two facts
would actually be equivalent if we were using coefficients in a
field, instead of the ring $\Z$.} by the existence of
homomorphisms $\partial_q:C_q(f) \rightarrow C_{q-1}(f)$ making
$\{C_*(f),\partial_*\}$ a chain complex, whose homology groups are
isomorphic to the singular $\Z$-homology groups of $M$:
\begin{equation}
\label{iso} H_q(\{C_*(f),\partial_*\}) = \frac{\ker
\partial_q}{\ran \partial_{q+1}} \cong H_q(M).
\end{equation}
A chain complex with the above properties is indeed provided by a
suitable cellular filtration of $M$. More precisely, if we fix a
Riemannian structure on $M$ such that the corresponding gradient
flow of $f$, i.e.\ the integral flow $\phi:\R \times M \rightarrow
M$ of the vector field $-\grad f$, is Morse-Smale\footnote{Here
one needs just that the unstable manifold $W^u(x)$
  and the stable manifold $W^s(y)$ have empty intersection, for every pair of
  distinct critical points $x,y$ with $m(x)\leq m(y)$.}, then the open
subsets
\[
M^q := \bigcup_{\substack{x\in \crit (f) \\ m(x)\leq q}}
\phi([0,+\infty[
  \times U_x), \quad q=0,1,\dots,\dim M,
\]
for $U_x$ a suitable small neighborhood of $x$, constitute a
cellular filtration of $M$, such that
\[
H_q(M^q,M^{q-1}) \cong C_q(f).
\]
So we get the boundary homomorphism
\begin{equation}
\label{cc}
\partial_q : C_q(f) \cong H_q(M^q,M^{q-1}) \rightarrow
H_{q-1}(M^{q-1},M^{q-2}) \cong C_{q-1}(f),
\end{equation}
and the classical isomorphism between the homology of the cellular chain
complex (\ref{cc}) and the singular homology of $M$ (see
\cite{dol80}, section V.1) implies (\ref{iso}).

The boundary homomorphism $\partial_q$ constructed above has also
the following combinatorial description, in terms of the
intersections between the unstable manifolds $W^u(x)$ and the
stable manifolds $W^s(y)$ of pairs of critical
points\footnote{Here one needs that $W^u(x)$ and $W^s(y)$
  meet transversally just when $m(x)-m(y)\leq 1$.}.
Since $\dim W^u(x) = m(x)$ and $\dim W^s(y)=\dim M - m(y)$, the
intersection $W^u(x) \cap W^s(y)$ is a submanifold of dimension
$m(x) - m(y)$. An arbitrary choice of an orientation for each
unstable manifold $W^u(x)$ determines a co-orientation (i.e.\ an
orientation of the normal bundle) for each stable manifold
$W^s(x)$, and thus an orientation for each intersection\footnote{
Indeed by transversality, a normal bundle of $W^u(x)\cap W^s(y)$
in $W^u(x)$ is also the restriction of a normal bundle of $W^s(y)$
in $M$, so it is oriented, and together with the orientation of
$W^u(x)$, it determines an orientation of $W^u(x)\cap W^s(y)$.
Notice that the
  manifold $M$ needs not be orientable.} $W^u(x)\cap W^s(y)$.
When $m(x)=q$ and $m(y)=q-1$, $W^u(x)\cap W^s(y)$ consists of
finitely many gradient flow lines, each of which can be counted as
$+1$ or as $-1$, depending on whether its orientation agrees with
the direction of the gradient flow or not. The algebraic sum of
these numbers gives an integer $n(x,y)$, and
$\partial_q$ can be expressed in terms of the generators of
$C_q(f)$ and $C_{q-1}(f)$ as
\begin{equation}
\label{fdb}
\partial_q x = \sum_{\substack{y\in \crit (f) \\ m(y)=q-1}}
n(x,y) \, y, \quad \mbox{for } x\in \crit (f), \;
m(x)=q.
\end{equation}
The Morse complex $\{C_*(f),\partial_*\}$ depends on the choice of
the Riemannian structure on $M$ (a different Riemannian structure
would produce a different gradient flow) and on the choice of the
orientations of the unstable manifolds, but the isomorphism class
of such a chain complex depends just on the function $f$.

The approach described above was essentially clear to the pioneers
of Morse theory, such as Thom \cite{tho49} and Milnor
\cite{mil63,mil65b}, and to people in dynamical systems, such as
Smale \cite{sma60,sma61,sma67} and Franks \cite{fra79,fra80}, but
it has received increasing attention after the works of Witten
\cite{wit82} and Floer \cite{flo89a}. See the systematic study by
Schwarz \cite{sch93}, and Weber's thesis \cite{web93}.
The observation on the invariance of the
isomorphism class of the Morse complex is due to Cornea and
Ranicki \cite{cr03}, together with more striking rigidity results.

Already in the sixties, Morse theory had been generalized to
infinite dimensional Hilbert manifolds (manifolds modeled on a
Hilbert space) by Palais \cite{pal63}, and Smale
\cite{sma64,sma64b}, and had been successfully applied to many
variational problems (see the expository papers of Bott \cite{bot82,bot88}, 
the books of Klingenberg
\cite{kli78,kli82}, of Mawhin and Willem \cite{mw89}, of Chang
\cite{cha93}, and references therein). Indeed, the compactness of
$M$ can be replaced by a compactness assumption on $f$, the well
known Palais-Smale condition ((PS) for short): any sequence
$(p_n)\subset M$ such that $f(p_n)$ is bounded and $\|Df(p_n)\|$
is infinitesimal must be compact. If $M$ is a Hilbert manifold
endowed with a complete Riemannian structure, and $f\in C^2(M,\R)$
is a Morse function, bounded below and satisfying (PS), then the
Morse relations (\ref{mr}) still hold, the difference being that
now (\ref{mr}) is an equality between formal power series, with
coefficients in $\N\cup \{\infty\}$.

However, (\ref{mr}) takes into account only critical points with
finite Morse index, the ultimate reason being that the closed ball
of an infinite dimensional Hilbert space is retractable onto its
boundary, so that critical points with infinite Morse index are
invisible to homotopy theory. It was Floer
\cite{flo88b,flo88a,flo88d,flo89a} who observed that the Morse
complex approach is suitable to deal with critical points of
infinite Morse index and co-index: even if the unstable and stable
manifolds are infinite dimensional, one may still hope the
dimension of their intersection to be finite. In this case, one
could try to see (\ref{fdb}) not as a description, but rather as
the definition of a chain complex. In this way, Floer was able to
develop the analogue of Morse theory in a case where the gradient
flow ODE is replaced by a Cauchy-Riemann type PDE, which does not
even determine a local flow (so that there are no stable and
unstable manifolds). The resulting theory, known as Floer
homology, plays now a central role in symplectic geometry (see
\cite{hz94,sal99} and references therein).

In the present paper, and in the forthcoming Part II, we introduce
and study the Morse complex for gradient-like flows on infinite
dimensional Hilbert manifolds. The
results we present are a far reaching generalization of a previous
work on a special class of functionals on Hilbert spaces
\cite{ama01}. See also \cite{av99} for a construction of the Morse
complex for the energy functional of an elliptic system, and
Chapter 6 in Jost's book \cite{jos02} for a general approach to
the Morse complex. More precisely, we give an answer to the
following questions.

\begin{enumerate}

\item When is $W^u(x)\cap W^s(y)$ a finite dimensional manifold?

\item How can we give coherent orientations to the manifolds
  $W^u(x)\cap W^s(y)$?

\item When is the closure of $W^u(x)\cap W^s(y)$ compact?

\item Having defined $\partial_q$ by (\ref{fdb}), how do we prove
that
  $\partial_{q-1}\circ \partial_q=0$?

\item Which form of transversality is generic?

\item How do we recover the classical infinite dimensional Morse
theory?

\item How can we compute the homology of the Morse complex?

\end{enumerate}

In the present paper we address questions (i), (ii), (iii), and
(iv), leaving questions (v), (vi), and (vii) to Part II. We wish
to emphasize the fact that these questions are only formally
analogue to corresponding issues in Floer homology. Indeed, since
in our case the gradient-like vector field
determines a $C^1$ local flow, some of the problems above can be
dealt by dynamical systems techniques. On the other hand, finite
dimensionality and compactness results do not come from elliptic
estimates, but involve different ideas. In particular, the study
of some infinite dimensional Grassmannians, of ordinary
differential operators on Hilbert spaces, and the use of
Hausdorff measures of non-compactness  turn out to be important tools.

We conclude this introduction by giving an informal description
of our results.

\paragraph{Finite dimensional intersections.} Let $f$ be a $C^2$ Morse
function on a paracompact Hilbert manifold $M$. Let $F$ be a $C^1$
Morse vector field on $M$, having $f$ as a non-degenerate Lyapunov
function: this means that $Df(p)[F(p)]<0$ for every $p\in M$ which is
not a rest point of $F$, that the Jacobian of $F$ at every rest point
$x$ - denoted by $\nabla F(x)$ - is a hyperbolic operator, and that
the quadratic form $D^2f(x)$ is coercive on $V^-(\nabla F(x))$,
the negative eigenspace of $\nabla F(x)$, while $-D^2 f(x)$ is
coercive on the positive eigenspace $V^+(\nabla F(x))$. Under these
assumptions, $x$ is a rest point of $F$ if and only if it is a 
critical point of $f$. Typically, $F=
- \grad f$, the negative gradient of $f$ with respect to some
Riemannian metric on $M$, or $F = - h\, \grad f$, for some positive
function $h$.

The unstable and stable manifolds of a critical point $x$ are $C^1$ 
submanifolds of dimension the Morse index and co-index of $x$. 
When the critical
points $x$ and $y$ have infinite index and co-index, respectively,
the intersection $W^u(x)\cap W^s(y)$ can be infinite dimensional:
consider for example the restriction of a continuous linear form
$f$ on a Hilbert space $H$ to the unit sphere $S$ of $H$. Its
critical points are a maximum point $x$ and a minimum point $-x$,
and $W^u(x) \cap W^s(-x) = S\setminus \{x,-x\}$.

What is more striking, if $x$ and $y$ are critical points of $f$ with 
infinite Morse index and co-index, the dimension of the intersection 
between their unstable and stable manifolds (with respect to the
negative gradient flow of $f$) depends on the metric on $M$: indeed, if all the
critical points of a Morse function $f$ have infinite Morse index
and co-index, and $a:\crit (f) \rightarrow \Z$ is any function, then
$M$ supports a metric $g$ - uniformly equivalent to a given one - such that
the corresponding negative gradient flow of $f$ has the property that
for every pair of critical points $x,y$ the intersection $W^u(x) \cap
W^s(y)$ is transverse and has dimension $a(x)-a(y)$ (see \cite{ama04a}).

Therefore, some extra structure on the manifold $M$ is needed: we
will assume the existence of a subbundle $\mathcal{V}$ of $TM$,
which can be used to make comparisons. More precisely, the object
of our study will be a quartet $(M,F,f,\mathcal{V})$, where $f$ is a
non-degenerate Lyapunov function for the Morse vector field $F$, and
the subbundle $\mathcal{V}$ of $TM$ is compatible to $F$,
meaning that:
\begin{description}
\item[(C1)] for every rest point $x$, $V^+(\nabla F(x))$,
  the positive eigenspace of
  the Jacobian of $F$ at $x$, is a compact perturbation of
  $\mathcal{V}(x)$ (this means that the corresponding orthogonal
  projectors have compact difference);
\item[(C2)] denoting by $\mathcal{P}$ a projector onto
$\mathcal{V}$,
  $(L_F \mathcal{P})(p) \mathcal{P}(p)$ is a compact 
  linear operator on
  $T_p M$, for every $p\in M$ (here $L_F \mathcal{P}$ denotes the Lie
  derivative of the tensor $\mathcal{P}$ along the vector field $F$).
\end{description}
Assumption (C1) allows us to define the {\em relative Morse index
of a rest point $x$ with respect to $\mathcal{V}$} to be the integer
\[
m(x,\mathcal{V}) := \dim (V^+(\nabla F(x)),\mathcal{V}(x)) = \dim
V^+(\nabla f(x))\cap \mathcal{V}(x)^{\perp} - \dim V^-(\nabla
F(x))\cap \mathcal{V}(x).
\]
Notice that $m(x,\mathcal{V})$ can be negative. A subbundle
$\mathcal{V}=\mathcal{P}(TM)$ is invariant for the differential of
the integral flow of a vector field $X$ if and only if $(L_X
\mathcal{P})\mathcal{P}=0$. Assumption (C2) says that
$\mathcal{V}$ is {\em essentially invariant} for the linearized
flow of $F$, meaning that the differential of the flow of $F$ maps
$\mathcal{V}$ into a compact perturbation of $\mathcal{V}$. 
Assumptions (C1) and (C2) are automatically
fulfilled when all the critical points have finite index, by
choosing $\mathcal{V}=(0)$: in this case $m(x,(0))$ is the usual
Morse index.

Our first result will be that if (C1) and (C2) hold, and
$W^u(x)$, $W^s(y)$ meet transversally, then their intersection is
finite dimensional, and
\begin{equation}
\label{ladim} \dim W^u(x) \cap W^s(y) = m(x,\mathcal{V}) -
m(y,\mathcal{V}),
\end{equation}
which is the first step for the construction of the Morse complex.
A useful tool, in the proof of this result and in transversality
questions, will be the study, presented in \cite{ama03}, of the
Fredholm properties of the differential operator
\[
\frac{d}{dt} - A(t) : C^1_0(\R,H) \rightarrow C^0_0(\R,H),
\]
where the subscript $0$ means vanishing at infinity, and $A$ is a
continuous path of bounded operators on the Hilbert space $H$,
converging to hyperbolic operators for $t\rightarrow \pm \infty$.

As we shall see, the usefulness of conditions (C1) and (C2) lies
in the fact that they are both {\em stable} and {\em convex}.

In many cases, the choice of the subbundle $\mathcal{V}$ for which
(C1) and (C2) hold, is suggested by the problem itself: for
example, this is the case of semi-linear equations, where $f$ is
a lower order perturbation of a non-degenerate quadratic form on
a Hilbert space, and of many functionals coming from geometric
problems, such as the energy of curves on a
semi-Riemannian manifold. In other cases, (C1) and (C2) just hold
locally: one finds an open covering $\set{U_j}{j\in J}$ of $M$
and subbundles $\mathcal{V}_j$ of $TU_j$, which satisfy (C1),
(C2), and are such that $\mathcal{V}_i|_{U_i\cap U_j}$ is a
compact perturbation of $\mathcal{V}_j|_{U_i\cap U_j}$, for any
$i,j\in J$. That is, (C1) and (C2) hold with respect to an {\em
essential
  subbundle}. In such a situation the intersection of the unstable and
stable manifolds are finite dimensional, but no formula like
(\ref{ladim}) can possibly hold. Indeed, we will show an example
of a Morse function on $S^1 \times H$, $H$ an infinite
dimensional Hilbert space, with two rest points $x,y$, such
that different components of the transverse intersection 
$W^u(x)\cap W^s(y)$ have different
dimension. This is a purely infinite dimensional phenomenon,
related to the fact that the general linear group of an infinite
dimensional Hilbert space is contractible (see \cite{kui65}).
Formula (\ref{ladim}) will hold in the intermediate situation in
which $\dim (\mathcal{V}_i,\mathcal{V}_j)=0$ for every $i,j\in
J$: in this case we will say that (C1) and (C2) hold with respect
to a {\em (0)-essential subbundle}.

These facts are closely related to Cohen, Jones, and Segal's use
of polarizations to understand the the homotopy theory which lies
behind Floer homology \cite{cjs95}.

\paragraph{Coherent orientations.}
As we have seen, when $M$ is finite dimensional - or more
generally when the rest points have finite Morse index -
$W^u(x) \cap W^s(y)$ is orientable. In the case of infinite Morse
indices and co-indices, however, $W^u(x) \cap W^s(y)$ needs not be
orientable: indeed we will provide an example showing that such a 
transverse intersection can be diffeomorphic to $Z\times \R$, where $Z$ 
is any manifold.

The existence of a subbundle $\mathcal{V}$ satisfying (C1) and
(C2) will imply that all the intersections $W^u(x) \cap W^s(y)$
are orientable, and it will allow us to define their orientations
in a coherent way. The starting point is the fact that Fredholm
pairs (i.e.\ pairs $(V,W)$ of closed linear subspaces of a Hilbert
space $H$ with $\dim V\cap W<\infty$, $\codim (V+W) <\infty$) can
be oriented: an orientation of $(V,W)$ is by definition an
orientation of the finite dimensional space $(V\cap W) \times
(H/(V+W))^*$. Actually, a determinant bundle can be defined on the
space of Fredholm pairs, extending the determinant bundle on the
space of Fredholm operators, defined by Quillen \cite{qui85}.
Together with the fact that the fundamental group of the space of
Fredholm pairs $(V,W)$ with $\dim V=\dim W=\infty$, is $\Z_2$,
this implies that the notion of orientation of a Fredholm pair
shares all the good properties of orientations of finite
dimensional spaces.

For every rest point $x$, one fixes an orientation of the
Fredholm pair $(T_x W^s(x),\mathcal{V}(x))$. Assumptions (C1) and
(C2) guarantee that $(T_p W^s(x),\mathcal{V}(p))$ is a Fredholm
pair, for every $p\in W^s(x)$. Hence, the orientation chosen at
$x$ propagates to all the stable manifold of $x$. The way of
orienting $W^u(x)\cap W^s(y)$ is then similar to what we have
described in the case of a finite dimensional $M$.

If conditions (C1) and (C2) hold with respect to a (0)-essential
subbundle, coherent orientations cannot be defined, and one
obtains just a Morse complex with $\Z_2$ coefficients. Bott
periodicity theorem \cite{bot59} can be used to find the
obstructions to have a Morse complex with integer coefficients:
they are given by the homotopy groups $\pi_i(M)$, with $i\equiv
1,2,3,5 \mod 8$.

\paragraph{Relative compactness of the intersections.} 
When the rest point $x$ has a finite Morse index, the
(PS) condition\footnote{In this contest, a (PS) sequence is a sequence
  $(p_n)\subset M$ such that $(f(p_n))$ is bounded and
  $(Df(p_n)[F(p_n)])$ is infinitesimal.}
and the completeness of the flow imply that the intersection $W^u(x)\cap
W^s(y)$ has compact closure in $M$. When the indices are infinite,
even if (C1-2) guarantee that $W^u(x)\cap W^s(y)$ is finite
dimensional, we cannot conclude that its closure is compact: for
instance, it may consist of infinitely many isolated curves, with
no cluster points besides $x$ and $y$.

The reason is that (C1-2) are local assumptions, while compactness
involves a global condition: we shall need a global version of
condition (C2). Let us assume for simplicity that the subbundle $\mathcal{V}$ 
of $TM$ has a {\em global presentation}, that is a
submersion $\mathcal{Q}:M\rightarrow N$ into a complete Riemannian 
Hilbert manifold $N$ such that $\mathcal{V}(p) = \ker
D\mathcal{Q}(p)$. We will denote by $\beta_X(A)$ the Hausdorff measure of
non-compactness of the subset $A$ of a metric space $X$, that is the 
infimum of all positive numbers $r$
such that $A$ can be covered by finitely many balls of radius $r$.
The new assumption is:
\begin{description}
\item[(C3)] (i) $D\mathcal{Q} \circ F$ is bounded;

(ii) for every $q\in N$ there exist $\delta>0$ and $c\geq 0$
  such that $\beta_{TN}(D\mathcal{Q}(F(A)))\leq c \, \beta_N(
  \mathcal{Q}(A))$, for any
  $A$ in a $\mathcal{Q}^{-1}(B_{\delta}(q))$.
\end{description}
This condition implies (C2) by differentiation.
Condition (C3) is also stable and convex, in a sense to be
specified.

We shall prove that conditions (C1) and (C3), together with (PS)
and the completeness of the flow, imply that $W^u(x)\cap W^s(y)$ has
compact closure in $M$, for every pair of critical points $x,y$.

This compactness result will be proved in the more general setting of
a flow which preserves an {\em essentially vertical family}
$\mathcal{F}$ of subsets of $M$, with respect to a {\em strong
  integrable structure} for an essential subbundle $\mathcal{E}$ of
$TM$. When $\mathcal{E}$ is the essential class of a subbundle
$\mathcal{V}$ with a global presentation $\mathcal{Q}$, one chooses 
$\mathcal{F}$ to be the family
of subsets $A\subset M$ such that $\mathcal{Q}(A)$ is
pre-compact. More-generally, one can deal with a suitable presentation
of $\mathcal{E}$ consisting of an open covering $\{M_i\}_{i\in I}$ of
$M$ and of semi-Fredholm maps with non-negative index
$\mathcal{Q}_i:M_i \rightarrow N_i$, such that $\mathcal{E}(p) = [\ker
D\mathcal{Q}(p) ]$ for every $p\in M_i$. 

\paragraph{The boundary homomorphism.} Assume that
$(M,F,f,\mathcal{V})$ satisfies (C1-3) and (PS), 
and that the stable and unstable
manifolds of rest points meet transversally. For $q\in \Z$,
we can define $C_q(F)$ to be the free Abelian group generated by
the rest points $x$ with $m(x,\mathcal{V})=q$. In order to
define the homomorphism $\partial_q : C_q(F) \rightarrow
C_{q-1}(F)$, we just need the last condition
\begin{description}
\item[(C4)] for any $q\in \Z$, $f$ is bounded below on the set of
critical points $x$ of relative Morse index
$m(x,\mathcal{V})=q$,
\end{description}
which guarantees that the sum appearing in (\ref{fdb}) is finite.

The boundary property $\partial_{q-1} \circ \partial_q=0$ comes
from the possibility of describing exactly the flow on
the closure of each component of $W^u(x)\cap W^s(y)$, when
$m(x,\mathcal{V}) - m(y,\mathcal{V})=2$: such a flow is either
topologically conjugated to the exponential flow $(t,z)\mapsto e^t
z$ on the Riemann sphere $\C\cup \{\infty\}$, or it is
topologically conjugated to the shift flow $(t,(u,v))\mapsto
(u+t,v+t)$ on $[-\infty,+\infty] \times [-\infty,+\infty]$. In the
latter case, the orientation of this component is the product
orientation of its sides. These results will be proved by
hyperbolic dynamical systems techniques, which in this case seem
more natural than the gluing method used in Floer homology.

The resulting complex $\{C_*(F),\partial_*\}$ is said the {\em
Morse complex} of $F$. If $F_1$ and $F_2$ are two Morse vector fields
having the same non-degenerate Lyapunov function $f$, the Morse
complexes of $F_1$ and of $F_2$ are isomorphic.
In particular, their homology depends only on the Lyapunov
function $f$, and it will be said the {\em Morse homology} of $f$
and denoted by $H_*(f)$.

\paragraph{Transversality.} The transversality of the intersection of
stable and unstable manifolds will be achieved by perturbing the
vector field $F$. Small perturbations in a suitable class of vector fields
keep the conditions (C1-4) and (PS) valid: in this sense, these 
conditions were said to be stable. 
When one restricts the attention to the class of gradient vector
fields, transversality can be achieved by using rank 2 perturbations
of the given Riemannian metric.
A difference with respect to the finite dimensional case is the
regularity requirement. Indeed, high regularity of $F$ is needed
to apply Sard-Smale theorem, and such a regularity cannot be
obtained by smoothing the vector field $F$, because $C^{k+1}$
functions on an infinite dimensional Hilbert space are not $C^k$
dense (see \cite{ns73,ll86}). As a consequence, we shall assume
$F\in C^2(M)$, and we will achieve transverse intersections of
$W^u(x)$ and $W^s(y)$ whenever $m(x,\mathcal{V}) -
m(y,\mathcal{V})\leq 2$, which is what we need for the
construction of the Morse complex.

\paragraph{Relationship with classical infinite dimensional Morse theory.}
In the case of $f$ bounded below, satisfying (PS), and with
critical points of finite Morse index, we shall prove that the
Morse homology of $f$ is isomorphic to the singular
homology of $M$, a result which agrees with the Morse relations
proved by Palais. This will be a simple generalization of the
cellular filtration argument described for the compact case.

From this fact, it is easy to determine the Morse complex of some
classes of vector fields having rest points of infinite Morse
index and co-index. For instance, if $M=M^- \times M^+$ is the
product of two infinite dimensional Hilbert manifolds, endowed
with a complete product Riemannian structure, and the Morse
function $f:M \rightarrow \R$ has the special form
\begin{equation}
\label{spfrm} f(p^-,p^+) = f^+(p^+) - f^-(p^-),
\end{equation}
where $f^+: M^+ \rightarrow \R$, $f^-: M^- \rightarrow \R$ 
are bounded below and satisfy (PS), then 
\[
F = -\left( \frac{\grad f^-}{1+\|\grad f^-\|^2} , 
\frac{\grad f^+}{1+\|\grad f^+\|^2} \right) 
\]
satisfies (C1-3) with respect to the subbundle 
$\mathcal{V} = TM^- \times (0)$, with global presentation
the submersion $\mathcal{Q}:M \rightarrow M^+$, $(p^-,p^+) \mapsto p^+$. 
Notice that $F$ has the form $F(p^-,p^+) = (F^-(p^-),F^+(p^+))$.
It is easy to see that the Morse complex of $F$ is
\[
C_q(F) = (C_*(F^+) \otimes C_{-*}(F^-))_q =
\bigoplus_{\substack{(q^-,q^+) \in \N^2 \\ q^+ - q^- = q}}
C_{q^+}(F^+) \otimes C_{q^-}(F^-), \quad \forall q\in \Z,
\]
and the Morse homology of $f$ is
\begin{equation}
\label{spmh}
H_q(f) \cong (H_*(M^+) \otimes H_{-*}(M^-))_q =
\bigoplus_{\substack{(q^-,q^+) \in \N^2 \\ q^+ - q^- = q}}
H_{q^+}(M^+) \otimes H_{q^-}(M^-), \quad \forall q\in \Z.
\end{equation}

\paragraph{Computation of the homology and functoriality.}
In the case of infinite Morse indices and co-indices, the topology
of $M$ is not immediately related to the Morse homology of $f$. 
However, the homology groups $H_q(f)$ are still
considerably stable.

The key ingredient to compute the Morse homology groups will be
the fact that Morse homology is a functor from the class of Morse
functions with a gradient-like vector field
satisfying (C1-4) and (PS), seen as a small category
with the usual order relation, to the category of Abelian groups:
to every inequality $f_0\geq f_1$ is associated a homomorphism
\[
\phi_{f_0 f_1} : H_*(f_0) \rightarrow H_*(f_1),
\]
in such a way that $\phi_{f_1 f_2} \circ \phi_{f_0 f_1} =
\phi_{f_0 f_2}$, and $\phi_{f f} = \mathrm{id}$ (actually,
$\phi_{\theta\circ f f}=\mathrm{id}$, for $\theta(s)\geq s$ a
strictly increasing smooth function). The idea for the definition
of $\phi_{f_0 f_1}$ comes from the following observation: every
chain homomorphism $\psi:\{C^0_*,\partial^0_*\}\rightarrow
\{C^1_*,\partial^1_*\}$ comes from a boundary operator $\partial_q
: C^0_q \oplus C^1_{q+1} \rightarrow C^0_{q-1} \oplus C^1_q$, the
cone of $\psi$, namely
\begin{equation}
\label{cone}
\partial_q = \left( \begin{array}{cc} \partial_q^0 & 0 \\ \psi & -
\partial_{q+1}^1 \end{array} \right).
\end{equation}
With this in mind, we will construct a Morse function $f:\R
\times M \rightarrow \R$, of the form
\[
f(s,p) = \chi(s) f_0(p) + (1-\chi(s)) f_1(p) + \varphi(s),
\]
with $\chi$ a monotone smooth function such that $\chi(s)=1$ for
$s\leq 0$, and $\chi(s)=0$ for $s\geq 1$, while $\varphi(s)= 2s^3
-3s^2 +1$ has a non-degenerate maximum at $0$ and a non-degenerate
minimum at $1$. The function $f$ is a non-degenerate Lyapunov function
for a Morse vector field on $\R \times M$ 
satisfying (C1-4) and (PS), and the boundary operator in the associated Morse 
complex has the form
(\ref{cone}). This allows us to define $\phi_{f_0f_1}$ as the
homomorphism induced by the chain homomorphism $\psi$.

We wish to emphasize that this functorial approach is possible
thanks to the fact that the conditions (C1-4), and (PS) naturally
pass from the functions $f_0$, $f_1$ to the {\em cone function}
$f$: in this sense, these conditions were said to be convex.

In particular, two functions $f_0$ and $f_1$ such that $c:=
\|f_1-f_0\|_{\infty}$ is finite, have always isomorphic Morse
homologies, as implied by the functoriality applied to the
inequalities
\[
f_0-c \leq f_1 \leq f_0 + c, \quad f_1 - c \leq f_0 \leq f_1 + c.
\]
For example, let $f: M^- \times M^+ \rightarrow \R$ be a Morse function
satisfying (PS) and such that $F=-\grad
f/(1+\|\grad f\|^2)$ satisfies (C1-4) with respect to the
subbundle $\mathcal{V}= TM^- \times (0)$. If $f$ has bounded
distance from a function of the form (\ref{spfrm}), still
satisfying the same assumptions, the Morse homology of $f$
is given by (\ref{spmh}). More generally, if there exists $c> 0$ such that
\[
\frac{1}{c} f^+(p^+) - c f^-(p^-) - c \leq f(p^-,p^+) \leq c
f^+(p^+) - \frac{1}{c} f(p^-) + c,
\]
we deduce the existence of a surjective homomorphism
\[
H_q(f) \rightarrow \bigoplus_{\substack{(q^-,q^+) \in \N \\ q^+ -
q^- = q}} H_{q^+}(M^+) \otimes H_{q^-}(M^-),
\]
which implies lower estimates on the number of critical points of
$f$ of a given relative Morse index.

\numberwithin{equation}{section}

\tableofcontents

\section{Essential subbundles of a Hilbert bundle}
\label{sbu}

In this section we will fix some basic facts about the
Grassmannian of a Hilbert space and some related constructions.
We refer to Appendix A for more details.

\paragraph{Hilbert Grassmannians.}
If $E$ and $F$ are Banach spaces, $\mathcal{L}(E,F)$ will denote
the space of bounded linear operators from $E$ to $F$, while
$\mathcal{L}_c(E,F)$ will denote the subspace consisting of
compact operators. In the case $F=E$, we will simply write
$\mathcal{L}(E)$ and $\mathcal{L}_c(E)$.

Let $H$ be an infinite dimensional separable real Hilbert space.
The orthogonal projection onto a closed linear subspace $V\subset
H$ will be denoted by $P_V$, while the orthogonal complement of
$V$ will be indicated by $V^{\perp}$. We will denote by $\mathrm{Gr}(H)$
the {\em Grassmannian} of $H$, that is the space of all closed
linear subspaces of $H$, endowed with the operator norm topology.
By $\mathrm{Gr}_{\infty,\infty}(H)$ we will denote the connected component
of $\mathrm{Gr}(H)$ consisting of subspaces of infinite dimension and
infinite codimension. The other connected components of $\mathrm{Gr}(H)$
are the subsets $\mathrm{Gr}_{n,\infty}(H)$, the set of linear subspaces of
$H$ of dimension $n$, and $\mathrm{Gr}_{\infty,n}(H)$, the set of linear
subspaces of $H$ of codimension $n$.

\paragraph{Compact perturbations and essential Grassmannians.}
Given $V,W\in \mathrm{Gr}(H)$, we will say that $V$ is a {\em compact
perturbation} of $W$ if $P_V-P_W$ is a compact operator.
In this case, the {\em relative dimension} of $V$ with respect to
$W$ is the integer
\[
\dim(V,W) = \dim V\cap W^{\perp} - \dim V^{\perp}\cap W.
\]
Given $m\in \N$, the {\em $(m)$-essential Grassmannian}
$\mathrm{Gr}_{(m)}(H)$ is the quotient space of $\mathrm{Gr}(H)$ by the equivalence
relation
\[
\set{ (V,W)\in \mathrm{Gr}(H)\times \mathrm{Gr}(H) }{\mbox{$V$ is a compact
perturbation of $W$, and $\dim(V,W)\in m\Z$}}.
\]
By $\mathrm{Gr}^*_{(m)}(H)$ we will denote the quotient of
$\mathrm{Gr}_{\infty,\infty}(H)$ by the same equivalence relation. The
space $\mathrm{Gr}_{(1)}(H)$ is called just the {\em essential
Grassmannian} of $H$. If $[W]\in \mathrm{Gr}_{(m)}(H)$ and $V\in \mathrm{Gr}(H)$ is
a compact perturbation of an element (hence every element) of the
class $[W]$, then $\dim(V,[W]):= \dim(V,W)$ is well-defined as an 
integer modulo $m$.

\paragraph{Essential subbundles.}
Fix some $k\in \N\cup \{\infty\}$. Let $B$ be a topological space
if $k=0$, or a $C^k$ Banach manifold if $k\geq 1$, and let
$\mathcal{H} \rightarrow B$ be an $H$-bundle on $B$, that is a
$C^k$ fiber bundle with base space $B$, total space $\mathcal{H}$,
typical fiber the Hilbert space $H$, and structure group $\mathrm{GL}(H)$.
Since the Hilbert space $H$ is infinite dimensional, the group
$\mathrm{GL}(H)$ is contractible (see \cite{kui65}), so the above bundle is
always trivial.

We can associate to the $C^k$ Hilbert bundle $\mathcal{H}
\rightarrow B$ the $C^k$ fiber bundles
\[
\mathrm{Gr}(\mathcal{H}) = \bigcup_{b\in B} \mathrm{Gr}(\mathcal{H}_b) \rightarrow
B, \quad \mathrm{Gr}_{(m)}(\mathcal{H}) = \bigcup_{b\in B}
\mathrm{Gr}_{(m)}(\mathcal{H}_b) \rightarrow B, \quad m\in \N.
\]
The spaces $\mathrm{Gr}(H)$ and $\mathrm{Gr}_{(m)}(H)$ admit natural analytic
structures, so the above bundles have $C^k$ structures. A $C^k$
section of $\mathrm{Gr}(\mathcal{H})\rightarrow B$ is just a $C^k$
subbundle of $\mathcal{H} \rightarrow B$. Similarly, a $C^k$
section of $\mathrm{Gr}_{(m)}(\mathcal{H}) \rightarrow B$ will be called a
$C^k$ {\em $(m)$-essential subbundle of} $\mathcal{H} \rightarrow
B$, or just a $C^k$ {\em essential subbundle} in the case $m=1$.

\paragraph{Lifting properties.}
The following questions arise naturally: when is an
$(m)$-essential subbundle, $m\in \N$, liftable to a true
subbundle? when is an $(m)$-essential subbundle, $m\geq 1$,
liftable to an $(hm)$-essential subbundle, for $h\in \N$ ? We
shall discuss these questions in the nontrivial case of subbundles
with infinite dimension and codimension.

Since the Hilbert bundle $\mathcal{H}\rightarrow B$ is trivial,
the $(m)$-essential subbundle we wish to lift can be identified
with a map
\[
\mathcal{E}:B \rightarrow \mathrm{Gr}^*_{(m)}(H),
\]
and we are looking at the lifting problems
\[
\xymatrix{ & \mathrm{Gr}_{\infty,\infty}(H) \ar[d] \\ B \ar@{-->}[ru]
\ar[r]^{\mathcal{E}} &  \mathrm{Gr}_{(m)}^*(H)} \quad\quad\quad \xymatrix{ &
\mathrm{Gr}_{(hm)}^*(H) \ar[d] \\ B \ar@{-->}[ru] \ar[r]^{\mathcal{E}} &
  \mathrm{Gr}_{(m)}^*(H)}
\]
In the first diagram, the vertical map is a fibration from a
contractible space, so the $(m)$-essential subbundle $\mathcal{E}$ is
liftable to a true subbundle if and only if the map $\mathcal{E}$ is
null-homotopic. It can be proved that $\mathrm{Gr}^*_{(0)}(H)$ is simply
connected, while the fundamental group of $\mathrm{Gr}_{(m)}(H)$ for
$m\geq 1$ is infinite cyclic. Furthermore,
$\pi_i(\mathrm{Gr}_{(m)}(H))\cong \pi_{i-1}(\mathrm{BO}(\infty))$ for $i\geq 2$,
where $\mathrm{BO}(\infty)$ denotes the classifying space of the infinite
real orthogonal group. Hence, using Bott periodicity theorem, we
deduce that $\mathcal{E}$ is null homotopic if and only the homomorphism
\[
\mathcal{E}_* : \pi_i(B) \rightarrow \pi_i(\mathrm{Gr}_{(m)}^*(H))
\]
vanishes for every $i\equiv 1,2,3,5 \mod 8$. In particular, every
$(m)$-essential subbundle is liftable to a true subbundle when
$\pi_i(B)=0$ for every $i\equiv 1,2,3,5 \mod 8$.

In the second diagram, the vertical arrow is a covering map, and
the image of the induced homomorphism
\[
\pi_1(\mathrm{Gr}_{(hm)}^*(H)) \rightarrow \pi_1(\mathrm{Gr}_{(m)}^*(H)) = \Z
\]
is the subgroup $h\Z$, so the $(m)$-essential subbundle $\mathcal{E}$ is
liftable to a $(hm)$-essential subbundle if and only if $\mathcal{E}_*
(\pi_1(B)) \subset h\Z$. In particular, every $(m)$-essential
subbundle is liftable to a $(0)$-essential subbundle when $B$ is
simply connected.

In this paper, we will be mainly interested in subbundles and
essential subbundles of the tangent bundle $TM$ of a Hilbert
manifold $M$ (that is a paracompact manifold modeled on the
Hilbert space $H$). Notice that, since $M$ is locally
contractible, any $(m)$-essential subbundle $\mathcal{E}$ is
locally liftable to a true subbundle, which will be called a {\em
local representative} of $\mathcal{E}$.

\paragraph{Integrable essential subbundles of $TM$.} An essential
subbundle $\mathcal{E}$ of $TM$ is called {\em integrable} if $M$
admits an atlas whose charts $\varphi: \dom(\varphi) \subset M
\rightarrow H$ map $\mathcal{E}$ into the essential subbundle
represented by a constant closed linear subspace $V\subset H$:
\begin{equation}
\label{integr0}
\forall p\in \dom (\varphi) \quad D\varphi(p) \mathcal{E}(p) = [V].
\end{equation}
If $\varphi$ and $\psi$ are two such charts, the transition map $\tau
= \varphi \circ \psi^{-1}: \dom (\tau) \subset H \rightarrow H$
satisfies
\[
D\tau(\xi) V \mbox{ is a compact perturbation of }V \quad \forall
\xi\in \dom (\tau).
\]
By Proposition \ref{appu}, the above fact is equivalent to 
\begin{equation}
\label{integr}
Q D\tau (\xi) (I-Q) \mbox{ is a compact operator} \quad \forall \xi
\in \dom(\tau),
\end{equation}
$Q$ being a projector with kernel $V$. An atlas $\mathcal{A}$ of $M$ 
satisfying (\ref{integr0}) and (\ref{integr}) form an {\em integrable
structure modeled on $(H,V)$} for the essential subbundle $\mathcal{E}$.  

Conversely, an atlas of $M$ whose transition maps satisfy
(\ref{integr}) defines an integrable essential subbundle of $TM$. Such
an essential subbundle is liftable to an $(m)$-essential subbundle
$m\in \N$, if and only if 
\[
\dim (D\tau(\xi)V,V) \equiv 0 \quad \mod m
\]
for every transition map $\tau$ and every $\xi \in \dom(\tau)$.

Considering integrable essential subbundles will be important starting
from section \ref{scomp}. Actually, we will be interested in the
following stronger version of integrability.

\begin{defn}
\label{sintr}
Let $V$ be a closed linear subspace of the Hilbert space $H$, and let
$Q\in \mathcal{L}(H)$ be a projector with kernel $V$.
A strong integrable structure modeled on $(H,V)$ for the essential subbundle
$\mathcal{E}$ of $TM$ is atlas $\mathcal{A}$ of $M$ such that:
\begin{enumerate}
\item for every $\varphi\in \mathcal{A}$ and every $p\in \dom
  (\varphi)$, $D\varphi(p) \mathcal{E}(p) = [V]$;
\item for every $\varphi,\psi\in \mathcal{A}$ the transition map $\tau
  = \varphi \circ \psi^{-1} : \dom (\tau) \subset H \rightarrow H$
  satisfies
\[
Q A \mbox{ is pre-compact if and only if } Q \tau(A)
\mbox{ is pre-compact,}
\]
for every bounded $A\subset \dom (\tau)$.
\end{enumerate}
\end{defn}

Since the set $QA$ is pre-compact if and only if the projection of
$A$ into the quotient space $H/V$ is pre-compact, the above definition
does not depend on the choice of the projector $Q$, but only on the
subspace $V$.

Let $\tau$ be a transition map satisfying condition (ii) of the above
definition, and let $\xi \in \dom(\tau)$. Then the restriction of 
the map $Q\tau$ to the set $\dom (\tau) \cap (\xi + V)$ is a
compact map (i.e. it maps bounded sets into pre-compact sets). 
Therefore its differential at $\xi$, namely the linear
operator $Q D\tau(\xi) |_V$
is compact, implying (\ref{integr}). Hence a strong integrable
structure is also an integrable structure. The notion of a strong
integrable structure is strictly more restrictive,
because a nonlinear map whose differential at every point is compact
need not be compact. 

\begin{rem}
\label{chv}
Assume that $W$ is a compact perturbation of the closed
linear subspace $V$. Notice that if $A\subset H$
\[
P_{W^{\perp}} A \subset P_{V^{\perp}} A + (P_V - P_W) A.
\]
Then if $A$ is bounded $P_{V^{\perp}} A$ is pre-compact if and only if
$P_{W^{\perp}} A$ is pre-compact. Therefore a strong integrable
structure modeled on $(H,V)$ is also a strong integrable structure
modeled on $(H,W)$. 
\end{rem}

\paragraph{Presentations of an essential subbundle.}
A natural way to construct an integrable subbundle of $TM$ is to consider the
kernel of a submersion, or of a family of submersions with matching
kernels. We wish to describe the essential version of this construction.

The following lemma can be considered the essential version of the fact
that in suitable charts a submersion is a linear projection.

\begin{lem}\label{lin}
Let $M$ and $N$ be manifolds modeled on the Hilbert spaces $H$ and
$E$, respectively. Let $\mathcal{Q}:M\rightarrow N$ be a
$C^k$, $k\geq 1$, semi-Fredholm map with constant non-negative index. 
Then there exists a projector $Q\in \mathcal{L}(H)$ such that,
denoting by $V$ its kernel, there holds: for every $p\in M$ there
exists a $C^k$ chart $\varphi:U\rightarrow H$, $p\in U\subset M$, such that
\begin{enumerate}
\item $\varphi(U)$ is bounded;
\item for every $\xi\in \varphi(U)$, $\ker D(\mathcal{Q}\circ
  \varphi^{-1})(\xi)$ is a compact perturbation of $V$, with
\[
\dim(\ker D(\mathcal{Q}\circ \varphi^{-1}) (\xi), V)=
\dim \coker D\mathcal{Q}(\varphi^{-1}(\xi));
\]
\item for every $A\subset \varphi(U)$, $\mathcal{Q}(\varphi^{-1}(A))$
is pre-compact if and only if $QA$ is compact.
\end{enumerate}
\end{lem}

In most applications, the index of $\mathcal{Q}$ will be $+\infty$.

\proof The matter being local,
we may assume that $M$ is an open subset of the Hilbert
space $H$, that $p=0$, that $N$ is an open subset of the Hilbert
space $E$, and that $\mathcal{Q}(0)=0$. Since $\ind D\mathcal{Q}(0)
\geq 0$, there is $T\in \mathcal{L}(H,E)$ with finite rank such that 
$D\mathcal{Q}(0)+T$ is surjective. By the open mapping theorem,
$D\mathcal{Q}(0)+T$ has a linear bounded right inverse $R\in
\mathcal{L}(E,H)$. Let $Q:=R(D\mathcal{Q}(0)+T)\in \mathcal{L}(H)$ be
the associated linear projection, and set $V:=\ker Q$. Since the index
of $\mathcal{Q}$ is constant (i.e.\ it does not take different values
on different connected components of $M$), by applying a linear 
conjugacy the same $Q$ and  $V$ can be used for every point $p\in M$. 

The map $R(\mathcal{Q}+T):M \rightarrow \ran Q$ is a
local submersion at $0$, with differential at $0$ equal to $Q$.  
Therefore, there exists a
neighborhood $U\subset M$ of $0$ and a $C^k$ diffeomorphism
$\varphi:U\rightarrow H$ such that $\varphi(0)=0$,
$D\varphi(0)=I$, $\varphi(U)$ and $\mathcal{Q}(U)$ bounded 
(so that (i) holds), and
\begin{equation}
\label{ift}
R(\mathcal{Q} + T)\circ \varphi^{-1}(\xi) = Q\xi, 
\quad \forall \xi \in \varphi(U).
\end{equation}
Differentiating we get $R D(Q\circ \varphi^{-1})(\xi) + RT
D\varphi^{-1}(\xi)  = Q$. Since $R$ is injective and since $T$ has
finite rank, Proposition \ref{trasf} implies that 
\[
\ker D(\mathcal{Q}\circ \varphi^{-1}) (\xi) = \ker R
D(\mathcal{Q}\circ \varphi^{-1}) (\xi) 
\]
is a compact perturbation of $\ker Q=V$, and
\begin{eqnarray*}
\dim (\ker D(\mathcal{Q}\circ \varphi^{-1}) (\xi),V) = - \dim (\ran
RD(\mathcal{Q}\circ \varphi^{-1}) (\xi) , \ran Q) \\
= - \dim ( \ran D(\mathcal{Q}\circ \varphi^{-1}) (\xi), E) =
\dim \coker D\mathcal{Q}(\varphi^{-1}(\xi)) D\varphi^{-1}(\xi) = 
\dim \coker D\mathcal{Q}(\varphi^{-1}(\xi)),
\end{eqnarray*}
proving (ii). By (\ref{ift}), for every $A\subset \varphi(U)$,
\[
QA \subset R\mathcal{Q}(\varphi^{-1}(A))+ \ran RT, \quad
\mathcal{Q}(\varphi^{-1}(A)) \subset (D\mathcal{Q}(0)+T) QA + \ran T,
\]
so claim (iii) follows from the fact that $\varphi(U)$ and
$\mathcal{Q}(U)$ are bounded, and from the fact that $T$ has finite rank.
\qed

\begin{prop} \label{pesfe}
Consider an open covering $\{M_i\}_{i\in I}$ of $M$, a family of
infinite dimensional Hilbert manifolds $\{N_i\}_{i\in I}$ modeled on $E$, and a
family of semi-Fredholm $C^k$, $k\geq 1$, maps 
$\mathcal{Q}_i:M_i \rightarrow N_i$ with the same constant
non-negative index, such that 
for any $i,j\in I$ and for any $A\subset M_i\cap M_j$, 
\begin{equation}
\label{lni}
\mathcal{Q}_i(A) \mbox{ is
pre-compact if and only if } \mathcal{Q}_j(A) \mbox{ is pre-compact.}
\end{equation}
Then the family $\set{\ker D \mathcal{Q}_i(p)}{p\in M_i}$, $i\in I$, defines
a $C^{k-1}$ essential subbundle $\mathcal{E}$ of $TM$. The atlas
$\mathcal{A}$ consisting of all the charts
$\varphi$ of $M$ satisfying properties (i), (ii), and (iii) of Lemma
\ref{lin} applied to all the maps $\mathcal{Q}_i$ is a strong
integrable structure modeled on $(H,V)$ for $\mathcal{E}$. This atlas
is such that for every $\varphi\in \mathcal{A}$ and every $A\subset
\dom(\varphi) \subset M_i$, 
\begin{equation}
\label{pip} Q\varphi(A) \mbox{ is
pre-compact if and only if } \mathcal{Q}_i(A) \mbox{ is pre-compact.}  
\end{equation}
Moreover, for every $p\in M_i\cap M_j$ the operator $D
\mathcal{Q}_i(p)\, D \mathcal{Q}_j(p)^*\in
\mathcal{L}(T_{\mathcal{Q}_j(p)} N_j,T_{\mathcal{Q}_i(p)} N_i)$ is 
Fredholm, and $\mathcal{E}$ is liftable to an $(m)$-essential 
subbundle, $m\in \N$, if and only if
\[
\ind (D \mathcal{Q}_i(p)\,D \mathcal{Q}_j(p)^*) \equiv 0 \mod m, 
\quad \forall i,j\in I, \; \forall p\in M_i\cap M_j.
\]
\end{prop}

\proof Let us prove that for any $p\in M_i\cap M_j$ the subspace
$\ker D\mathcal{Q}_i(p)$ is a compact perturbation of 
$\ker D\mathcal{Q}_j(p)$. Since $D\mathcal{Q}_i(p)$ has finite corank, 
we can find a $C^1$ embedded finite dimensional open disk
$D\subset N_i$ with $\overline{D}$ compact, such that 
$\mathcal{Q}_i(p)\in D$ and the map $\mathcal{Q}_i$ is transverse to
$D$. Then $\mathcal{Q}_i^{-1}(D)$ is a $C^1$ submanifold of $M$, and
by our assumption the map $\mathcal{Q}_j|_{\mathcal{Q}_i^{-1}(D)\cap
M_j}$ is compact. Therefore its differential, namely the restriction of 
$D\mathcal{Q}_j$ to the subspace $T_p \mathcal{Q}^{-1}_i(D)\supset 
\ker D\mathcal{Q}_i(p)$ is compact. In particular, the restriction of 
$D\mathcal{Q}_j(p)$ to $\ker D\mathcal{Q}_i(p)$ is compact, and
similarly the restriction of $D\mathcal{Q}_i(p)$ to 
$\ker D\mathcal{Q}_j(p)$ is 
compact. Hence Proposition \ref{kersf} implies
that $\ker D\mathcal{Q}_i(p)$ is a compact perturbation of $\ker
D\mathcal{Q}_j(p)$, as we wished to prove, and that
\begin{equation}
\label{kuno} \ind (D\mathcal{Q}_i(p)\, D\mathcal{Q}_j(p)^*) = 
\dim \coker D\mathcal{Q}_j(p) - \dim \coker D\mathcal{Q}_i(p) + \dim 
(\ker D\mathcal{Q}_i(p),\ker D\mathcal{Q}_j(p)).
\end{equation}
Now let $\varphi$ and $\psi$, $\dom(\varphi)\subset M_i$, $\dom
(\psi)\subset M_j$, be two charts satisfying conditions (i), (ii),
and (iii) of Lemma \ref{lin} applied to $\mathcal{Q}_i$ and
$\mathcal{Q}_j$, respectively (possibly $i=j$). Let $\tau=\varphi\circ
\psi^{-1}$ be the transition map. If $A\subset \dom(\tau) = \psi(\dom
(\varphi) \cap \dom (\psi))$, by Lemma \ref{lin}.(iii) $QA$ is
pre-compact if and only if $\mathcal{Q}_j(\psi^{-1}(A))$ is
pre-compact, by (\ref{lni}) if and only if $\mathcal{Q}_i(\psi^{-1}(A)) =
\mathcal{Q}_i(\varphi^{-1}(\tau(A)))$ is pre-compact, and again by
Lemma \ref{lin}.(iii) if and only if $Q\tau(A)$ is pre-compact. This
proves condition (ii) of Definition \ref{sintr}, and proves that the
atlas $\mathcal{A}$ satisfies (\ref{pip}).

Finally, let $p\in M_i$. By Lemma \ref{lin}.(ii), there is a
neighborhood $U_p$ of $p$ and a $C^k$ submersion $\widetilde{\mathcal{Q}}_p:=
Q\varphi:U_p \rightarrow \ran Q$ into a Hilbert space such that for
any $q\in U_p$, $\ker D\mathcal{Q}_i(q)$ is a compact perturbation of
$\ker D\widetilde{\mathcal{Q}}_p(q)$, and 
\begin{equation}
\label{kdue} \dim (\ker D \mathcal{Q}_i(q), \ker 
D\widetilde{\mathcal{Q}}_p (q)) = \dim \coker D\mathcal{Q}_i(q).
\end{equation}
Then the family $\set{\ker D\mathcal{Q}_i}{i\in I}$ defines the same
$C^{k-1}$ essential subbundle of $TM$ as the one defined by the
family
\begin{equation}
\label{ktre} \set{\ker D\widetilde{\mathcal{Q}}_p}{p\in M}.
\end{equation}
If $q\in U_p \cap U_{p^{\prime}}\cap M_i \cap M_j$, by
(\ref{kuno}) and (\ref{kdue}) we obtain (see formula
(\ref{astar}))
\begin{eqnarray*}
\dim (\ker D\widetilde{\mathcal{Q}}_p(q), \ker
\widetilde{\mathcal{Q}}_{p^{\prime}}(q)) =
\dim (\ker D\widetilde{\mathcal{Q}}_p(q), \ker D\mathcal{Q}_i(q)) + 
\dim (\ker D\mathcal{Q}_i(q),
\ker D\mathcal{Q}_j(q)) \\ + \dim (\ker D\mathcal{Q}_j(q), \ker
D\widetilde{\mathcal{Q}}_{p^{\prime}}(q)) = \ind (D\mathcal{Q}_i(q)\,
D\mathcal{Q}_j(q)^*),
\end{eqnarray*}
so (\ref{ktre}) defines an $(m)$-essential subbundle of $TM$ if and only
if
\[
\ind (D \mathcal{Q}_i(q)\, D\mathcal{Q}_j(q)^*) \equiv 0 \mod m, 
\quad \forall i,j\in I, \; \forall q\in M_i\cap M_j.
\]
\qed

The above proposition suggests the following:

\begin{defn}
\label{spres}
A strong presentation of the essential subbundle $\mathcal{E}$ of $TM$
consists of an open covering $\{M_i\}_{i\in I}$ of $M$, a family of 
manifolds $N_i$, $i\in I$, modeled on the Hilbert space $E$, 
a family of semi-Fredholm $C^1$ maps $\mathcal{Q}_i:M_i \rightarrow
N_i$ with the same constant non-negative index such that: 
\begin{enumerate}
\item for every $i\in I$ and every $p\in M_i$, the kernel of
  $D\mathcal{Q}_i(p)$ belongs to the essential class $\mathcal{E}(p)$;
\item for every $i,j\in I$ and every $A\subset M_i\cap M_j$, 
$\mathcal{Q}_i(A)$ is pre-compact if and only if $\mathcal{Q}_j(A)$ is
pre-compact.
\end{enumerate}
\end{defn}

Proposition \ref{pesfe} states among other facts that a strong
presentation of $\mathcal{E}$ determines a strong integrable structure
for $\mathcal{E}$.

\section{Morse vector fields and subbundles}

\paragraph{Definitions and basic facts.}
Let $M$ be a paracompact manifold of class $C^2$, modelled on the
infinite dimensional separable real Hilbert space $H$. Let $F$ be
a tangent vector field of class $C^1$ on $M$. This field
determines a local flow on $M$,
\[
\phi\in C^1(\Omega(F),M), \quad \partial_t \phi(t,p) =
F(\phi(t,p)), \quad \phi(0,p)=p,
\]
where $\Omega(F)\subset \R\times M$ is the maximal set of
existence for the solutions of this ordinary differential
equation. We will also use the notation $\phi_t(p)=\phi(t,p)$.

A {\em rest point} of $F$ is a point $x\in M$ such that $F(x)=0$.
The set of rest points of $F$ is denoted by $\rest(F)$. If $x\in
\rest(F)$, the Jacobian of $F$ at $x$, $\nabla F(x)$, is the
bounded linear operator on $T_x M$ defined as
\[
\nabla F(x) \xi = L_X F (x), \quad \mbox{for $X$ a tangent vector
field
  such that $X(x)=\xi\in T_x M$,}
\]
where $L_X F$ denotes the Lie derivative of $F$ along $X$.
Indeed, the fact that $F(x)=0$ implies that $L_X F(x)$ depends
only on the value of $X$ at $x$.

We recall that an operator $L\in \mathcal{L}(H)$ is said {\em
  hyperbolic} if $\sigma(L)\cap i\R = \emptyset$. In this case, the
  decomposition of the spectrum of $L$ into the subset with positive
  real part and the one with negative real part determines an
  $L$-invariant splitting $H=V^+(L) \oplus V^-(L)$.
A point $x\in \rest(F)$ is said {\em hyperbolic} if the operator
$\nabla F(x)$ is hyperbolic. In this case, the {\em linear
unstable
  space} $H^u_x$ and the {\em linear stable space} $H^s_x$, are defined as
\[
H^u_x := V^+(\nabla F(x)), \quad H^s_x := V^-(\nabla F(x)).
\]
A vector field $F$ all of whose rest points are hyperbolic is said
a {\em Morse vector field}.

A {\em Lyapunov function} for the vector field $F$ is a function
$f\in C^1(M)$ such that
\begin{equation}
\label{decr} Df(p) [F(p)] <0 , \quad \forall p\in M \setminus
\rest(F).
\end{equation}
In particular, $t\mapsto f(\phi(t,p))$ is strictly decreasing if
$p\notin \rest(F)$. Note that every critical point of $f$ must be
a rest point of $F$. If $x$ is a hyperbolic rest point for $F$,
then it is a critical point of $f$, as it easily follows from a
first order expansion of $F$ at $x$.

If the vector field $F$ is Morse, we shall ask the Lyapunov
function to be {\em non-degenerate}: $f$ is twice differentiable
at every rest point $x$ and, denoting by $D^2 f(x)$ the second
differential of $f$ at $x$, seen as a symmetric bounded bilinear
form, we have that $\xi \mapsto D^2 f(x) [\xi,\xi]$ is coercive on $H^s_x$,
while $\xi \mapsto -D^2 f(x) [\xi,\xi]$ is coercive on $H^u_x$. The Morse
vector field $F$ is said {\em gradient-like} if it has a
non-degenerate Lyapunov function.

\paragraph{The relative Morse index.}
For $\mathcal{V}$ a subbundle of $TM$ of class $C^1$, consider
the following compatibility condition between $F$ and
$\mathcal{V}$:
\begin{description}
\item[(C1)] for every $x$ rest point of $F$, the linear unstable
space $H^u_x$ is a compact perturbation of $\mathcal{V}(x)$.
\end{description}

If (C1) holds, the {\em relative Morse index} of $x\in \rest(F)$
is the integer
\[
m(x,\mathcal{V}) := \dim (H^u_x,\mathcal{V}(x)),
\]
and the sets
\[
\mathrm{rest}_q(F):=\set{x\in\rest(F)}{m(x,\mathcal{V})=q}, \quad q\in
\Z,
\]
constitute a partition of $\rest (F)$.

Condition (C1) clearly depends only on the essential class of
$\mathcal{V}$. Therefore it makes sense to talk about vector
fields which satisfy (C1) with respect to an essential subbundle.
More precisely, the $C^1$ Morse vector field $F$ satisfies (C1)
with respect to the essential subbundle $\mathcal{E}$ if for every
rest point $x$ of $F$ the unstable space $H^u_x$ belongs to the
essential class $\mathcal{E}(x)$. In this more general situation,
there is no relative Morse index. However, if the essential
subbundle $\mathcal{E}$ comes from an $(m)$-essential subbundle -
still denoted by $\mathcal{E}$ - then the relative Morse index of
$x\in \rest(F)$ is an integer modulo $m$ - denoted by
$m(x,\mathcal{E})$. In particular, if $\mathcal{E}$ is a
$(0)$-essential subbundle, the relative Morse index is still
integer valued.

\paragraph{Essentially invariant subbundles.}
We shall say that the $C^1$ subbundle $\mathcal{V}$ is {\em
invariant with respect to $F$ at $p\in M$} if
\begin{equation}
\label{huno} (L_F \mathcal{P})(p) \mathcal{P}(p) = 0,
\end{equation}
where $\mathcal{P}$ is a projector onto $\mathcal{V}$ of $TM$:
$\mathcal{P}$ is a $C^1$ section of the Banach bundle of linear
endomorphisms of $TM$ such that for
every $p\in M$, $\mathcal{P}(p)\in \mathcal{L}(T_p M)$ is a
projector onto $\mathcal{V}(p)$. This notion does not depend on
the choice of the projector $\mathcal{P}$, but only on the
subbundle $\mathcal{V}$. Indeed, if $\mathcal{P}$ and
$\mathcal{Q}$ are two projectors onto $\mathcal{V}$, we have the
identity
\begin{equation}
\label{hdue} (L_F \mathcal{Q}) \mathcal{Q} = (I-\mathcal{Q}) (L_F
\mathcal{P}) \mathcal{P} \mathcal{Q},
\end{equation}
which can be verified by taking the Lie derivative of the
identities $\mathcal{PQ} = \mathcal{Q} = \mathcal{Q}^2$. This
definition is motivated by the well known fact that (\ref{huno})
holds for any $p\in M$ if and only if the subbundle $\mathcal{V}$
is invariant under the action of the local flow $\phi$, that is
$D \phi_t (p) \mathcal{V}(p) = \mathcal{V}(\phi_t(p))$ for every
$(t,p)\in \Omega(F)$.

Similarly, we shall say that $\mathcal{V}$ is {\em essentially
invariant with respect to $F$ at $p$} if $(L_F \mathcal{P})(p)
\mathcal{P}(p)$ is a compact endomorphism of $T_p M$. Again,
(\ref{hdue}) shows that this notion depends only on
$\mathcal{V}$. By Proposition \ref{appu}, $\mathcal{V}$ is
essentially invariant with
respect to $F$ at every $p\in M$ if and only if $D\phi_t (p)
\mathcal{V}(p)$ is a compact perturbation of
$\mathcal{V}(\phi_t(p))$, for every $(t,p)\in \Omega(F)$. The
second compatibility condition between $F$ and $\mathcal{V}$ is:
\begin{description}
\item[(C2)] $\mathcal{V}$ is essentially invariant with respect
to $F$ at any point $p\in M$.
\end{description}
Also this condition can be stated for an essential subbundle.
Indeed, an essential subbundle $\mathcal{E}$ of $TM$ will be said
{\em
  invariant with respect to $F$ at $p$} if a local representative of
$\mathcal{E}$ at $p$ is essentially invariant with respect to $F$
at $p$. This notion does not depend on the choice of the local
representative of $\mathcal{E}$ at $p$: if $\mathcal{V}$ and
$\mathcal{W}$ are two such local representatives on some
neighborhood $U$ of $p$, and $\mathcal{P}$, $\mathcal{Q}$ are the
orthogonal projectors onto $\mathcal{V}$, $\mathcal{W}$, with
respect to some Riemannian structure on $M$, we have that
$\mathcal{P}(q) - \mathcal{Q}(q)\in \mathcal{L}_c(T_q M)$ for any
$q\in U$, so $(L_F (\mathcal{P} - \mathcal{Q}))(p) \in
\mathcal{L}_c(T_p M)$, and the identity
\[
(L_F \mathcal{P}) \mathcal{P} - (L_F \mathcal{Q}) \mathcal{Q} =
(L_F \mathcal{P}) (\mathcal{P} - \mathcal{Q}) + (L_F (\mathcal{P}
- \mathcal{Q})) \mathcal{Q},
\]
shows that $(L_F \mathcal{P})\mathcal{P}$ is compact if and only
if $(L_F \mathcal{Q})\mathcal{Q}$ is compact. Hence, we shall say
that the $C^1$ vector field $F$ satisfies (C2) with respect to the
$C^1$ essential subbundle $\mathcal{E}$ of $TM$ if $\mathcal{E}$
is invariant with respect to $F$ at every $p\in M$.

\begin{prop}
Let $\mathcal{E}$ be a $C^1$ essential subbundle of $TM$. Then
the set of $C^1$ vector fields on $M$ which satisfy (C2) with
respect to $\mathcal{E}$ is a $C^1(M)$-module.
\end{prop}

\proof Everything follows from the formulas
\[
(L_{X+Y} \mathcal{P}) \mathcal{P} = (L_X \mathcal{P})\mathcal{P}
+ (L_Y \mathcal{P})\mathcal{P}, \quad (L_{h X} \mathcal{P})
\mathcal{P} \xi = h (L_X \mathcal{P}) \mathcal{P} \xi +
Dh[\mathcal{P} \xi] (\mathcal{P} - I)X, \quad \forall
\xi\in TM,
\]
where $h\in C^1(M)$. \qed

\paragraph{Examples.} We conclude this section with some simple examples.

\begin{ex} \label{ex1}
{\em (Vector fields whose rest points have finite Morse index or
  finite Morse co-index)}
Consider the classical situation of a Morse vector field $F$ all
of whose rest points have finite Morse index. Then (C1) and (C2)
hold with respect to the trivial subbundle $\mathcal{V}=(0)$.
With such a $\mathcal{V}$ indeed, (C2) is fulfilled by any vector
field, while (C1) is equivalent to asking the unstable space of
every rest point to be finite dimensional. In this case,
$m(x,(0))$ is the usual Morse index of the rest point $x$.

Similarly, a Morse vector field all of whose rest points have
finite Morse co-index satisfies (C1) and (C2) with respect to the
trivial subbundle $\mathcal{V}=TM$, and $-m(x,TM)$ is the
co-index of the rest point $x$.
\end{ex}

\begin{ex} \label{ex2}
{\em (Perturbations of a non-degenerate quadratic form)} Assume
that $M=H$ is a Hilbert space, and consider a function of the form
\[
f(\xi) = \frac{1}{2} \langle L\xi,\xi \rangle + b(\xi),
\]
where $L\in \mathcal{L}(H)$ is self-adjoint invertible, and $b\in
C^2(H)$. Let $F$ be the (negative) gradient vector field of $f$,
\[
F(\xi) = - \grad f(\xi) = - L \xi - \grad b(\xi),
\]
and consider the constant subbundle $V=V^-(L)$. In this case,
condition (C2) means asking that
\[
(L_{\grad f} P_{V})(\xi) P_V = [P_V, \hess f(\xi)] P_V = [P_V,
\hess b(\xi)] P_V
\]
should be compact for every $\xi\in H$. In particular, if we
assume that the Hessian of $b$ at every point is compact,
condition (C2) holds.  Since the negative eigenspace of a compact
perturbation of $L$ is a compact perturbation of $V$ (Proposition
\ref{abuno}), also condition (C1) holds.
\end{ex}

\begin{ex}{\em (Product manifolds)} \label{ex3}
Assume that $M=M^- \times M^+$ is the product of two Hilbert
manifolds, and consider the subbundle $\mathcal{V}=TM^- \times
(0)$ of $TM$. Fix some Riemannian structure on $M^-$ and on $M^+$,
and consider the product Riemannian structure on $M$. Let
$F=-\grad f$ be the negative gradient of a Morse function on $M$.
Then $F$ satisfies (C1) with respect to $\mathcal{V}$ if and only
if for every critical point $x$ the Hessian of $f$ at $x$
decomposes as $\hess f(x) = L_x + K_x$, where $L_x$ is
self-adjoint, invertible, and $V^-(L_x)=\mathcal{V}(x)$, while
$K_x$ is a compact endomorphism of $T_x M$. Moreover, $F$
satisfies (C2) with respect to $\mathcal{V}$ if and only if for
every $p\in M$ the operator
\[
(L_F \mathcal{P})(p) \mathcal{P}(p) = (\nabla_{\grad f} \mathcal{P}
(p) + [\mathcal{P}(p), \hess f(p)]) \mathcal{P}(p)
\]
is compact, where $\mathcal{P}$ denotes the orthogonal projection
onto $\mathcal{V}$.
\end{ex}

\begin{ex}{\em (Semi-Riemannian geodesics \cite{ama04b})} Let
  $Q$ be an $n$-dimensional manifold, endowed with a semi-Riemannian
  structure $h$, that is a symmetric non-degenerate bilinear form on
  $TQ$. Denote by $(n^+,n^-)$ the signature of $h$, $n^+ + n^- = n$. 
  The semi-Riemannian structure $h$ induces a Levi-Civita
  covariant derivation $\nabla$, and the geodesics - i.e. the
  $1$-periodic solutions $q$ of the second order ODE $\nabla_{\dot{q}}
  \dot{q} = 0$ - joining two fixed points $q_0,q_1\in Q$ 
  are the critical points of the energy functional
\[
f(q) = \frac{1}{2} \int_0^1 h(\dot{q}(t),\dot{q}(t)) \, dt,
\]
on the Hilbert manifold $M:= \set{q\in W^{1,2}([0,1],Q)}{q(0)=q_0, 
q(1)=q_1}$ consisting of paths in
$Q$ of Sobolev class $W^{1,2}$ joining $q_0$ and $q_1$. 
When $n^+\neq 0$ and $n^-\neq 0$, all
the critical points of $f$ have infinite Morse index and
co-index. Assume that $TQ$ has an {\em integrable} subbundle $V$ of
dimension $n^-$ such that $h$ is strictly negative on $V$, and set
\[
\mathcal{V}(q) = \set{\zeta\in T_q M = q^*(TQ)}{\zeta(t) \in V(q(t)) \;
  \forall t\in [0,1]}, \quad \forall q\in M.
\]
The integrability of $V$ is reflected into the integrability of
$\mathcal{V}$, and this fact can be used to build a class of
Riemannian structures on $M$ - equivalent to the standard $W^{1,2}$
metric - such that $\grad f$ satisfies (C1) and (C2) with respect to
$\mathcal{V}$. In this situation, it can also be proved that the
relative Morse index $m(q,\mathcal{V})$ of the geodesic $q$
coincides with the Maslov index of a suitable path of Lagrangian
subspaces, obtained by looking at the Hamiltonian
system on the cotangent bundle of $Q$ generated by the Legendre
transform $H:T^*Q \rightarrow \R$ of the Lagrangian $L:TQ \rightarrow
\R$, $L(\zeta) = 1/2 \,h(\zeta,\zeta)$.
\end{ex}      

\section{Finite dimension of $\mathbf{W^u(x) \cap W^s(y)}$}
\label{sstab}

\paragraph{Stable and unstable manifolds.}
The unstable and stable manifolds of a hyperbolic rest point $x$
are the sets
\begin{eqnarray*}
W^u(x) & := & \set{p\in M}{]-\infty,0] \times \{p\} \subset
\Omega(F) \mbox{ and } \phi(t,p)\rightarrow x \mbox{ for }
t\rightarrow -\infty}, \\
W^s(x) & := & \set{p\in M}{[0,+\infty[ \times \{p\} \subset
\Omega(F) \mbox{ and } \phi(t,p)\rightarrow x \mbox{ for }
t\rightarrow +\infty},
\end{eqnarray*}
and classical results in the theory of dynamical systems imply
that $W^u(x)$ and $W^s(x)$ are the images of injective $C^1$
immersions of $H^u_x$ and $H^s_x$, respectively, and that
\[
T_x W^u(x) = H^u_x, \quad T_x W^s(x) = H^s_x.
\]
In general, they need not be embedded submanifolds. Starting from
section \ref{scomp} however, we will restrict our attention to
gradient-like vector fields, for which $W^u(x)$ and $W^s(x)$ are
embedded submanifolds (see also Appendix C). 

\begin{prop}
\label{stab} Let $\mathcal{E}$ be an essential subbundle of $TM$,
and let $x$ be a hyperbolic rest point of the $C^1$ vector field
$F$ on $M$. Then the following facts are equivalent:
\begin{enumerate}
\item $H^u_x$ belongs to the essential class $\mathcal{E}(x)$,
  and $\mathcal{E}$ is invariant with respect to $F$ at every $p\in W^u(x)$;
\item the tangent space $T_p W^u(x)$ belongs to the essential
class $\mathcal{E}(p)$ for every $p\in W^u(x)$.
\end{enumerate}
If either (i) or (ii) holds, and if $\mathcal{E}$ is liftable to
an $(m)$-essential subbundle - still denoted by $\mathcal{E}$ -
then we have the identity between integers modulo $m$
\[
\dim(T_p W^u(x),\mathcal{E}(p)) = m(x,\mathcal{E}), \quad
\forall p\in W^u(x).
\]
\end{prop}

\proof Let $p\in W^u(x)$ and define $u:[-\infty,0] \rightarrow M$
by $u(t):=\phi_t(p)$ for $t>-\infty$, and $u(-\infty) =x$. If
$\psi:U\rightarrow H$, $x\in U$, is a local chart mapping the open
set $U$ diffeomorphically into the Hilbert space $H$, then for $T$
large $\psi\circ \phi_{-T}:\phi_{-T}^{-1}(U)\rightarrow H$ is a
local chart whose domain contains $u([-\infty,0])$. Therefore,
since both the assertions of the theorem are invariant with
respect to differentiable conjugacy, we may assume that $M$ is an
open subset of $H$. The set $\phi([-\infty,0] \times \{p\})$ has a
contractible neighborhood $U$, and we can find a $C^1$ map
$\mathcal{P}:U \rightarrow \mathcal{L}(H)$ such that
$\mathcal{P}(\xi)$ is a projector onto a subspace in the essential
class $\mathcal{E}(\xi)$, for every $\xi\in U$.

Set $P:= \mathcal{P}(x)$, and let $R:[-\infty,0] \rightarrow
\mathrm{GL}(H)$ be such that $R(t)P=\mathcal{P}(u(t))R(t)$, $R(-\infty)=I$,
and $R^{\prime}(t)\rightarrow 0$ for $t\rightarrow -\infty$. Set
\[
X(t):= R(t)^{-1}D\phi_t(p)R(0).
\]
Then $X$ solves $X^{\prime}=AX$, $X(0)=I$, where
\[
A(t) = R(t)^{-1}R^{\prime}(t) + R(t)^{-1}DF(u(t)) R(t)\in
\mathcal{L}(H)
\]
converges to the hyperbolic operator $A(-\infty)=DF(x)$ for
$t\rightarrow -\infty$. Let
\[
W^u_A := \set{\xi\in H}{\lim_{t\rightarrow -\infty} X(t)\xi=0},
\]
be the linear unstable space of the path of operators 
$A$ (see Appendix B). Then
\begin{equation}
\label{er} T_{u(t)} W^u(x) = R(t) X(t) W^u_A, \quad T_{x} W^u(x)
= V^+(A(-\infty)).
\end{equation}
Differentiating $R(t)P=\mathcal{P}(u(t))R(t)$ we obtain the
identity
\[
R^{\prime}(t)P=D\mathcal{P}(u(t))[F(u(t))]R(t) +
\mathcal{P}(u(t))R^{\prime}(t),
\]
from which an easy computation gives
\[
[A(t),P]P = R(t)^{-1} \Bigl( D\mathcal{P}(u(t)) [F(u(t))] +
  [\mathcal{P}(u(t)), DF(u(t))] \Bigr) \mathcal{P}(u(t)) R(t).
\]
So by the usual expression for the Lie derivative,
\begin{equation}
\label{er3} [A(t),P]P = R(t)^{-1} (L_F \mathcal{P})(u(t))
\mathcal{P}(u(t)) R(t),
\end{equation}
and the equivalence of (i) and (ii) follows form (\ref{er}),
(\ref{er3}), and Proposition \ref{utl}.

Assume now that $\mathcal{E}$ comes from an $(m)$-essential
subbundle. Since $W^u(x)$ is connected and the relative dimension
is a continuous function, for every $p\in W^u(x)$ we have the
following identity between integers modulo $m$
\[
\dim(T_p W^u(x),\mathcal{E}(p)) = \dim(T_x W^u(x), \mathcal{E}(x)) =
\dim (H^u_x,\mathcal{E}(x)) = m(x,\mathcal{E}).
\]
\qed

Recall that a pair of closed subspaces $(V,W)$ of the Hilbert
space $H$ is said a {\em Fredholm pair} if $V\cap W$ has finite
dimension and $V+W$ has finite codimension, in which case we
define the index of $(V,W)$ to be
\[
\ind (V,W) = \dim V\cap W - \codim (V+W).
\]
The space of Fredholm pairs of $H$, denoted by $\mathrm{Fp}(H)$, is an open
subspace of $\mathrm{Gr}(H)\times \mathrm{Gr}(H)$, and the index is a continuous
function. If $\mathcal{H}\rightarrow B$ is a $C^k$ Hilbert bundle,
there is an associated $C^k$ bundle
\[
\mathrm{Fp}(\mathcal{H}) = \bigcup_{b\in B} \mathrm{Fp}(\mathcal{H}_b) \rightarrow
B.
\]
The above proposition has the following corollary.

\begin{cor}
\label{coro} Assume that the Morse vector field $F$ satisfies
(C1-2) with respect to a subbundle $\mathcal{V}$ of $TM$. Then
for every rest point $x$:
\begin{enumerate}
\item for any $p\in W^u(x)$, $T_p W^u(x)$ is a compact
perturbation of $\mathcal{V}(p)$, and $\dim(T_p W^u(x), \mathcal{V}(p)) =
  m(x,\mathcal{V})$;
\item for any $p\in W^s(x)$, $(T_p W^s(x),\mathcal{V}(p))$ is a
  Fredholm pair of index $-m(x,\mathcal{V})$.
\end{enumerate}
\end{cor}

\proof Assertion (i) follows immediately from Proposition
\ref{stab} and from the continuity of the relative dimension. By
(C1) and Proposition \ref{indice}, $(T_x W^s(x),\mathcal{V}(x))=
(H^s_x,\mathcal{V}(x))$ is a Fredholm pair of index
\[
\ind (T_x W^s(x),\mathcal{V}(x)) = \ind (H^s_x,H^u_x) + \dim
(\mathcal{V}(x),H^u_x) = - m(x,\mathcal{V}).
\]
Therefore, $(T_p W^s(x),\mathcal{V}(p))$ is a Fredholm pair of the
same index for any $p$ in a neighborhood $U$ of $x$ in the
intrinsic topology of the immersed submanifold $W^s(x)$. The
backward evolution of $U$ by $\phi$ is the whole $W^s(x)$, so
assertion (ii) follows from the fact that the tangent bundle of
$W^s(x)$ is invariant, and $\mathcal{V}$ is essentially invariant
under the action of $\phi$. 
\qed

\paragraph{Intersections.}
Recall that two immersed submanifolds $N,O\subset M$ have a {\em
transverse intersection} if for every $p\in N\cap O$ there holds
$T_p N + T_p O = T_p M$. In this case, $N\cap O$ is an immersed
submanifold of $M$, and $T_p (N\cap O) = T_p N \cap T_p O$.
Similarly, $N,O\subset M$ have a {\em Fredholm intersection} if
for every $p\in N \cap O$, $(T_p N,T_p O)$ is a Fredholm pair of
subspaces of $T_p M$. We are now ready to state the result about
the dimension of the intersection of the unstable and the stable
manifolds.

\begin{thm}
\label{ubu} Let $k\in \N$, let $\mathcal{E}$ be a
$(k)$-essential subbundle of $TM$, and assume that the Morse
vector field $F$ satisfies (C1-2) with respect to $\mathcal{E}$.
Let $x,y$ be two rest points of $F$. Then $W^u(x)$ and $W^s(y)$
have Fredholm intersection, with the number
\[
\ind (T_p W^u(x),T_p W^s(y)), \quad p\in W^u(x)\cap W^s(y),
\]
depending only on the homotopy class of the curve $t\mapsto
\phi(t,p)$ in the space of continuous paths
$u:\overline{\R}\rightarrow M$ such that $u(-\infty)=x$,
$u(+\infty)=y$.  Furthermore
\begin{equation}
\label{lequa} \ind (T_p W^u(x),T_p W^s(y)) \equiv
m(x,\mathcal{E}) - m(y,\mathcal{E}) \quad \mod k,
\end{equation}
for every $p\in W^u(x) \cap W^s(y)$.

In particular, if  $W^u(x)$ and $W^s(y)$ have non-empty
transverse intersection, then $W^u(x)\cap W^s(y)$ is an immersed
finite dimensional submanifold of $M$, and the dimension of the
component $W_p$ of $W^u(x)\cap W^s(y)$ containing $p$ depends only on
the homotopy class of the curve $t\mapsto \phi(t,p)$, and
\[
\dim W_p \equiv m(x,\mathcal{E}) -
m(y,\mathcal{E}) \quad \mod k.
\]
\end{thm}

In particular, when $F$ satisfies (C1-2) with respect to a
$(0)$-essential subbundle $\mathcal{E}$, then all the components
of the transverse intersection $W^u(x)\cap W^s(y)$ have the same
dimension $m(x,\mathcal{E})-m(y,\mathcal{E})$.

\medskip

\proof Let us fix two points $p_0,p_1\in W^u(x)\cap W^s(y)$ such
that their orbits are homotopic in the space of paths
$u:\overline{\R}\rightarrow M$ with $u(-\infty)=x$,
$u(+\infty)=y$. In other words, there exists a continuous map
\[
h:\overline{\R} \times [0,1] \rightarrow M,
\]
such that $h(-\infty,s)=x$, $h(+\infty,s)=y$,
$h(t,i)=\phi(t,p_i)$, for $i=0,1$.

The map $\overline{\R} \times [0,1]\rightarrow \mathrm{Gr}_{(k)}(TM)$,
$(t,s) \mapsto \mathcal{E}(h(t,x))$, is liftable to a map
$\mathcal{W}:\overline{\R} \times [0,1]\rightarrow \mathrm{Gr}(TM)$ such
that $\mathcal{W}(-\infty,\cdot)$ and
$\mathcal{W}(+\infty,\cdot)$ are constant maps. By Proposition
\ref{stab}, $T_{h(t,s)} W^u(x)$ is a compact perturbation of
$\mathcal{W}(t,s)$ and $\dim (T_{h(t,s)} W^u(x),
\mathcal{W}(t,s))= \dim(H^u_x,\mathcal{W}(-\infty,\cdot))$, for
any $t<+\infty$. Using an argument analogous to the proof of
Corollary \ref{coro} (ii), it is easy to see that $(T_{h(t,s)}
W^s(y), \mathcal{W}(t,s))$ is a Fredholm pair of index $-\dim
(H^u_y,\mathcal{W}(+\infty,\cdot))$, for any $t>-\infty$. Then by
Proposition \ref{indice}, $(T_{h(t,s)} W^s(y), T_{h(t,s)} W^u(x))$
is a Fredholm pair of index
$\dim(H^u_x,\mathcal{W}(-\infty,\cdot)) -   \dim
(H^u_y,\mathcal{W}(+\infty,\cdot))$. In particular, $(T_{p_0}
W^s(y), T_{p_0} W^u(x))$ and $(T_{p_1} W^s(y), T_{p_1} W^u(x))$
are Fredholm pairs of the same index
\begin{eqnarray*}
\ind (T_{p_0} W^s(y), T_{p_0} W^u(x)) = \ind (T_{p_1} W^s(y),
T_{p_1} W^u(x))\\ = \dim(H^u_x,\mathcal{W}(-\infty,\cdot)) -
\dim (H^u_y,\mathcal{W}(+\infty,\cdot)),
\end{eqnarray*}
and the above formula implies (\ref{lequa}). Finally, the
statements which assume transversality follow from the fact that,
under such an assumption,
\[
\ind (T_p W^u(x),T_p W^s(y)) = \dim T_p W^u(x) \cap T_p W^s(y).
\]
\qed

\begin{rem}
We wish to remark that (C1-2) are asymmetric conditions: if $F$
satisfies (C1-2) with respect to a subbundle $\mathcal{V}$, there
need not exist a subbundle $\mathcal{W}$ such that $-F$ satisfies
(C1-2) with respect to $\mathcal{W}$. Such an asymmetry is
reflected into Corollary \ref{coro}, which states that $T W^u(x)$
is a compact perturbation of $\mathcal{V}$ - a closed condition - while $T
W^s (x)$ is in Fredholm pair with $\mathcal{V}$ - an open
condition. If we symmetrize (C1-2) we obtain the following
stronger assumptions: if $\mathcal{P}$ is a projector on $TM$,
with image $\mathcal{V}$ and kernel $\mathcal{W}$, there holds
(C1'): $H^u_x$ is a compact perturbation of $\mathcal{V}(x)$,
$H^s_x$ is a compact perturbation of $\mathcal{W}(x)$ for every $x\in \rest
(F)$, and (C2'): $(L_F \mathcal{P})(p)$ is compact for any $p\in M$.
This setting - actually its essential version - is closer to the
setting of a polarized manifold (see \cite{cjs95}).
\end{rem}

\section{Which manifolds can be obtained as $\mathbf{W^u(x)\cap W^s(y)}$}

\paragraph{Arbitrary gradient-like vector fields.}
Let $F$ be a gradient-like Morse-Smale vector field on the Hilbert
manifold $M$, with Lyapunov function $f$. 
If $x,y\in \rest(F)$ and $f(y)<c<f(x)$, the set $Z=
W^u(x)\cap W^s(y)  \cap f^{-1}(\{c\})$ is a submanifold (non
necessarily closed), being the transverse intersection in
$f^{-1}(\{c\})$ of the submanifolds $W^u(x)\cap f^{-1}(\{c\})$ and
$W^s(y) \cap f^{-1}(\{c\})$, and $\phi$ subordinates a
diffeomorphism from $\R \times Z$ onto $W^u(x)\cap W^s(y)$.

When $M$ is finite dimensional, there are limitations on the
topological type of $Z$, e.g.\ if $M$ is compact and there are no
rest points $z$ with $f(y)<f(z)<f(x)$, then $Z$ is the transverse
intersection in $f^{-1}(\{c\})$ of two spheres. When $M$ is
infinite dimensional, and the rest points $x,y$ have infinite
Morse index and co-index, there are no limitations at all on the
topological type of $Z$, the reason being that any manifold can be
obtained as the transverse intersection of two infinite
dimensional spheres.

More precisely, for any Hilbert manifold $Z$ (finite dimensional or not) 
there is a gradient
like Morse vector field $F$ on the Hilbert space $H$, with a
non-degenerate Lyapunov function $f$, having exactly two rest
points $x,y$ with $f(y)<0<f(x)$, such that $W^u(x)\cap
f^{-1}(\{0\})$ and $W^s(y)\cap f^{-1}(\{0\})$ are infinite
dimensional spheres intersecting transversally in $f^{-1}(\{0\})$
at a closed submanifold diffeomorphic to $Z$. Notice that in this
case, the closure of $W^u(x)\cap W^s(y)$ is $(W^u(x)\cap
W^s(y))\cup \{x,y\}$, which is homeomorphic to the suspension of
$Z$.

The construction relies on the following lemma.

\begin{lem}
Let $Z$ be a closed infinite codimensional submanifold of a
Hilbert manifold $M$. Then there exists a smooth map $\varphi:M
\rightarrow H$ such that $0$ is a regular value and
$Z=\varphi^{-1}(\{0\})$.
\end{lem}

\proof A suitable tubular neighborhood $U$ of $Z$ is diffeomorphic
to the normal bundle of $Z$. Since $Z$ has infinite codimension,
such a bundle is trivial. Therefore, there exists a submersion
$\psi:U\rightarrow H$ such that $\psi^{-1}(\{0\})=Z$. Since
$H\setminus \overline{B}$, $B$ denoting the open unit ball of $H$,
is diffeomorphic to $H$ (see \cite{bes66}), it is easy to define a
smooth map $\varphi:M\rightarrow H$ which agrees with $\psi$ on a
neighborhood $U^{\prime}\subset U$ of $Z$ and such that
$\varphi(M\setminus U^{\prime})\subset H\setminus \overline{B}$, so
that $\varphi^{-1}(\{0\})=Z$. \qed

Let $F_0$ be the vector field on $H\times H$
\[
F_0(\xi,\eta) = (\xi,-\chi(\|\xi\|) \eta),
\]
where $\chi\in C^{\infty}(\R)$ is decreasing, $\chi(s)=1$ for
$s\leq 1/3$ and $\chi(s)=0$ for $s\geq 2/3$. The vector field $F_0$ has a
unique rest point $o=(0,0)$, with $W^u(o)= H \times \{0\}$, and
has a non-degenerate Lyapunov function
\[
f_0(\xi,\eta) = 1 - \|\xi\|^2 + \chi(\|\xi\|) \|\eta\|^2.
\]
Let $B$ be the open unit ball of $H$ and let $S$ be its boundary.
We can embed $Z$ as a closed infinite codimensional submanifold of
$S$. By the above lemma, there exists a smooth map $\varphi:S
\rightarrow H$ such that $0$ is a regular value and
$Z=\varphi^{-1}(\{0\})$. Let $C_1$ and $C_2$ be two copies of the
Hilbert manifold with boundary $\overline{B} \times H$, and
consider the Hilbert manifold $M:= C_1 \cup_{\Phi} C_2$, where the
gluing map $\Phi$ is
\[
\Phi: \partial C_1 = S \times H \rightarrow S \times H = \partial
C_2, \quad (\xi,\eta) \mapsto (\xi,\eta+\varphi(\xi)).
\]
Let $F$ be the vector field on $M$ coinciding with $F_0$ on $C_1$
and with $-F_0$ on $C_2$, and let $f:M \rightarrow \R$ be the
function coinciding with $f_0$ on $C_1$ and with $-f_0$ on $C_2$.
It is readily seen that $F$ and $f$ are well-defined and smooth,
and that $f$ is a non-degenerate Lyapunov function for $F$. By
construction, $C_1$ is negatively invariant for the flow of $F$,
$C_2$ is positively invariant, and there are exactly two rest
points $x=(0,0)\in C_1$ and $y=(0,0)\in C_2$. Moreover,
$f^{-1}(\{0\})= \partial C_1 = \partial C_2$, and
\begin{eqnarray*}
W^s(y) \cap \partial C_2 = S \times \{0\}, \\
W^u(x) \cap \partial C_2 = \Phi (W^u(x) \cap \partial C_1) = \Phi
(S\times \{0\}) = \graf \varphi.
\end{eqnarray*}
Hence
\[
W^s(y) \cap W^u(x) \cap \partial C_2 = (S \times \{0\}) \cap \graf
\varphi = \varphi^{-1}(\{0\}) \times \{0\} = Z \times \{0\},
\]
the intersection being transversal, as required. Finally, since
the gluing map $\Phi$ is isotopic to the identity map on $S \times
H$, $M$ is diffeomorphic to $(\overline{B} \times H)
\cup_{\mathrm{id}} (\overline{B} \times H)= (\overline{B}
\cup_{\mathrm{id}} \overline{B}) \times H$. Being an infinite
dimensional sphere, $\overline{B} \cup_{\mathrm{id}} \overline{B}$
is diffeomorphic to $H$ (again by \cite{bes66}), hence $M$ is
diffeomorphic to $H$.

\paragraph{Gradient-like vector fields satisfying (C1-2).}
In particular, if $Z$ has components of different dimension, the
above example shows that in the case of infinite Morse indices and
co-indices, the components of $W^u(x) \cap W^s(y)$ may have
different dimension. Actually, the discussion of section
\ref{sstab} suggests that this phenomenon may happen also when $F$
satisfies (C1-2) with respect to an essential subbundle, provided
$M$ is not simply connected. Indeed this is the case, as it is
shown by the following example, where the vector field is actually
the gradient of a smooth function.

We recall some pieces of notation from Appendix B. If
$A:[-\infty,+\infty]\rightarrow \mathcal{L}(H)$ is a path of
bounded linear operators with $A(-\infty)$ and $A(+\infty)$
hyperbolic, we denote by $X_A:\R \rightarrow \mathrm{GL}(H)$ the solution
of the linear Cauchy problem $X_A^{\prime}(t)=A(t)X_A(t)$,
$X_A(0)=I$, and we consider the closed linear subspaces
\[
W^s_A = \set{\xi\in H}{\lim_{t\rightarrow +\infty} X_A(t)\xi =0},
\quad W^u_A = \set{\xi\in H}{\lim_{t\rightarrow -\infty} X_A(t)\xi
=0}.
\]
Let $M=H \times \T^1$, with $\T^1=\R/2\pi\Z$. Let $H=H^-\oplus
H^+$ be an orthogonal splitting such that $H^-,H^+\in
\mathrm{Gr}_{\infty,\infty}(H)$, with associated projectors $P^-,P^+$. Let
$k\geq 1$ be an integer. By Proposition \ref{costra} there exists
$A\in C^{\infty}(\R,\mathrm{GL}(H) \cap \mathrm{Sym}(H))$ with $A(t)=P^+-P^-$ for
$t\notin ]0,1[$, such that $W^s_A+W^u_A=H$ and $\dim W^s_A \cap
W^u_A = k$. Consider the smooth tangent vector field on $M$
\[
F(\xi,s) = \begin{cases} ((P^+ -P^-)\xi,\sin s) & \mbox{for }
s\notin
  [\pi/2,\sigma(1)], \\ (A(\tau(s))\xi,\sin s) & \mbox{for } s\in
  ]0,\pi[, \end{cases} \quad (\xi,s)\in H \times \T^1,
\]
where $\tau(s) = \log \tan (s/2)$ for $0<s<\pi$, and $\sigma(t) =
\tau^{-1}(t) = 2\arctan e^t$, $t\in \R$. So $\sigma^{\prime} =
\sin \sigma$ and $\tau^{\prime} = \cosh \tau$.

The rest points of $F$, $x=(0,0)$ and $y=(0,\pi/2)$, are
hyperbolic, with stable and unstable linear spaces
\begin{equation}
\label{ut} H^s_x = H^- \times (0), \quad H^u_x = H^+ \times \R,
\quad\quad H^s_y = H^- \times \R, \quad H^u_y = H^+ \times (0).
\end{equation}
The flow of $F$, $\phi:\R \times M \rightarrow M$, is readily
seen to be
\begin{equation}
\label{dt} \phi_t(\xi,s) = \begin{cases} (e^{t(P^+-P^-)}
\xi,-\sigma(t+\tau(s))) &
  \mbox{for } -\pi<s<0, \\ (e^{t(P^+-P^-)} \xi,s)  &
  \mbox{for $s=0$ or $s=\pi$}, \\ (X_A(t+\tau(s)) X_A(\tau(s))^{-1} \xi,
  \sigma(t+\tau(s))) & \mbox{for } 0<s<\pi. \end{cases}
\end{equation}
As a consequence, the unstable manifold of $x$ and the stable
manifold of $y$ are the sets
\begin{eqnarray*}
W^u(x) = (H^+ \times ]-\pi,0]) \cup \bigcup_{0<s<\pi} X_A(\tau(s))
      W^u_A \times \{s\}, \\
W^s(x) = (H^- \times [-\pi,0[) \cup \bigcup_{0<s<\pi} X_A(\tau(s))
      W^s_A \times \{s\},
\end{eqnarray*}
with tangent spaces
\begin{eqnarray*}
T_{(\xi,s)} W^u(x) = \begin{cases} H^+ \times \R & \mbox{for }
(\xi,s)
  \in H^+ \times ]-\pi,0], \\ X_A(\tau(s)) W^u_A \oplus \R F(\xi,s) &
  \mbox{otherwise}, \end{cases} \\
T_{(\xi,s)} W^s(y) = \begin{cases} H^- \times \R & \mbox{for }
(\xi,s)
  \in H^- \times [-\pi,0[, \\ X_A(\tau(s)) W^s_A \oplus \R F(\xi,s) &
  \mbox{otherwise}. \end{cases}
\end{eqnarray*}
Therefore, the unstable manifold of $x$ and the stable manifold
of $y$ meet transversally, and their intersection
\[
W^u(x) \cap W^s(y) = (\{0\}\times ]-\pi,0[)\cup \bigcup_{0<s<\pi}
X_A(\tau(s))
      (W^u_A\cap W^s_A) \times \{s\},
\]
consists of two connected components, one of which
one-dimensional, the other one of dimension $k+1$.

We are going to show that the vector field $F$ satisfies
conditions (C1) and (C2) with respect to a (non-liftable)
essential subbundle of $TM$. Consider the two subbundles of
$T(H\times [-\pi,0])$ and $T(H\times ]0,\pi[)$,
\begin{eqnarray*}
V_1(\xi,s) = H^+ \times (0) \quad \mbox{for } (\xi,s)\in H\times
[-\pi,0], \\ V_2(\xi,s) = X_A(\tau(s)) Y \times (0) \quad
\mbox{for } (\xi,s)\in H\times ]0,\pi[,
\end{eqnarray*}
where $Y$ is a closed supplement of $W^s_A \cap W^u_A$ in
$W^u_A=H^+$. Since $A(\tau(s))=P^+-P^-$ for $0<s\leq \pi/2$,
$V_2(\xi,s)=Y \times (0)$ for any $(\xi,s)\in H\times
    ]0,\pi/2]$. Moreover, since $H=W^s_A \oplus Y$, by Theorem
\ref{abue} (iii),
\[
\dist (V_2(\xi,s),H^+ \times (0)) =
\dist(X_A(\tau(s))Y,V^+(A(+\infty))) \rightarrow 0 \quad
\mbox{for } s\rightarrow \pi-.
\]
Therefore $V_1$ and $V_2$ define a $C^0$ essential subbundle
$\mathcal{E}$ of $TM$. In order to show that $\mathcal{E}$ is of
class $C^1$, we have to verify that
\begin{equation}
\label{tt} \frac{d}{ds} P_{X_A(\tau(s))Y} = \tau^{\prime}(s)
Q^{\prime}(\tau(s)) = \cosh \tau(s) Q^{\prime}(\tau(s))
\rightarrow 0 \quad \mbox{for } s\rightarrow \pi-,
\end{equation}
where $Q(t) = P_{X_A(t)Y}$. For $t_0\geq 1$ large,
$X_A(t_0)Y\subset H^+ \times H^-$ is the graph of some operator
$L\in \mathcal{L}(H^+,H^-)$, so
\[
X_A(t) Y = X_{P^+ -P^-} (t-t_0) X_A(t_0) Y = \Bigl\{ \left(
\begin{array}{cc}
    e^{t-t_0} & 0 \\ 0 & e^{t_0-t} \end{array} \right) \left(
    \begin{array}{c} \xi \\ L \xi \end{array} \right) \, \Big| \,
    \xi\in H^+ \Bigr\},
\]
from which we deduce that $Q(t) - P^+ = O(e^{-2t})$ for
$t\rightarrow +\infty$. By identity (\ref{ricceq}), $Q$ solves
the Riccati equation
\[
Q^{\prime} = (I-Q)AQ + QA(I-Q),
\]
and since $A(t) = P^+-P^-$ for $t\geq 1$, we obtain
\[
Q^{\prime} (t) = 2 (Q(t)-P^+) P^- (Q(t)-I) +2 (Q(t)- I) P^- (Q(t)
-P^+) = O(e^{-2t}) \quad \mbox{for } t\rightarrow +\infty,
\]
which proves (\ref{tt}).

By (\ref{ut}) the vector field $F$ satisfies (C1) with respect to
the essential subbundle $\mathcal{E}$. By (\ref{dt}),
\[
D\phi_t(\xi,s)[(\eta,0)] = \begin{cases} (e^{t(P^+-P^-)} \eta,0) &
  \mbox{for } \pi\leq s \leq 0, \\ (X_A(t+\tau(s))X_A(\tau(s))^{-1}
  \eta,0) & \mbox{for }0<s<\pi, \end{cases}
\]
for every $t\in \R$, $(\xi,s)\in M$, $\eta\in H$. Therefore, the
subbundle $V_1$ is invariant with respect to $F$. Since
$X_A(t+\tau)X_A(\tau)^{-1} = X_{A(\cdot+\tau)}(t)$, also the
subbundle $V_2$ is invariant with respect to $F$. Hence $(L_F
P_{V_i})P_{V_i} =0$, for $i=1,2$, and $F$ satisfies (C2) with
respect to the essential subbundle $\mathcal{E}$.

The smooth function
\[
f(\xi,s) = \begin{cases} -\frac{1}{2} \langle (P^+-P^-) \xi,\xi
  \rangle + \cos s & \mbox{for } s\notin [\pi/2,\sigma(1)], \\
-\frac{1}{2} \langle A(\tau(s)) \xi,\xi
  \rangle + \cos s & \mbox{for } s\in ]0,\pi/2[, \end{cases}
\]
satisfies
\[
Df(\xi,s)[F(\xi,s)] = \begin{cases} -\|\xi\|^2 - \sin^2 s &
\mbox{for } s\notin [\pi/2,\sigma(1)], \\ - \|A(\tau(s)) \xi\|^2
-\sin^2 s - \frac{1}{2} \cosh \tau(s) \langle A^{\prime}(\tau(s))
\xi,\xi \rangle
 & \mbox{for } s\in ]0,\pi/2[. \end{cases}
\]
Since $A(t)$ is invertible for any $t$ and $A(t)$ is constant for
$t\notin (0,1)$, we find
\begin{equation}
\label{le} Df(\xi,s) [F(\xi,s)] \leq -\delta \|\xi\|^2 - \sin^2 s
\quad \mbox{for } \|\xi\|<r,
\end{equation}
for suitable $\delta>0$, $r>0$, so $f$ is a non-degenerate
Lyapunov function for $F$ on the open set $B_r(0) \times \T^1$.
Actually, on $B_r(0)\times \T^1$ the vector field $F$ is the
gradient of $-f$ with respect to a smooth metric of the form
$\alpha_p (\zeta_1,\zeta_2) = \langle T(p) \zeta_1,\zeta_2
\rangle$, $p\in M$, $\zeta_1,\zeta_2\in T_p M = H \oplus \R$. Here
$T\in C^{\infty}(B_r(0)\times \T^1, \mathrm{Sym} \cap \mathrm{GL}(H \oplus \R))$ is
positive, and verifies
\[
T(\xi,s) = I \quad \mbox{for } s\notin [\pi/2,\sigma(1)], \quad
T(p) F(p) = - \grad f(p)\quad \mbox{for } p\in B_r(0) \times \T^1,
\]
where $\grad f$ denotes the gradient of $f$ with respect to the
flat metric on $H\times \T^1$. Such a map $T$ can be easily found
because by (\ref{le}), $\langle F(p),-\grad f(p) \rangle >0$ for
every $p\in B_r(0) \times [\pi/2,\sigma(1)]$.

\section{Orientation of $\mathbf{W^u(x)\cap W^s(y)}$}
\label{oriente}

The first example of the previous section shows that the
transverse intersection of an unstable and a stable manifold of
two rest points with infinite Morse index and co-index, even if
finite dimensional, needs not be orientable. This intersection
will be orientable when the vector field satisfies (C1-2) with
respect to a subbundle of $TM$.

\paragraph{Orientation of Fredholm pairs.}
We need to recall some facts about the orientation bundle on the
space of Fredholm pairs. See Appendix A for more details. For $H$
a real Hilbert space and $n\in \N$, we denote by
$\mathrm{Or}(\mathrm{Gr}_{n,\infty}(H))\rightarrow \mathrm{Gr}_{n,\infty}(H)$ the orientation
bundle of the Grassmannian of $n$-dimensional subspaces of $H$:
the fiber of $X\in \mathrm{Gr}_{n,\infty}(H)$ consists of the two
orientations of $X$. Similarly, $\mathrm{Or}(\mathrm{Fp}(H))\rightarrow \mathrm{Fp}(H)$
denotes the orientation bundle of the space of Fredholm pairs: the
fiber of $(V,W)\in \mathrm{Fp}(H)$ consists of the two orientations of the
finite dimensional vector space $(V\cap W) \times (H/(V+W))^*$.
This bundle is actually a double covering of $\mathrm{Fp}(H)$. If
$\mathcal{H}\rightarrow B$ is a Hilbert bundle, we obtain the
bundles $\mathrm{Or}(\mathrm{Gr}_{n,\infty}(\mathcal{H})) \rightarrow B$,
$\mathrm{Or}(\mathrm{Fp}(\mathcal{H}))\rightarrow B$, where
\[
\mathrm{Or}(\mathrm{Gr}_{n,\infty}(\mathcal{H}) = \bigcup_{b\in B} \mathrm{Or}
(\mathrm{Gr}_{n,\infty}(\mathcal{H}_b)), \quad \mathrm{Or}(\mathrm{Fp}(\mathcal{H}) =
\bigcup_{b\in B} \mathrm{Or} (\mathrm{Fp}(\mathcal{H}_b)),
\]
and the maps
\[
\mathrm{Or}(\mathrm{Gr}_{n,\infty}(\mathcal{H})) \rightarrow
\mathrm{Gr}_{n,\infty}(\mathcal{H}), \quad \mathrm{Or}(\mathrm{Fp}(\mathcal{H})) \rightarrow
\mathrm{Fp}(\mathcal{H})
\]
are double coverings.

Let $n\in \N$. If $\mathcal{S}(n,\mathrm{Fp})$ denotes the open set
consisting of all $(X,(V,W))$ in $\mathrm{Gr}_{n,\infty}(H) \times \mathrm{Fp}(H)$
such that $X\cap V=(0)$, the map $\mathcal{S}(n,\mathrm{Fp}) \rightarrow
\mathrm{Fp}(H)$, $(X,(V,W)) \mapsto (X\oplus V,W)$, is continuous, and it
lifts to a continuous map - the product of orientations -
\[
(o_X, o_{(V,W)}) \mapsto o_X {\textstyle \bigwedge} o_{(V,W)},
\]
from the corresponding open subset of $\mathrm{Or}(\mathrm{Gr}_{n,\infty}(H)) \times
\mathrm{Or}(\mathrm{Fp}(H))$ to $\mathrm{Or}(\mathrm{Fp}(H))$. If $\alpha:B\rightarrow
\mathrm{Gr}_{n,\infty}(\mathcal{H})$, $\beta,\gamma: B \rightarrow
\mathrm{Fp}(\mathcal{H})$ are continuous sections such that
\[
\alpha(b)\cap \beta_1(b)=(0) \mbox{ and } (\alpha(b)\oplus
\beta_1(b),\beta_2(b)) = \gamma(b),\quad \forall b\in B,
\]
for any choice of liftings of two of $\alpha$, $\beta$, $\gamma$
to the orientation bundles, there is a unique lifting of the third
one such that $\hat{\alpha}(b) {\textstyle \bigwedge}
\hat{\beta}(b) = \hat{\gamma}(b)$ for every $b\in B$.

\paragraph{Orientation of $\mathbf{W^u(x)\cap W^s(y)}$.}
Let $\mathcal{V}$ be a $C^1$ subbundle of $TM$, and let us assume
that the Morse vector field $F$ satisfies (C1-2) with respect to
$\mathcal{V}$. By (C1) for every rest point $x$ the pair
$(H^s_x,\mathcal{V}(x))$ is a Fredholm pair. Let us fix arbitrarily
an orientation $o_x$ of $(H^s_x,\mathcal{V}(x))$. Let $x,y\in
\rest(F)$ be such that $W^u(x)$ and $W^s(y)$ have a non-empty and
transverse intersection. By Theorem \ref{ubu}, $W^u(x)\cap W^s(y)$
is an immersed submanifold of dimension $n=m(x,\mathcal{V}) -
m(y,\mathcal{V})$. In this section we will prove that
$W^u(x)\cap W^s(y)$ is orientable, and we will show how an
orientation of such a manifold is determined by the orientations
$o_x$ and $o_y$.

Let $h_x^u:H^u_x \rightarrow M$ and $h^s_y:H^s_y \rightarrow M$ be
injective $C^1$ immersions onto $W^u(x)$ and $W^s(y)$ such that
$h^u_x(0)=x$ and $h^s_y(0)=y$. Then $W^u(x)\cap W^s(y)$ is the
image of the injective immersion $h=h^u_x \circ p^u = h^s_y \circ
p^s : W\rightarrow M$ coming from the fiber product square of the
transverse maps $h^u_x$ and $h^s_y$:
\[
\xymatrix{W \ar^{p^s}[r] \ar_{p^u}[d] \ar^{h}[rd] & H^s_y
  \ar^{h^s_y}[d] \\ H^u_x \ar^{h^u_x}[r] & M,} \quad W =
  \set{(\xi,\eta)\in H^u_x \times H^s_y }{h^u_x(\xi) = h^s_y
  (\eta)}.
\]
Giving an orientation to $W^u(x)\cap W^s(y)$ is equivalent to
lifting the section
\[
\tau:W \rightarrow \mathrm{Gr}_{n,\infty}(h^*(TM)), \quad w\mapsto T_{h(w)}
(W^u(x) \cap W^s(y)),
\]
to a section $\hat{\tau}:W \rightarrow
\mathrm{Or}(\mathrm{Gr}_{n,\infty}(h^*(TM)))$.

Since $H^s_y$ is contractible, the section
\[
H^s_y \rightarrow \mathrm{Fp}({h^s_y}^*(TM)), \quad \eta \mapsto
(T_{h^s_y(\eta)} W^s(y), \mathcal{V}(h^s_y(\eta))),
\]
has a unique lifting $H^s_y \rightarrow \mathrm{Or}(\mathrm{Fp}({h^s_y}^*(TM)))$
whose value at $0$ is $o_y$. By composition with the projection
$p^s:W\rightarrow H^s_y$, we obtain the section
\[
\omega:W \rightarrow \mathrm{Fp}(h^*(TM)), \quad w\mapsto (T_{h(w)}
W^s(y),\mathcal{V}(h(w))),
\]
and a lifting of $\omega$, $\hat{\omega}:W \rightarrow
\mathrm{Or}(\mathrm{Fp}(h^*(TM)))$.

Let $Y:W \rightarrow \mathrm{Gr}(h^*(TM))$ be a continuous section of
linear supplements of $T(W^u(x) \cap W^s(y))$ in $TW^s(y)$, that
is
\[
T_{h(w)} W^s(y) = T_{h(w)} (W^u(x) \cap W^s(y)) \oplus Y(w), \quad
\forall w\in W.
\]
By the transversality of the intersection  of $W^u(x)$ and
$W^s(y)$, $Y(w)$ is a linear supplement of $T_{h(w)} W^u(x)$ in
$T_{h(w)} M$, so by Theorem \ref{abue} (iii),
\[
\lim_{t\rightarrow -\infty} D\phi_t(h(w)) Y(w) = H^s_x,
\]
and the limit is locally uniform in $W$. Therefore the map $A:
[-\infty,0] \times W \rightarrow \mathrm{Fp}(TM)$ defined by
\[
A(t,w) = \begin{cases} (D \phi_t(h(w))
  Y(w),\mathcal{V}(\phi_t(h(w)))) & \mbox{for } t>-\infty, \\
  (H^s_x,\mathcal{V}(x)) & \mbox{for } t=-\infty, \end{cases}
\]
is continuous. Let $\hat{A}:[-\infty,0] \times W \rightarrow
\mathrm{Or}(\mathrm{Fp}(TM))$ be the unique lifting of $A$ such that
$\hat{A}(-\infty,w) = o_x$ for any $w\in W$. By restriction to
$\{0\} \times W$, we obtain the section
\[
\alpha:W \rightarrow \mathrm{Fp}(h^*(TM)), \quad w\mapsto (Y(w),
\mathcal{V}(h(w))),
\]
and a lifting of $\alpha$, $\hat{\alpha}:W \rightarrow
\mathrm{Or}(\mathrm{Fp}(h^*(TM)))$.

By construction,
\[
\tau(w) \cap \alpha_1(w) = (0) \mbox{ and } (\tau(w)\oplus
\alpha_1(w), \alpha_2(w)) = \omega(w), \quad \forall w\in W,
\]
so $\tau$ has a unique lifting $\hat{\tau}: W \rightarrow
\mathrm{Or}(\mathrm{Gr}_{n,\infty}(h^*(TM)))$ such that
\[
\hat{\tau}(w) {\textstyle \bigwedge}  \hat{\alpha}(w) =
\hat{\omega}(w) , \quad \forall w\in W,
\]
which provides us with an orientation of $W^u(x)\cap W^s(y)$.

Since the set of linear supplements of  $T_p (W^u(x) \cap W^s(y))$
in $T_p W^u(x)$ is connected, the orientation we have defined
does not depend on the choice of $Y$. The construction shows that
it does not depend on the choice of the immersions $h^u_x$ and
$h^s_y$.

\section{Compactness of $\mathbf{W^u(x)\cap W^s(y)}$}
\label{scomp}

\paragraph{The Palais-Smale condition.}
Let $F$ be a gradient-like Morse vector field on $M$. Then the stable 
and unstable manifolds
are (embedded) submanifolds, and so are their transverse
intersections. We would like to state a set of assumptions which
imply that $W^u(x)\cap W^s(y)$ is pre-compact, i.e.\ it 
has compact closure in $M$. The
first assumption is a version of the well known Palais-Smale
condition:

\begin{defn}
\label{psd}
Let $F$ be a $C^1$ vector field on $M$, and $f\in C^1(M)$ be a
Lyapunov function for $F$. A {\em (PS) sequence} for $(F,f)$ is a
sequence $(p_n)\subset M$ with $f(p_n)$ bounded and
$Df(p_n)[F(p_n)]\rightarrow 0$. The pair $(F,f)$ satisfies the
{\em (PS) condition} if every (PS) sequence is compact. We shall
say that $F$ satisfies (PS) if $(F,f)$ satisfies (PS) for some
Lyapunov function $f$.
\end{defn}

When $F$ is the negative gradient of a function $f$ with respect to
some Riemannian metric on $M$, one finds the usual notion of a
Palais-Smale sequence: $f(p_n)$ bounded and $\|Df(p_n)\|$ infinitesimal.

Since $Df(p)[F(p)]< 0$ for $p\notin \rest(F)$, the limit points of
the (PS) sequences are rest points of $F$. If $(F,f)$ satisfies
(PS), then the set $\rest(F) \cap \{a\leq f \leq b\}$ is compact
for every $a,b\in \R$. If moreover $F$ is a Morse vector field,
this set is also discrete, hence finite.

\begin{rem} {\em (Genesis of (PS) sequences)}
\label{psr} Let $(t_n,p_n)\in \Omega(F)$ be such that
$t_n\rightarrow +\infty$, and $f(p_n)$, $f(\phi(t_n,p_n))$ are
bounded. Then by the mean value theorem there exists 
$s_n\in [0,t_n]$ such that $(\phi(s_n,p_n))$ is a (PS) sequence 
for $(F,f)$.
\end{rem}

The second assumption is the completeness of $F$, that is 
the fact that the local flow $\phi(t,\cdot)$ of $F$ is defined for 
every $t$, i.e.\ $\Omega(F)= \R \times M$.

\begin{rem}
Notice that multiplying $F$ by a positive function does not change
$W^u(x)\cap W^s(y)$, whereas it may have an effect on the validity of
the (PS) and the completeness assumption. For instance, multiplication
by a suitably small function makes the vector field complete, while
multiplication by a suitable large function makes (PS) true, when
$\rest (F) \cap f^{-1}([a,b])$ is compact for every $a,b\in \R$. The
two assumptions are meaningful here only when considered together.
\end{rem}

\paragraph{Essentially vertical families.}
As we shall see, the (PS) condition and the completeness
imply that $W^u(x)\cap
W^s(y)$ has compact closure, when either all the rest points of $F$ have
finite Morse index, or they have finite Morse co-index. However, they
are not sufficient in the general case. 

The first assumption we need in the general case is that
the essential subbundle $\mathcal{E}$ of $TM$ should have a strong
integrable structure $\mathcal{A}$ modeled on $(H,V)$ (see 
Definition \ref{sintr}). In this case, denoting by $Q$ a linear
projector with kernel $V$, we can introduce the following:

\begin{defn}
\label{evf}
A family $\mathcal{F}$ of subsets of $M$ is called an
essentially vertical family for the strong integrable structure
$\mathcal{A}$ of $\mathcal{E}$ if it satisfies:
\begin{enumerate}
\item it is an ideal of $\mathcal{P}(M)$, meaning that it is closed
  for finite unions and if $A\in \mathcal{F}$ then any subset of $A$
  is also in $\mathcal{F}$;
\item every point $p$ has a neighborhood $U$ which is the domain of a 
chart $\varphi\in \mathcal{A}$ such that every set
$A\subset U$ with $\varphi(A)$ bounded 
belongs to $\mathcal{F}$ if and only if $Q\varphi(A)$ is pre-compact. 
\end{enumerate}
\end{defn}

Once an essentially vertical family $\mathcal{F}$ has been fixed, its
elements will be called essentially vertical sets. 
Clearly, there can be many different essentially vertical families 
associated to the same strong integrable structure of $\mathcal{E}$,
because only the ``small'' essentially vertical subsets are determined. 

The family $\mathcal{F}$ will be called {\em positively invariant} 
if it is closed under the positive action of the flow $\phi$: 
for every $A\in \mathcal{F}$ and every $t\geq 0$, the set
$\phi([0,t]\times A)$ is in $\mathcal{F}$.

The main result of this section is the following compactness theorem.

\begin{thm}
\label{tcomp} Assume that the Morse gradient-like vector field $F$ is 
complete, satisfies (C1) with respect to an essential
subbundle $\mathcal{E}$ of $TM$, and satisfies (PS). 
Assume also  that $\mathcal{E}$ has a strong integrable structure 
$\mathcal{A}$ modeled on $(H,V)$ and an essentially vertical family
$\mathcal{F}$, which is positively invariant for the flow of $F$.

Let  $(p_n)\subset M$, $(s_n)\subset (-\infty,0]$,
$(t_n)\subset [0,+\infty)$, be such that $(\phi(s_n,p_n))$
converges to a rest point $x$, while $(\phi(t_n,p_n))$ converges to a 
rest point $y$. Then the sequence $(p_n)$ is compact.
\end{thm}
 
An immediate consequence is the following corollary.

\begin{cor}
Under the assumptions of Theorem \ref{tcomp}, for every $x,y\in
\rest(F)$ the intersection $W^u(x)\cap W^s(y)$ is pre-compact.
\end{cor}

In order to prove the above theorem, we need to establish the
following lemma.

\begin{lem}
\label{lmlm}
Let $x$ be a rest point of $F$. Then $x$ has a fundamental system of
neighborhoods $U$ such that:
\begin{enumerate}
\item the set $W^u(x)\cap U$ is essentially vertical;
\item if $A\subset U$ is essentially vertical, then $A \cap
    W^s(x)$ is pre-compact.
\end{enumerate}
Furthermore, if $f$ is a non-degenerate Lyapunov function for $F$,
for any sequence $(p_n)\subset U$ converging to $x$ there holds:
\begin{enumerate}
\setcounter{enumi}{2}
\item if $t_n\geq 0$ is such that $\phi(t_n,p_n)\in \partial U$ then the set
$\set{\phi(t_n,p_n)}{n\in \N}$ is essentially vertical, and
\[
\limsup_{n\rightarrow \infty} f(\phi(t_n,p_n)) < f(x).
\]
\item if $t_n\leq 0$ is such that $\phi(t_n,p_n)\in \partial U$ then
the set $\set{\phi(t_n,p_n)}{n\in \N}$ has a pre-compact intersection
with any essentially vertical subset of $M$, and
\[
\liminf_{n\rightarrow \infty} f(\phi(t_n,p_n)) > f(x).
\]       
\end{enumerate}
\end{lem}

\proof
By choosing a chart $\varphi\in \mathcal{A}$ satisfying property (ii)
in Definition \ref{evf}, we can identify a neighborhood $U$ of $x$ in $M$
with a bounded neighborhood - still denoted by $U$ - of $0$ in $H$, in such a
way that $x$ corresponds to $0$, the essential subbundle $\mathcal{E}$ is
represented by the constant subbundle $V$, the kernel of a projector
$Q$, and the
essentially vertical subsets $A\subset U$ are those for which
$QA$ is pre-compact.

Let $H=H^u \oplus H^s$, with projections $P^u$, $P^s$, be the
splitting into the linear unstable and the stable spaces of the hyperbolic
rest point $0$. Endow $H$ with an adapted norm $\|\cdot\|$ (see
Appendix C), and denote by $H^u(r)$, $H^s(r)$ the closed $r$-balls of
$H^u$ and $H^s$, respectively. Up to reducing the neighborhood $U$, we
may assume that $U=H^u(r) \times H^s(r)$, where $r>0$ is so
small that all the results of Appendix C hold.

By (C1), $H^u$ is a compact perturbation 
of $V$. Therefore, we may assume that $Q$ is a compact perturbation of
$P^s$. By Remark \ref{chv} and by the boundeness of $U$, a subset $A\subset U$
is essentially vertical if and only if $P^s A$ is pre-compact. In
particular, the graph of a map $\sigma:H^u(r) \rightarrow H^s(r)$
is essentially vertical if and only if the map $\sigma$ is compact.

Let $\sigma_0:H^u(r) \times H^s(r)$ be a 1-Lipschitz map. By the graph
transform method (Proposition \ref{grafici}), for every $t\geq 0$ the set
\[
\set{\phi(t,\xi)}{\xi \in \graf \sigma_0 \mbox{ and } \phi([0,t] \times
  \{\xi\}) \subset H^u(r) \times H^s(r)}
\]
is the graph of a 1-Lipschitz map $\sigma_t:H^u(r) \rightarrow
H^s(r)$, and $\sigma_t$ converges uniformly to a map $\sigma^u$ for
$t\rightarrow +\infty$, with $\graf \sigma^u =
W^u_{\mathrm{loc},r}(0)$, the local unstable manifold of $0$. If $\sigma_0$ is a compact map - for
example $\sigma_0=0$ - the fact that the family $\mathcal{F}$ is
positively invariant implies that $\sigma_t$ is a compact map for
every $t\geq 0$. By the uniform convergence, $\sigma^u$ is also
compact, so the local unstable manifold $W^u_{\mathrm{loc},r}(0)$ is
an essentially vertical set. By Theorem \ref{thesta},
$W^u_{\mathrm{loc},r}(0) = W^u(x) \cap U$, proving (i).

\medskip

The local stable manifold $W^s_{\mathrm{loc},r}(0)$ is the graph of a
1-Lipschitz map $\sigma^s : H^s(r) \rightarrow H^u(r)$. 
Let $A\subset U$ be an essentially vertical subset, that is $P^sA$ is
pre-compact. Then 
\[
A \cap W^s_{\mathrm{loc},r}(0) = \graf \sigma^s|_{P^s A}
\]
is also pre-compact. By Theorem \ref{thesta}, $W^s_{\mathrm{loc},r}(0)
= W^s(x) \cap U$,  proving (ii).

\medskip

Let $(p_n)\subset U$ be a sequence converging to $0$, and $t_n\geq
0$ such that $\phi(t_n,p_n) \in \partial U$. By Lemma \ref{hg},
\begin{eqnarray*}
\limsup_{n\rightarrow \infty} f(\phi(t_n,p_n)) < f(0), \\
\lim_{n\rightarrow \infty} 
\dist (\phi(t_n,p_n) , W^u_{\mathrm{loc},r}(0)\cap \partial U) = 0.
\end{eqnarray*}
The first limit proves part of assertion (iii). By the second limit,
we can find a sequence $(q_n)\subset W^u_{\mathrm{loc},r}(0)$ such
that $\|\phi(t_n,p_n)-q_n\|$ is infinitesimal. In particular $\|P^s
\phi(t_n,p_n) - P^s q_n\|$ is infinitesimal. By (i), the sequence
$(P^s q_n)$ is compact. So also the sequence $(P^s \phi(t_n,p_n))$ is
compact. This proves that the set $\set{\phi(t_n,p_n)}{n\in \N}$ is
essentially vertical, concluding the proof of (iii).

\medskip

The fact that $\sigma^s$ is 1-Lipschitz implies that
\begin{equation}
\label{pitagora}
\|P^u \xi - \sigma^s(P^s \xi)\| \leq \sqrt{2} \, \dist (\xi,
\graf \sigma^s) \quad \forall \xi \in U.
\end{equation}
Indeed, if $\xi\in U$ and $c>1$ we can find $\eta\in \graf \sigma^s$ 
such that
\[
\|\xi -\eta\| \leq c \, \dist( \xi, \graf \sigma^s).
\]
Since $P^u \eta = \sigma^s(P^s \eta)$ and since $\sigma^s$ is
1-Lipschitz,
\begin{eqnarray*}
\| P^u \xi - \sigma^s(P^s \xi) \| \leq \|P^u \xi - P^u \eta\| + \|
\sigma^s(P^s \eta) - \sigma^s(P^s \xi)\| \\
\leq \|P^u \xi - P^u \eta\| + \|P^s \eta - P^s \xi\| \leq \sqrt{2}
\|\xi -\eta\| \leq c\sqrt{2}\, \dist( \xi, \graf \sigma^s),
\end{eqnarray*}
and since $c>1$ is arbitrary, (\ref{pitagora}) follows.

Now assume that $t_n\leq 0$ are such that $\phi(t_n,p_n)\in \partial
U$. By Lemma \ref{hg} applied to $-F$,
\begin{eqnarray*}
\liminf_{n\rightarrow \infty} f(\phi(t_n,p_n)) > f(0), \\
\lim_{n\rightarrow \infty} 
\dist (\phi(t_n,p_n) , W^s_{\mathrm{loc},r}(0)\cap \partial U) = 0.
\end{eqnarray*} 
The first limit proves part of assertion (iv). 
Let $A\subset U$ be an essentially vertical set, that is $P^s A$ is 
pre-compact. If the set $A\cap \set{\phi(t_n,p_n)}{n\in \N}$ is
infinite (otherwise there is nothing to prove), its elements form a
subsequence $(\phi(t_{n_k},p_{n_k}))$ such that the sequence $(P^s
\phi(t_{n_k},p_{n_k}))$ is compact. By the continuity of $\sigma^s$,
also the sequence
$(\sigma^s(P^s \phi(t_{n_k},p_{n_k})))$ is compact. By (\ref{pitagora}),
\begin{eqnarray*}
\|P^u \phi(t_{n_k},p_{n_k}) - \sigma^s(P^s \phi(t_{n_k},p_{n_k})) \|
\leq \sqrt{2}\, \dist (\phi(t_{n_k},p_{n_k})),\graf \sigma^s) \\
\leq \sqrt{2}\, \dist (\phi(t_{n_k},p_{n_k})),W^s_{\mathrm{loc},r}(0)\cap
\partial U)
\end{eqnarray*}
is infinitesimal, so also $(P^u \phi(t_{n_k},p_{n_k}))$ is compact. We
deduce that $(\phi (t_{n_k},p_{n_k}))$ is compact, concluding the
proof of (iv).
\qed

\bigskip

{\em Proof }[of Theorem \ref{tcomp}].
Let $f$ be a non-degenerate Lyapunov function for $F$ such that 
$(F,f)$ satisfies (PS).
Up to taking a subsequence of $(p_n)$ and changing $(s_n)$ and
$x$, we may assume that for no choice of a sequence $(r_n)\subset
]-\infty,0]$, the sequence $(\phi(r_n,p_n))$ has a subsequence
which converges to a rest point $z$ with $f(z)<f(x)$. Indeed,
since there are finitely many rest points $z$ such that $f(y)\leq
f(z) \leq f(x)$, the set
\begin{eqnarray*}
\mathcal{Z} := \Bigl\{ z\in \rest(F) \, \Big| \, f(z)\leq f(x)
\mbox{ and there
    exists $(n_k)\subset \N$ increasing and $s^{\prime}_k\leq
    0$} \\ \mbox{such that } \lim_{k\rightarrow \infty}
    \phi(s_k^{\prime},p_{n_k})=z\Bigr\}
\end{eqnarray*}
is finite, and non-empty because it contains $x$. Let $x^{\prime}=
\lim_{k\rightarrow \infty} \phi(s_k^{\prime},p_{n_k})$ be a point
of $\mathcal{Z}$ where $f$ attains its minimum. Then the latter
requirement is verified with $(p_{n_k})$, $(s_k^{\prime})$, and
$x^{\prime}$.

Similarly, by taking a further subsequence of $(p_n)$, and by
changing $(t_n)$ and $y$, we may assume that for no choice of a
sequence $(r_n)\subset \R$, the sequence $(\phi(r_n,p_n))$ has a
subsequence which converges to a rest point $z$ with
$f(y)<f(z)<f(x)$. If either $x$ or $y$ is a cluster point for
$(p_n)$ there is nothing to prove, so we may assume that $(p_n)$
is bounded away from $x$ and $y$.

By Lemma \ref{lmlm} (iii), there exists a closed 
neighborhood $U\subset M$ of $x$ such that $p_n\notin U$, 
and choosing $s_n^{\prime}\in ]s_n,0[$ such that
$\phi(s_n^{\prime},p_n) \in \partial U$ (for $n$ large), we have
\begin{eqnarray}
\label{stel}
\set{\phi(s_n^{\prime},p_n)}{n\in \N} \mbox{ is essentially vertical,}\\
\label{chioc} \limsup_{n\rightarrow \infty}
f(\phi(s_n^{\prime},p_n)) < f(x).
\end{eqnarray}
By Lemma \ref{lmlm} (iv), there exists a closed
neighborhood $V\subset M$ of $y$ such that $p_n \notin V$, and 
choosing $t_n^{\prime}\in ]0,t_n[$ such that
$\phi(t_n^{\prime},p_n)\in \partial V$ (for $n$ large), we have
\begin{eqnarray}
\label{vic}
\set{\phi(t_n^{\prime},p_n)}{n\in \N} \cap A \mbox{ is pre-compact
  for every essentially vertical set }A\subset M\mbox{,}\\
\label{formi} \liminf_{n\rightarrow \infty}
f(\phi(t_n^{\prime},p_n)) > f(y).
\end{eqnarray}
The (PS) condition implies that $(t_n^{\prime}-s_n^{\prime})$ is
bounded: otherwise by Remark \ref{psr}, (\ref{chioc}) and
(\ref{formi}), we would obtain a sequence $(r_n)\subset \R$ such
that $(\phi(r_n,p_n))$ has a subsequence converging to a rest
point $z$, with $f(y)<f(z)<f(x)$, contradicting our previous
assumption. Therefore $t_n^{\prime}-s_n^{\prime}\leq T$ for every
$n\in \N$.

Since the essentially vertical family $\mathcal{F}$ is positively invariant, 
(\ref{stel}) implies that the set
\[
\set{\phi(t_n^{\prime},p_n)}{n\in \N} \subset \phi([0,T] \times
\set{\phi(s_n^{\prime},p_n)}{n\in \N})
\]
is essentially vertical. But then we can choose
$A=\set{\phi(t_n^{\prime},p_n)}{n\in \N}$ in (\ref{vic}), and we
obtain that the sequence $(\phi(t_n^{\prime},p_n))$ is compact. By the
boundeness of $t_n^{\prime}$ and by the fact that the vector field $F$
is complete, we conclude that also the sequence $(p_n)$ is compact.
\qed

\begin{rem}
An argument similar to the one used above shows that, if $F$
satisfies the assumptions of Theorem \ref{tcomp}, $x\in \rest(F)$
and $a\in \R$, then the set $W^u(x) \cap \{f\geq a\}$ is essentially
vertical. 
\end{rem}

\paragraph{Examples.}
Let us see what Theorem \ref{tcomp} says in the cases of finite
Morse indices or co-indices.

\begin{ex} {\em (Vector fields whose rest points have finite Morse
    index or co-index)}
Notice that the trivial subbundle $\mathcal{E}=(0)$ (relevant in
the case of rest points with finite Morse index, see Example
\ref{ex1}) has a strong integrable structure (choose any atlas of $M$).
The family consisting of all pre-compact subsets of
$M$ is a family of essentially vertical subsets for $\mathcal{E}=(0)$,
and it is obviously closed under the action of the flow.

Similarly, the trivial subbundle $\mathcal{E}=TM$ (relevant in
the case of rest points with finite Morse co-index) has a strong
integrable structure (again, consider an arbitrary atlas of $M$). 
The family consisting of all subsets of $M$ is a family of essentially
vertical subsets for $\mathcal{E}= TM$, clearly closed under the
action of the flow.

We have already seen that in the case $\mathcal{E}=(0)$ (resp.\
$\mathcal{E}=TM$) (C1) is
equivalent to the fact that all the rest points of $F$ have finite
Morse index (resp.\ co-index).

Therefore the
conclusion of Theorem \ref{tcomp} holds when (i) the $C^1$
vector field $F$ is Morse and gradient-like, (ii) either all the rest points
of $F$ have finite Morse index, or they have finite Morse co-index, 
(iii) $F$ satisfies (PS), and (iv) $F$ is complete.
\end{ex}

Now let us look back at Example \ref{ex2}.

\begin{ex}
{\em (Perturbations of a non-degenerate quadratic form)} Assume
that $M=H$ is a Hilbert space, and consider a function of the form
\[
f(\xi) = \frac{1}{2} \langle L\xi,\xi \rangle + b(\xi),
\]
where $L\in \mathcal{L}(H)$ is self-adjoint invertible, and 
the gradient of the function $b\in C^2(H)$ is a compact map. 
Let $F$ be the (negative) gradient vector field of $f$,
\[
F(\xi) = - \grad f(\xi) = - L \xi - \grad b(\xi).
\]
If $\grad b$ has linear growth - i.e. $\|\grad b(\xi)\| \leq
c(1+\|\xi\|)$ for every $\xi\in H$ - then $F$ is complete, its flow
$\phi$ maps bounded subsets of $\R\times H$ into bounded subsets of
$H$, and it satisfies 
\begin{equation}
\label{ide}
\phi(t,\xi) = e^{-tL} \xi - \int_0^t e^{(s-t)L} \grad b(\phi(s,\xi))\, ds.
\end{equation}
Consider the constant subbundle $V=V^-(L)$, and the orthogonal
projection $Q$ with kernel $V$. This bundle has the
trivial strong integrable structure modeled on $(H,V)$ 
consisting of the identity map: $\mathcal{A}=\{I\}$. 
The family $\mathcal{F}$ consisting of
all bounded subsets $A$ of $H$ such that $QA$ is pre-compact is an
essentially vertical family for $\mathcal{A}$. Moreover the identity
(\ref{ide}) together with the fact that $\grad b$ is a compact map
implies that $\mathcal{F}$ is invariant for $\phi$. 

The assumption that $\grad b$ has linear growth can be easily
dropped. Indeed, the vector field
\[
\tilde{F} (\xi) = - h(\xi) \grad f (\xi), \quad \mbox{where } h(\xi) =
\frac{1}{1+ \|\grad f(\xi)\|^2},
\]
is bounded, hence complete, and its flow $\tilde{\phi}$ maps bounded
subsets of $\R \times H$ into bounded subsets of $H$. Notice that $f$
is a non-degenerate Lyapunov function for $\tilde{F}$, and since $Df
[\tilde{F}] = - \|\grad f\|^2/ (1+\|\grad f\|^2)$, the Palais-Smale
sequences for $(\tilde{F},f)$ (in the sense of Definition \ref{psd})
are exactly the Palais-Smale sequences for $f$ (in the usual
sense). The flow $\tilde{\phi}$ satisfies
\[
\tilde{\phi}(t,\xi) = e^{-\tau(t,\xi)L} \xi - \int_0^t
h(\tilde{\phi}(t,\xi)) e^{(\tau(s,\xi)-\tau(t,\xi))L} \grad
b(\tilde{\phi}(s,\xi))\, ds,
\]
where $\tau:\R \times H\rightarrow \R$ is the function
\[
\tau(t,\xi) = \int_0^t h(\tilde{\phi}(s,\xi))\, ds.
\]
Then $|\tau(t,\xi)|\leq |t|$, and the fact that $\grad b$ is a compact 
map again implies that the family $\mathcal{F}$ is invariant for
$\tilde{\phi}$.    

We conclude that the thesis of the compactness Theorem \ref{tcomp}
holds, when $L$ is invertible and self-adjoint, $b\in C^2(H)$ has
compact gradient, and $f$ satisfies the Palais-Smale condition.
\end{ex}

In the case of a non-trivial subbundle $\mathcal{E}$ the question is
how to find an essentially vertical family which is
closed under the action of the flow of $F$. This question will be
addressed in the next section.  
  
\section{Flow-invariant essentially vertical families} 

\paragraph{Hausdorff measure of non-compactness.}
We recall that the Hausdorff distance of two subsets $A,B$ of a
metric space $X$ is the number
\[
\mathrm{dist}_{\mathcal{H}} (A,B) = \max \Bigl\{ \sup_{a\in A}
\inf_{b\in
  B} \dist (a,b),  \sup_{b\in B} \inf_{a\in A} \dist (a,b)\Bigr\}\in
  [0,+\infty].
\]
We denote by $\mathcal{H}(X)$ the family of all closed subsets of
$X$, and by $\mathcal{H}_b(X)$ the subfamily consisting of bounded
subsets, which is a metric space with the Hausdorff distance. A
related concept is the notion of measure of non-compactness. If
$A$ is a subset of a metric space $X$, its {\em Hausdorff measure
of non-compactness} is
\[
\beta_X(A) := \inf \set{r>0}{A \mbox{ can be covered by finitely
many
    balls of radius } r}\in [0,+\infty].
\]
Equivalently, $\beta_X(A)$ is the distance from the set of compact
subsets of $X$:
\begin{equation}
\label{hmnc} \beta_X(A) = \inf \set{\mathrm{dist}_{\mathcal{H}}
(A,K)}{K\subset X
  \mbox{ compact}}.
\end{equation}
It has the following properties (see \cite{dei85}, section 2.7.3):
\renewcommand{\theenumi}{\alph{enumi}}
\renewcommand{\labelenumi}{(\theenumi)}
\begin{enumerate}
\item $\beta_X(A)<+\infty$ if and only if $A$ is bounded; \item
$\beta_X(A)=0$ if and only if $A$ is totally bounded; \item if
$A_1 \subset A_2$ then $\beta_X(A_1) \leq \beta_X(A_2)$; \item
$\beta_X(A)\leq \beta_A(A) \leq 2 \beta_X(A)$; \item $\beta_X(A_1
\cup A_2) = \max \{ \beta_X(A_1),\beta_X(A_2) \}$; \item
$\beta_X(A) = \beta_X(\overline{A})$; \item $\beta_X$ is
continuous with respect to the Hausdorff distance.
\end{enumerate}
If $X$ is a normed vector space, denoting by $\co A$ the convex
hull of $A\subset X$, we also have:
\begin{enumerate}
\setcounter{enumi}{7} \item $\beta_X(\lambda A)=|\lambda|
\beta_X(A)$, and $\beta_X(A_1+A_2) \leq
  \beta_X(A_1) + \beta_X(A_2)$;
\setcounter{enumi}{9} \item $\beta_X(\co A) = \beta_X(A)$.
\end{enumerate}
\renewcommand{\theenumi}{\roman{enumi}}
\renewcommand{\labelenumi}{(\theenumi)}

\paragraph{Admissible presentations.}
We shall require that the essential subbundle $\mathcal{E}$ of $TM$ has a
strong presentation (see Definition \ref{spres}) which satisfies the 
following finiteness and a uniformity conditions.

\begin{defn}
\label{pres}
A strong presentation $\{M_i,N_i,\mathcal{Q}_i\}_{i\in I}$ is called
an admissible presentation if the Hilbert manifolds $N_i$ are endowed
with complete Riemannian metrics, and 
\begin{enumerate}
\item the covering $\{M_i\}_{i\in I}$ is star-finite (i.e. every $M_i$
  has non-empty intersection with finitely many $M_j$'s);

\item there is $r>0$ such that for every $p\in M$ there exists 
$i\in I$ such that
\[
\overline{\mathcal{Q}_i^{-1} (B_r(\mathcal{Q}_i(p)))} \subset M_i.
\]
\end{enumerate}
\end{defn}

An admissible presentation for $\mathcal{E}$ determines a strong
integrable structure $\mathcal{A}$ (see Proposition \ref{pesfe}).
Moreover, it determines a useful family of
essentially vertical subsets of $M$. Indeed, let $\mathcal{F}$ be the family
of subsets $A\subset M$ such that:
\begin{eqnarray}
\label{stla}
A \mbox{ can be covered by finitely many }M_i\mbox{'s;} \\
\label{crce}
\mbox{for every }i\in I, \;\mathcal{Q}_i(A\cap M_i) \mbox{ is pre-compact.}
\end{eqnarray}  
Proposition \ref{pesfe} implies that this is a family of essentially
vertical sets for the strong integrable structure $\mathcal{A}$.

Given an admissible presentation of $\mathcal{E}$ as above,
we shall assume the following condition on the vector field $F$:

\begin{description}
\item[(C3)] (i) there is $b>0$ such that $\|D\mathcal{Q}_i \circ
  F\|_{\infty}\leq b$ for every $i\in I$;

(ii) for every $i\in I$ and $q\in N_i$, there exists $\delta=\delta(q)>0$ and
  $c=c(q)\geq 0$ such that
\begin{equation}
\label{betalip} \beta_{TN_i} (D\mathcal{Q}_i(F(A))) \leq c \, \beta_{N_i}
(\mathcal{Q}_i(A)), \quad \forall A\subset \mathcal{Q}_i^{-1}(B_{\delta}(q)).
\end{equation}
\end{description}

Here the tangent bundle $TN_i$ is given the standard metric
induced by the Riemannian structure of $N_i$. Notice that no
Riemannian metric on $M$ is involved in this condition.

\begin{rem}
\label{rqbc}
If (C3)-(ii) holds, we can replace the point $q\in N_i$ by a compact
set $K\subset N_i$ in (\ref{betalip}): for every $i\in I$ and every
compact set $K\subset N_i$, there exists $\delta=\delta(K)>0$ and
$c=c(K)\geq 0$ such that
\[
\beta_{TN_i} (D\mathcal{Q}_i(F(A))) \leq c \, \beta_{N_i}
(\mathcal{Q}_i(A)), \quad \forall A\subset
\mathcal{Q}_i^{-1}(N_{\delta}(K)),
\]   
where $N_{\delta}(K)$ denotes the $\delta$-neighborhood of $K$.
\end{rem}

\begin{rem}
\label{Qcond} 
If $E$ is a Hilbert space and $\mathcal{Q}:M\rightarrow E$ is a $C^1$
map, one often makes no distinction between the tangential map
$D\mathcal{Q} : TM \rightarrow TE = E \times E$, $(p,\xi)\mapsto
(\mathcal{Q}(p), D\mathcal{Q}(p)[\xi])$, and its second component
$D\mathcal{Q} : TM \rightarrow E$, $(p,\xi)\mapsto D\mathcal{Q}
(p)[\xi]$. When $N_i=E$ is a Hilbert space, we are allowed
to replace the tangential map of $\mathcal{Q}_i$ by its second
component in (\ref{betalip}), writing $\beta_E(D\mathcal{Q}_i(F(A)))$
instead of $\beta_{E\times E}(D\mathcal{Q}_i (F(A)))$ on the left-hand
side of the inequality. Indeed, if $S \subset
TE=E\times E$, and $P_1,P_2:E\times E\rightarrow E$ are the
projections onto the first and the second factor, we have
\[
\max\{\beta_E(P_1 S),\beta_E(P_2 S)\} \leq \beta_{E\times E}(S)
\leq \beta_E(P_1 S) + \beta_E(P_2 S).
\]
\end{rem}

The main result of this section is the following proposition.

\begin{prop}
\label{clfam}
Let $\{M_i,N_i,\mathcal{Q}_i\}_{i\in I}$ be an admissible presentation for the
essential subbundle $\mathcal{E}$ of $TM$.
Assume that the vector field $F$ is complete and satisfies condition
(C3). Then the essentially vertical family $\mathcal{F}$ defined by
(\ref{stla}) and (\ref{crce}) is positively invariant for the flow of $F$.
\end{prop}

We start with the following local result.

\begin{lem}
\label{locale} Let $\mathcal{Q}:M\rightarrow E$ be a $C^1$ map into a
Hilbert space. Let $A\subset M$ be such that $\mathcal{Q}(A)$ is
pre-compact, and let $t^*\geq 0$ be such that $[0,t^*]\times A\subset
\Omega(F)$. Assume that there exists $c\geq 0$ such that
\[
\beta_E(D\mathcal{Q}(F(A^{\prime}))) \leq 
c\, \beta_E(\mathcal{Q}(A^{\prime})), \quad \forall A^{\prime}\subset
\phi([0,t^*]\times A).
\]
Then $\mathcal{Q}(\phi([0,t^*]\times A))$ is pre-compact.
\end{lem}

\proof Let $n=\lfloor ct^* \rfloor +1$, and set $\tau=t^*/n$, so
that $\tau c<1$. For $k\in \N$, $0\leq k\leq n$, set $A_k =
\phi([0,k\tau]\times A)$. Since
\[
\mathcal{Q}(\phi(t,p)) = \mathcal{Q}(p) + t \cdot \frac{1}{t} \int_0^t
D\mathcal{Q}(\phi(s,p))[F(\phi(s,p))]\, ds,
\]
we have
\[
\mathcal{Q}(A_{k+1}) = 
\mathcal{Q}(\phi([0,\tau]\times A_k) \subset \mathcal{Q}(A_k) + [0,\tau]
\overline{\mathrm{co}}\, (D\mathcal{Q}(F(A_{k+1}))).
\]
So, by the properties (c), (h), and (j) of the Hausdorff measure
of non-compactness, for $0\leq k\leq n-1$ we have,
\begin{eqnarray*}
\beta_E(\mathcal{Q}(A_{k+1})) \leq \beta_E(\mathcal{Q}(A_k)) + \tau \beta_E(
\overline{\mathrm{co}} \, (D\mathcal{Q}(F(A_{k+1})))) \\ \leq
\beta_E(\mathcal{Q}(A_k)) + \tau \beta_E(D\mathcal{Q}(F(A_{k+1}))) \leq
\beta_E(\mathcal{Q}(A_k)) + \tau c \beta_E(\mathcal{Q}(A_{k+1})).
\end{eqnarray*}
Since $\tau c<1$,
\[
\beta_E(\mathcal{Q}(A_{k+1})) \leq \frac{1}{1-\tau c} \beta_E(\mathcal{Q}(A_k)), \quad
k=0,1,\dots,n-1,
\]
and the fact that $\beta_E(\mathcal{Q}(A_0))=0$ implies that
$\beta_E(\phi([0,t^*]\times A))=\beta_E(\mathcal{Q}(A_n))=0$, as claimed.
\qed

\begin{ex} \label{exc}
The conclusion of the above lemma is not implied by the weaker
assumption that $\overline{\mathcal{Q}(F(S))}$ should be compact for every
set $S$ such that $\overline{\mathcal{Q}(S)}$ is compact, as the following
example shows.

Let $H=\ell_2(\Z)$, let $\set{e_k}{k\in \Z}$ be its standard
orthonormal basis, let $H^-=\overline{\Span} \set{e_k}{k\leq 0}$,
$H^+=\overline{\Span} \set{e_k}{k> 0}$, and let $Q$ be the
orthogonal projector onto $H^-$. Then there exists a smooth
bounded vector field $F:H\rightarrow H$ whose restriction to any
set of the form
\begin{equation}
\label{sof} \set{\xi \in H}{|\xi-\xi_0|<r} + H^-, \quad \xi_0\in
H^+, \; r<1,
\end{equation}
has finite rank, and whose flow $\phi$ has the property that
\[
\overline{Q\phi_1(\set{\xi\in H^+}{|\xi|\leq 1})}
\]
is not compact. In particular, $\overline{F(A)}$ is compact and
finite dimensional whenever $\beta_{H^-}(QA)<1$.

To construct such a vector field, for $k\in \N^*$ choose two
functions $f_k,g_k\in C^{\infty}(\R)$ such that
$f_k(s)=\sqrt{s+1/k}$ for $s\in [0,2]$, $\|f_k\|_{\infty}\leq 2$,
$g_k(1)=1$, $g_k(s)=0$ for $s\leq 1-1/k$, $0\leq g_k\leq 1$. Let
$\chi\in C^{\infty}(\R)$ be a function with compact support such
that $\chi(s)=1$ for $|s|\leq 2$, and set
\[
F(\xi) := \chi(|\xi|) \sum_{k=1}^{\infty} g_k(\xi \cdot e_{-k})
f_k(\xi \cdot e_k) e_k, \quad \xi\in H.
\]
The restriction of the vector field $F$ to a set of the kind
(\ref{sof}) has image contained in the finite dimensional
subspace $\Span\set{e_k}{k\in \N, \; \overline{\xi} \cdot e_k +
1/k >1-r}$. On the other hand, an easy computation shows that
\[
\phi(t,e_{-k}) = e_{-k} + \left(\frac{t^2}{4} + \frac{t}{\sqrt{k}}
\right) e_k, \quad k\geq 1, \; 0\leq t\leq 1,
\]
so the set $\set{Q \phi(t,e_{-k})}{k\geq 1}$ does not have compact
closure, for any $t\in [0,1]$.
\end{ex}

\begin{lem}
\label{aa1}
Let $\{M_i,N_i,\mathcal{Q}_i\}_{i\in I}$ be an admissible presentation of
the essential subbundle $\mathcal{E}$ of $TM$. Let $F$ be a $C^1$
complete vector field satisfying
\[
\|D\mathcal{Q}_i \circ F\|_{\infty} \leq b \quad \forall i\in I,
\]
for some $b\geq 0$. 
\begin{enumerate}
\item If $\overline{\mathcal{Q}_i^{-1}(B_r(\mathcal{Q}_i(p)))} \subset M_i$, then
  $\phi(s,p)\in M_i$ for every $|s|\leq r/b$, and
\begin{equation}
\label{ee1}
\dist (\mathcal{Q}_i(\phi(s,p)),\mathcal{Q}_i(p)) \leq b |s|.
\end{equation}

\item If a set $A\subset M$ can be covered by finitely many $M_i's$,
  then $A= \bigcup_{i\in I_0} A_i$, where 
\begin{equation}
\label{ee2}
A_i = \set{p\in A \cap M_i}{\overline{\mathcal{Q}_i^{-1}(B_r(\mathcal{Q}_i(p)))} \subset
  M_i},
\end{equation}
and $I_0 \subset I$ is finite.

\item If a set $A\subset M$ can be covered by finitely many $M_i$'s,
  then $\phi([0,t]\times A)$ can be covered by finitely many $M_i$'s,
  for every $t\geq 0$.
\end{enumerate}
\end{lem}

\proof
(i) Let $J$ be the maximal interval of numbers $s$ for which $\phi(s,p)\in
M_i$. Then
\begin{eqnarray*}
\dist(\mathcal{Q}_i(\phi(s,p)),\mathcal{Q}_i(p)) \leq \left| \int_0^s \Bigl|
  \frac{d}{d\sigma} \mathcal{Q}_i(\phi(\sigma,p)) \Bigr| \, d\sigma \right| \\
= \left| \int_0^s | D\mathcal{Q}_i\circ F(\phi(\sigma,p)) |\, d\sigma \right|
  \leq b |s| \quad \forall s\in J.
\end{eqnarray*}
Together with the fact that the closure of $\mathcal{Q}_i^{-1}(B_r(\mathcal{Q}_i(p)))$ is
contained in $M_i$, this implies that $]-r/b,r/b[\subset J$ and
(\ref{ee1}).

\medskip

(ii) Since $A$ is covered by finitely many $M_i$'s and the covering
$\{M_i\}_{i\in I}$ is star-finite, the indices $i\in I$ for which
$A_i\neq \emptyset$ form a finite subset $I_0$. By the uniformity
property of the presentation (Definition \ref{pres}, (iv)),
$A=\bigcup_{i\in I_0} A_i$. 

\medskip

(iii) If $A_i$ are the sets defined in (\ref{ee2}), statement (i)
implies that $\phi([0,r/b[\times A_i)\subset M_i$ for every $i\in
I_0$. Therefore $\phi([0,r/b[\times A)$ is covered by the finite
covering $\{M_i\}_{i\in I_0}$, and the conclusion follows by
induction.
\qed

\medskip

{\em Proof} [of Proposition \ref{clfam}].
By Lemma \ref{aa1} (iii), $\phi(0,t]\times A$ is covered by finitely
many $M_i's$, so it is enough to show that the interval
\[
\mathcal{T}(A) = \set{ t\geq 0}{\mathcal{Q}_i(\phi([0,t]\times A)) \mbox{ is
    pre-compact in } N_i, \;\forall i\in I}
\]
coincides with $[0,+\infty[$. Since $0\in \mathcal{T}(A)$, we can
argue by connectedness proving that $\mathcal{T}(A)$ is both open and
closed in $[0,+\infty[$.

\medskip

We claim that $\mathcal{T}(A)$ is open in $[0,+\infty[$. Let $t\in
\mathcal{T}(A)$. By Lemma \ref{aa1} (ii), $\phi([0,t]\times A) =
\bigcup_{i\in I_0} A_i$, where 
\begin{equation}
\label{diciasete}
A_i = \set{p\in \phi([0,t]\times A) \cap M_i}{\overline{
    \mathcal{Q}_i^{-1}(B_r(\mathcal{Q}_i(p)))} \subset M_i},
\end{equation}
and $I_0\subset I$ is finite.
Clearly, $\mathcal{T}(A)=[0,t] + \bigcap_{i\in I_0} \mathcal{T}(A_i)$,
so it is enough to prove that $\sup \mathcal{T}(A_i)>0$, for every
$i\in I_0$. 

Let $i\in I_0$. Since $\mathcal{Q}_i(A_i)$ is pre-compact, (C3)-(ii),
together with Remark \ref{rqbc}, implies
that there exist $c\geq 0$ and $\delta>0$ such that
\begin{equation}
\label{pesce}
\beta_{TN_i} (D\mathcal{Q}_i \circ F(S)) \leq c \beta_{N_i} (\mathcal{Q}_i(S)) \quad
\forall S \subset \mathcal{Q}_i^{-1}(N_{\delta}(\mathcal{Q}_i(A_i))).
\end{equation}
Moreover, $\overline{\mathcal{Q}_i(A_i)}$ is covered by finitely many 
coordinate neighborhoods: there exist $q_1,\dots,q_n\in 
\mathcal{Q}_i(A_i)$, $0<\rho\leq
\min \{\delta,r\}$, $\mathcal{Q}_i(A_i) \subset \bigcup_{j=1}^n B_{\rho/2}
(q_j)$, and local charts
\[
\psi_j : \dom (\psi_j) \rightarrow E,
\]
with $\overline{B_{\rho}(q_j)}\subset \dom (\psi_j)$,
$\overline{\psi_j(B_{\rho}(q_j))} \subset \dom (\psi_j^{-1})$, and
$\psi_j^{-1}$, $D\psi_j$ Lipschitz. Then $A_i = \bigcup_{j=1}^n
A_i^j$, with $A_i^j = A_i \cap \mathcal{Q}_i^{-1} (B_{\rho/2}
(q_j))$. 
Again, it
suffices to show that $\sup \mathcal{T}(A_i^j)>0$.

Let $1\leq j\leq n$, and set $U=\mathcal{Q}_i^{-1}(B_{\rho}(q_j))$. Since
$\rho\leq r$ and $q_j\in \mathcal{Q}_i(A_i)$, by the definition of $A_i$ we have
\begin{equation}
\label{stmar}
\overline{U} \subset \overline{ \mathcal{Q}_i^{-1} (B_r(q_j))} \subset M_i.
\end{equation}
Let $p\in A_i^j\subset M_i$. Let $[0,\tau(p)[$, $0<\tau(p)\leq
+\infty$, be the maximal interval of positive numbers $s$ for which
$\phi(s,p)\in U$. By Lemma \ref{aa1} (i),
\[
\dist (\mathcal{Q}_i(\phi(s,p)),q_j) \leq \dist (\mathcal{Q}_i(\phi(s,p)),\mathcal{Q}_i(p)) + \dist
(\mathcal{Q}_i(p),q_j) \leq bs + \rho/2
\]
for every $s\in [0,\tau(p)[$. Together with (\ref{stmar}) this implies
that $\tau(p)\geq \rho/(2b)$. Therefore
\begin{equation}
\label{triang}
\phi([0,\rho/(2b)[\times A_i^j) \subset U.
\end{equation}

Let $\mathcal{Q}:= \psi_j \circ \mathcal{Q}_i:U\rightarrow E$. Since $\rho\leq \delta$
and $q_j\in \mathcal{Q}_i(A_i)$, $U$ is contained in
$\mathcal{Q}_i^{-1}(N_{\delta}(\mathcal{Q}_i(A_i)))$. By (\ref{pesce}), for any $S\subset
U$,
\begin{eqnarray*}
\beta_E((D\mathcal{Q}_i\circ F)(S)) \leq \beta_{E\times E}((D\mathcal{Q}\circ F)(S)) =
\beta_{E\times E} (D\psi_j \circ D\mathcal{Q}_i(F(S))) \\ \leq \mathrm{lip}
(D\psi_j) \beta_{TN_i} (D\mathcal{Q}_i \circ F(S)) \leq c\, \mathrm{lip} (D\psi_j)
\beta_{N_i} (\mathcal{Q}_i(S)) \\ = c\, \mathrm{lip} (D\psi_j) \beta_{N_1}
(\psi_j^{-1}(\mathcal{Q}(S))) \leq c \, \mathrm{lip}(D\psi_j) \, \mathrm{lip}
(\psi_j^{-1}) \beta_E(\mathcal{Q}(S)).
\end{eqnarray*}
By (\ref{triang}), we can take $S=\phi([0,\rho/(2b)[\times A_i^j)$ in
the above inequality, and Lemma \ref{locale} implies that
$\mathcal{Q}(\phi([0,\rho/2b[\times A_i^j))$ is pre-compact in $E$.

Since $\overline{B_{\rho}(q_j)} \subset \dom (\psi_j)$ and
$\overline{\psi_j(B_{\rho}(q_j))} \subset \dom (\psi_j^{-1})$, the set
\[
\overline{\mathcal{Q}_i(\phi([0,\rho/(2b)[\times A_i^j))} = \psi_j^{-1} (
\overline{\mathcal{Q}(\phi([0,\rho/(2b)[\times A_i^j)))}
\]
is a compact subset of $N_i$. Therefore $\sup \mathcal{T}(A_i^j) \geq
\rho/(2b) >0$, as we wished to prove.

\medskip

There remains to show that the interval $\mathcal{T}(A)$ is
closed. Let $t= \sup \mathcal{T}(A)$. By Lemma \ref{aa1} (ii),
$\phi([0,t] \times A) = \bigcup_{i\in I_0} A_i$, where $A_i$ is
defined in (\ref{diciasete}) and $I_0\subset I$ is finite. It is
enough to prove that $\mathcal{Q}_i(A_i)$ has compact closure in $N_i$, for any
$i\in I_0$.

Fix $i\in I_0$, and let $q_k=\mathcal{Q}_i(\phi(t,p_k))$, where $p_k\in A$ and
$\phi(t,p_k)\in A_i$, be a sequence in $\mathcal{Q}_i(A_i)$. By Lemma \ref{aa1}
(i), $\phi(t-\tau,p_k)\in M_i$ for any $0\leq \tau<r/b$, and
\begin{equation}
\label{ladist}
\dist(\mathcal{Q}_i(\phi(t-\tau,p_k)),q_k) \leq b \tau.
\end{equation}
Since $t=\sup \mathcal{T}(A)$, the sequence
$(\mathcal{Q}_i(\phi(t-\tau,p_k)))_{k\in \N}$ is compact for any $0<\tau
<r/b$. Then (\ref{ladist}) and the completeness of $N_i$ imply that
also the sequence $(q_k)$ is compact, proving that $\mathcal{Q}_i(A_i)$ is
pre-compact.
\qed 

\paragraph{Properties of condition (C3).}
Condition (C3) is stronger than (C2), and like (C2) it is a convex
condition. Indeed, the following result holds.

\begin{prop}
Let $\{M_i,N_i,\mathcal{Q}_i\}_{i\in I}$ be an admissible
presentation of the essential subbundle $\mathcal{E}$. The
following facts hold:
\begin{enumerate}
\item condition (C3) implies condition (C2); 
\item the vector fields $F$ satisfying (C3) form a module over the
  ring $C^1(M)\cap C^0_b(M)$
\end{enumerate}
\end{prop}

\proof (i) Let $i\in I$, $p\in M_i$, and set $q:=\mathcal{Q}_i(p)$.
Up to the composition with $C^1$ local charts
\[
\varphi:\dom (\varphi)\subset M_i \rightarrow H, \; p\in \dom (\varphi),
\quad \psi: \dom (\psi)\subset N_i \rightarrow
E, \; q\in \dom (\psi),
\]
such that $\psi$, and $D\psi$ are bi-Lipschitz, we may
assume that $\mathcal{Q}_i$ is a $C^1$ semi-Fredholm map with 
$\ind \mathcal{Q}_i\geq 0$ from an open set of the Hilbert space $H$ 
into the Hilbert space $E$. By (C3)-(ii) and Remark \ref{Qcond}, there 
exist $\delta>0$ and $c\geq 0$ such that
\begin{equation}
\label{nc} \beta_E(D\mathcal{Q}_i(F(A))) \leq c\beta_E(\mathcal{Q}_i(A)) 
\quad \forall A\subset \mathcal{Q}^{-1}_i(B_{\delta}(q)).
\end{equation}
Let $T\in \mathcal{L}(H,E)$ be a linear map with finite rank such
that $D\mathcal{Q}_i(p)+T$ is surjective. Since $T$ has finite rank,
\begin{equation}
\label{ac} \beta_E(\mathcal{Q}_i(A)) = \beta_E((\mathcal{Q}_i+T)(A)), \quad
\beta_E(D\mathcal{Q}_i(F(A))) = \beta_E(D(\mathcal{Q}_i+T)(F(A))).
\end{equation}
The map $\mathcal{Q}_i+T$ is a local submersion at $p$, so up to considering a
change of variable at $p$, we may assume that the restriction of
$\mathcal{Q}_i$ to a neighborhood $U$ of $p$ coincides with a linear
surjective map $Q$ from $H$ to $E$, which by (\ref{nc}) and (\ref{ac})
verifies
\begin{equation}
\label{nncc} \beta_E(QF(A)) \leq c \beta_E(QA) \quad \forall
A\subset Q^{-1}(B_{\delta}(q))\cap U.
\end{equation}
By composing with a linear right inverse of $Q$, we may also
assume that $E$ is a closed subspace of $H$ and that $Q$ is a
linear projector onto $E$. In these coordinates, the essential
subbundle $\mathcal{E}$ is locally represented by the constant
subbundle $\ker Q$, with projector $P=I-Q$. By (\ref{nncc}), the
map $(I-P)FP$ is compact in a neighborhood of $p$, so its 
differential at $p$,
\[
D((I-P)FP)(p) = (I-P) DF(p) P = (L_F P)(p) P,
\]
is a compact operator, proving (C2).

\medskip

(ii) Let $F_1$ and $F_2$ be $C^1$ tangent vector fields on $M$,
and let $h_1,h_2 \in C^1(M) \cap C^0_b(M)$. Let $i\in I$. Clearly,
if $F_1$ and $F_2$ satisfy (C3)-(i) with constants $b_1$ and $b_2$,
$h_1 F_1 + h_2 F_2$ satisfy (C3)-(i) with constant $\|h_1\|_{\infty}
b_1 + \|h_2\|_{\infty} b_2$.

Let $p\in M_i$, and set $q:=\mathcal{Q}_i(p)$. Let $\psi:U \rightarrow
E$, $q\in U\subset N_i$, be a local chart such that 
\[
D \psi : TU \rightarrow \psi(U) \times E \subset E \times
E
\]
is bi-Lipschitz of constant 2. By property (h) of the Hausdorff
measure of non-compactness, if $A\subset\mathcal{Q}_i^{-1}(U)$,
\begin{eqnarray*}
\beta_{TN_i} (D\mathcal{Q}_i \circ (h_1 F_1+h_2 F_2)
  (A)) \leq 2 \beta_{E \times E} ( D \psi \circ D\mathcal{Q}_i 
\circ (h_1 F_1+h_2 F_2) (A))
\\ \leq 2 \|h_1\|_{\infty} \beta_{E \times E} ( D\psi \circ D\mathcal{Q}_i
\circ F_1 (A)) + 2 \|h_2\|_{\infty} \beta_{E \times E} ( D\psi
\circ D\mathcal{Q}_i \circ F_2 (A)) \\ \leq 4 \|h_1\|_{\infty} \beta_{TN_i}
(D\mathcal{Q}_i \circ F_1(A)) + 4 \|h_2\|_{\infty}  \beta_{TN_i} 
(D\mathcal{Q}_i \circ F_2(A)).
\end{eqnarray*}
Therefore, if $F_1$ and $F_2$ satisfy (\ref{betalip}) with
constants $\delta_1,c_1$ and $\delta_2,c_2$, then $h_1 F_1+h_2 F_2$
satisfies (\ref{betalip}) with constants $\delta =
\min\{\delta_1,\delta_2\}$ and $c=4\|h_1\|_{\infty} c_1+4\|h_2\|_{\infty}
c_2$. This proves that $h_1 F_1 + h_2 F_2$ satisfies (C2)-(ii). 
\qed

It seems useful to find sufficient conditions implying (C3)-(ii), which do
not make use of the measures of non-compactness but are stated only in
terms of Lipschitzianity and compactness of some maps. 
To this purpose, assume that $M$ is endowed with a Riemannian
metric.   

The following proposition says that under mild
Lipschitz assumptions on $\mathcal{Q}_i^{-1}$ and $F$, condition
(C3)-(ii) holds if and only if the maps 
$D\mathcal{Q}_i \circ F :\mathcal{Q}_i^{-1}(\{q\})
\rightarrow T_q N_i$ have pre-compact image. Example \ref{exc} 
suggests that without Lipschitz assumptions on $F$ the last
condition is not sufficient for the conclusion of Proposition
\ref{clfam} to hold.

\begin{prop}
\label{dati}
Assume that every map $N_i \rightarrow \mathcal{H}(M_i)$,
$q\mapsto \mathcal{Q}_i^{-1}(\{q\})$, is locally Lipschitz, and that for every
$q\in N_i$ there exists $\delta>0$ such that the map
$D\mathcal{Q}_i \circ F$ is Lipschitz on
$\mathcal{Q}_i^{-1}(B_{\delta}(q))$. Then (C3)-(ii) holds if and only if
$D\mathcal{Q}_i \circ F$ maps every fiber $\mathcal{Q}_i^{-1}(\{q\})$
into a pre-compact set.
\end{prop}

\begin{rem}
The assumption on the local Lipschitzianity of $\mathcal{Q}_i^{-1}$,
required in the above proposition involves a uniform lower bound 
on the non-zero singular values of $D\mathcal{Q}_i$. 
More precisely, let $\mathcal{Q}:M\rightarrow N$
be a $C^1$ submersion between Riemannian Hilbert manifolds, with
$M$ complete and $N$ connected. If there is $\alpha>0$ such that
\begin{equation}
\label{singval} \inf \left( \sigma(D\mathcal{Q}(p)^* D\mathcal{Q}(p)) 
\setminus \{0\} \right)\geq \alpha \quad \forall p\in M,
\end{equation}
then the map $N\rightarrow \mathcal{H}(M)$, $q\mapsto
\mathcal{Q}^{-1}(\{q\})$, is $1/\sqrt{\alpha}$-Lipschitz.
\end{rem}

In order to prove this statement, let $q_0,q_1\in N$ and
$k>1/\sqrt{\alpha}$. Let $p_0\in \mathcal{Q}^{-1}(\{q_0\})$ and let
$v:[0,1]\rightarrow N$ be a $C^1$ curve such that $v(0)=q_0$,
$v(1)=q_1$. Since $\mathcal{Q}$ is a submersion, for any $t_0\in [0,1]$ and
every $p\in \mathcal{Q}^{-1}(\{v(t_0)\})$, there exists a $C^1$ local
lifting $u$ of $v$ verifying $u(t_0)=p$ and $u^{\prime}(t_0)\in
(\ker T_p \mathcal{Q})^{\perp}$. Assumption (\ref{singval}) easily implies
that $|u^{\prime}(t_0)|\leq (1/\sqrt{\alpha}) |v^{\prime}(t_0)|$,
so
\begin{equation}
\label{llipp} |u^{\prime}(t)| \leq k |v^{\prime}(t)|,
\end{equation}
for any $t$ in a neighborhood of $t_0$. By a standard maximality
argument, it follows that there exists a Lipschitz lifting
$u:[0,1]\rightarrow M$ of $v$ with $u(0)=p_0$ and verifying
(\ref{llipp}) a.e.\ in $[0,1]$. Therefore,
\[
\dist (p_0,\mathcal{Q}^{-1}(\{q_1\})) \leq \int_0^1 |u^{\prime}(t)|\, dt
\leq k \int_0^1 |v^{\prime}(t)|\, dt.
\]
Hence taking the infimum over $k$ and $v$ we obtain $\dist (p_0,
\mathcal{Q}^{-1}(\{q_1\})) \leq 1/\sqrt{\alpha} \, d(q_0,q_1)$, and by
symmetry $\mathrm{dist}_{\mathcal{H}}(\mathcal{Q}^{-1}(\{q_0\}),
\mathcal{Q}^{-1}(\{q_1\})) \leq 1/\sqrt{\alpha} \, d(q_0,q_1)$, as claimed.

\medskip

\proof (of Proposition \ref{dati}) 
Let $i\in I$. Since $\beta_{N_i}(\{q\})=0$,
condition (C3)-(ii) trivially implies that $D\mathcal{Q}_i
\circ F$ maps every fiber $\mathcal{Q}_i^{-1}(\{q\})$ into a
pre-compact set, for every $q\in N_i$.

Let us prove the converse statement. Let $q\in N_i$ and let $\delta>0$
be so small that the maps
\begin{eqnarray*}
B_{\delta}(q) \rightarrow
\mathcal{H}(\mathcal{Q}_i^{-1}(B_{\delta}(q)), 
\quad q^{\prime} \mapsto \mathcal{Q}_i^{-1}( \{q^{\prime}\}), \\
\mathcal{Q}_i^{-1}(B_{\delta}(q)) \rightarrow TN_i, \quad p \mapsto
D\mathcal{Q}_i \circ F(p),
\end{eqnarray*}
are Lipschitz. Then also the maps
\begin{eqnarray*}
\mathcal{H}(B_{\delta}(q)) \rightarrow
\mathcal{H}(\mathcal{Q}_i^{-1}(B_{\delta}(q)), 
\quad \Sigma \mapsto \mathcal{Q}_i^{-1}(\Sigma), \\
\mathcal{H}(\mathcal{Q}_i^{-1}(B_{\delta}(q))) \rightarrow 
\mathcal{H}(TN_i), \quad A \mapsto
\overline{D\mathcal{Q}_i \circ F(A)},
\end{eqnarray*}
are Lipschitz. Let $c$ be the Lipschitz constant of their composition
\[
\mathcal{H}(B_{\delta}(q))
\stackrel{\mathcal{Q}_i^{-1}}{\longrightarrow} \mathcal{H}
(\mathcal{Q}_i^{-1}(B_{\delta}(q)) \stackrel{\overline{D\mathcal{Q}_i
    \circ F}}{\longrightarrow} \mathcal{H}(TN_i).
\]
Let $A\subset \mathcal{Q}_i^{-1}(B_{\delta}(q))$ be the set for which
we wish to prove (\ref{betalip}). We may assume that
$A=\mathcal{Q}_i^{-1}(\Sigma)$ for some closed subset $\Sigma\subset
B_{\delta}(q)$. 
If $\Sigma_0\subset B_{\delta}(q)$ is a finite set, our assumption
implies that $D\mathcal{Q}_i \circ F(\mathcal{Q}_i^{-1}(\Sigma_0))$ is 
pre-compact. So by (\ref{hmnc}),
\begin{eqnarray*}
\beta_{TN_i} (D\mathcal{Q}_i(F(A))) =
\beta_{TN_i} (\overline{D\mathcal{Q}_i(F(
\mathcal{Q}_i^{-1} (\Sigma)))}) \\ \leq \mathrm{dist}_{\mathcal{H}}
(\overline{D\mathcal{Q}_i \circ F} \circ \mathcal{Q}_i^{-1} (\Sigma), 
\overline{D\mathcal{Q}_i\circ
  F} \circ \mathcal{Q}_i^{-1} (\Sigma_0)) 
\leq c \, \mathrm{dist}_{\mathcal{H}}
(\Sigma,\Sigma_0).
\end{eqnarray*}
By the density of the space of finite sets in the space of compact
sets, by (\ref{hmnc}), and by property (d) of $\beta$, we obtain
\begin{eqnarray*}
\beta_{TN_i}(D\mathcal{Q}_i(F(A))) \leq c\, \inf_{\substack{\Sigma_0 
\subset B_{\delta}(q))\\
    \Sigma_0 \;\mathrm{finite}}} \mathrm{dist}_{\mathcal{H}}
    (\Sigma,\Sigma_0) = c\, \inf_{\substack{\Sigma_0 \subset B_{\delta}(q)\\
    \Sigma_0 \;\mathrm{compact}}} \mathrm{dist}_{\mathcal{H}}
    (\Sigma,\Sigma_0) \\ = c \, \beta_{B_{\delta}(q)} 
    (\Sigma) \leq 2c\, \beta_{N_i} (\Sigma) =
    2c \, \beta_{N_i}(\mathcal{Q}_i(A)),
\end{eqnarray*}
which proves (C3)-(ii). \qed

\begin{ex}
{\em (Product manifolds)} Let us consider again the situation of
Example \ref{ex3}: $M=M^- \times M^+$ is given a product complete
Riemannian structure, and $\mathcal{V}=TM^- \times (0)$. Consider
the projection onto the second factor $\mathcal{Q}:M\rightarrow M^+$,
$(p^-,p^+)\mapsto p^+$, and the admissible presentation of the
subbundle $\mathcal{V}$, 
\[
\set{\mathcal{Q}|_{B^-\times B^+}}{B^- \times
    B^+\subset M^- \times M^+
\mbox{ is bounded}}.
\] 
Writing the tangent vector field $F$ as
$F(p^-,p^+)=(F^-(p^-,p^+),F^+(p^-,p^+)) \in T_{p^-} M^- \times
T_{p^+} M^+$, there holds $D\mathcal{Q}\circ F (p^-,p^+) = F^+(p^-,p^+)$.
Assume that (i) $F^+$ is bounded, and (ii)
for every $p^+\in M^+$ and for every bounded set
$B^-\subset M^-$ there exists $\delta>0$ such that $F^+$ is Lipschitz on 
$B^- \times B_{\delta}(p^+)$, and that the map $M^-\rightarrow T_{p^+}
M^+$, $p^-\mapsto F^+(p^-,p^+)$ is compact. Then Proposition
\ref{dati} implies that $F$ satisfies (C3).
\end{ex}

\section{Broken flow lines}

Let $x$ and $y$ be rest points of the gradient-like Morse vector
field $F$. Let us assume that $W^u(x)\cap W^s(y)$ has compact
closure. Consider a sequence of flow lines from the rest point $x$
to the rest point $y$, and the sequence of their closures:
\[
S_n= \overline{\phi(\R \times \{p_n\})} = \phi(\R\times
\{p_n\})\cup \{x,y\}, \quad p_n\in W^u(x)\cap W^s(y).
\]
Since $\overline{W^u(x)\cap W^s(y)}$ is compact, up to a
subsequence we may assume that $p_n\rightarrow p$, and the
continuity of $\phi$ would give us the convergence
\[
\phi(\cdot,p_n) \rightarrow \phi(\cdot,p)
\]
uniformly on compact subsets of $\R$. However, it may happen that
$p\notin  W^u(x)$, or $p\notin W^s(y)$, so $\phi(\cdot,p)$ could
be a flow line connecting two other rest points, and the
convergence would not be uniform on $\R$. We will show that in
this case a subsequence of $(S_n)$ converges to a broken flow line
in the Hausdorff distance. The discussion is independent on
conditions (C1-3), and involves only the compactness of
$\overline{W^u(x)\cap W^s(y)}$.

\begin{defn}
\label{bgfl} Let $x,y\in \rest(F)$. A {\em broken flow line} from
$x$ to a $y$ is a set
\[
S = S_1 \cup \dots \cup S_k,
\]
where $k\geq 1$, $S_i$ is the closure of a flow line from
$z_{i-1}$ to $z_i$, where  $x=z_0\neq z_1\neq \dots \neq z_{k-1}
\neq z_k = y$ are rest points.
\end{defn}

When $k=1$, a broken flow line is just the closure of a genuine
flow line. Let us fix a Lyapunov function $f$ for $F$. If $S$ is
a broken flow line as in the above definition, the following
inequalities must hold:
\begin{equation}
\label{livfun} f(x) > f(z_1) > \dots > f(z_{k-1}) > f(y).
\end{equation}
It is easy to check that a compact set $S\subset M$ is a broken
flow line from $x$ to $y$ if and only if (i) $x,y\in S$, (ii) $S$
is $\phi$-invariant, (iii) the intersection
\[
S \cap \set{p\in M}{f(p)=c}
\]
consists of a single point if $c\in [f(y),f(x)]$, and it is empty
otherwise. Now we can state the compactness result for the
gradient flow lines.

\begin{prop}
\label{comp} Assume that the Morse vector field $F$ has a
Lyapunov function $f$, and that $x,y$ are rest points such that
$\overline{W^u(x) \cap W^s(y)}$ is compact. Let $(p_n)\subset
W^u(x)\cap W^s(y)$, and set $S_n := \phi(\R \times \{p_n\}) \cup
\{x,y\}$. Then $(S_n)$ has a subsequence which converges to a
broken flow line from $x$ to $y$, in the Hausdorff distance.
\end{prop}

\proof The space of compact subsets of a compact metric space is
compact with respect to the Hausdorff distance, so $(S_n)$ has a
subsequence $(S_n^{\prime})$ which converges to a compact set
$S\subset \overline{W^u(x)\cap W^s(y)}$. Then $x,y\in S$, and
since $S_n^{\prime}$ is $\phi$-invariant, so is $S$. Since
$S_n^{\prime} \subset f^{-1} ([f(y),f(x)])$, we obtain that the
set
\begin{equation}
\label{juve} S \cap \set{p\in M}{f(p)=c}
\end{equation}
is empty for every $c\notin [f(y),f(x)]$.

Let $c\in [f(y),f(x)]$. Then $(S_n^{\prime})$ has a subsequence
$(S_n^{\prime\prime})$ such that
\[
S_n^{\prime\prime} \cap \set{p\in M}{f(p)=c}
\]
converges to some point in (\ref{juve}), which is then non-empty.
If the set (\ref{juve}) contains two points $p,q$, the fact that
$(S_n)$ is a sequence of flow lines allows us to find a sequence
$(p_n)$ converging to $p$, and numbers $t_n\in \R$ such that
$\phi(t_n,p_n)\rightarrow q$. By reversing if necessary the role
of $p$ and $q$, we may assume that $t_n\geq 0$ for every $n$, and
we deduce the convergence:
\[
\int_0^{t_n} Df(\phi(t,p_n))[F(\phi(t,p_n))] \, dt =
f(\phi(t_n,p_n)) - f(p_n) \rightarrow f(q) - f(p) =0.
\]
Then the fact that the rest points of $f$ are isolated easily
implies that either $t_n \rightarrow 0$, or the sequence of sets
$\phi([0,t_n]\times \{p_n\})$ converges to a rest point. In both
cases, we obtain that $p=q$. This shows that the set (\ref{juve})
consists of a single point. Hence $S$ is a broken flow line from
$x$ to $y$. \qed

\section{Intersections of dimension 1 and 2}

Assume that the gradient-like Morse vector field  $F$ satisfies
(C1-2) with respect to a $(0)$-essential subbundle $\mathcal{E}$
of $TM$, so that the relative index $m(x,\mathcal{E})$ is a
well-defined integer, for any $x\in \rest(F)$.

We say that $F$ satisfies the {\em Morse-Smale property up to
order} $k$, if $W^u(x)$ and $W^s(y)$ have transverse intersection
for every pair of rest points $x,y$ such that
$m(x,\mathcal{E})-m(y,\mathcal{E})\leq k$. The Morse-Smale
condition up to order 0 implies that, for a broken flow line as in
Definition \ref{bgfl},
\begin{equation}
\label{livind} m (x,\mathcal{E}) > m (z_1,\mathcal{E}) > \dots
> m (z_{k-1},\mathcal{E})
> m (y,\mathcal{E}).
\end{equation}
In this section we shall assume that $F$ has the Morse-Smale
property up to order 2, and we shall describe the intersections
$W^u(x)\cap W^s(y)$ when $m(x,\mathcal{E})-m(y,\mathcal{E})$
is either 1 or 2. As in the last section, we shall assume that
such an intersection has compact closure. By Theorem \ref{ubu},
$W^u(x)\cap W^s(y)$ is a submanifold of dimension 1, respectively
2. The flow $\phi$ defines a free action of the group $\R$ onto
$W^u(x)\cap W^s(y)$, so the quotient, that is the set of the flow
lines from $x$ to $y$, is a manifold of dimension 0, respectively
1.

Assume that $x,y$ are rest points with
\begin{equation}
\label{dif1} m (x,\mathcal{E}) - m (y,\mathcal{E}) = 1.
\end{equation}
We claim that $W^u(x)\cap W^s(y)$ consists of {\em finitely many}
connected components. Indeed each connected component is a flow
line from $x$ to $y$, and the set $C$ of their closures is
discrete in the Hausdorff distance. On the other hand,
(\ref{livind}) and (\ref{dif1}) imply that these are the only
broken flow lines from $x$ to $y$. So by Proposition \ref{comp},
$C$ is also compact, hence finite.

Note that the restriction of the flow $\phi$ to the closure of a
component of $W^u(x)\cap W^s(y)$ is conjugated to the shift flow
on $\overline{\R}=[-\infty,+\infty]$:
\[
\R\times \overline{\R} \ni (t,u) \mapsto u+t \in \overline{\R}.
\]

Now assume that $x,z$ are rest points with
\begin{equation}
\label{dif2} m (x,\mathcal{E}) - m (z,\mathcal{E}) = 2.
\end{equation}
The quotient of each connected component $W$ by this action,
$W/\R$, being a connected one-dimensional manifold, is either the
circle or the open interval. In other words a connected component
$W$ of $W^u(x)\cap W^s(z)$ is described by a one-parameter family
of flow lines $u_{\lambda}$, where $\lambda$ ranges in $S^1$ or in
$]0,1[$.

In the first case one can easily verify that $\overline{W}=W\cup
\{x,z\}$ is homeomorphic to a 2-sphere, and that the restriction
of $\phi$ to $\overline{W}$ is conjugated to the exponential flow
on the Riemann sphere $S^2=\C\cup \{\infty\}$:
\[
\R\times S^2 \ni (t,\zeta) \mapsto e^t \zeta \in S^2.
\]
In the second case, by Proposition \ref{comp},
$\overline{W}\setminus W$ contains broken flow lines, which have
just one intermediate rest point, by (\ref{livind}) and
(\ref{dif2}). Then the flow $\phi$ restricted to $\overline{W}$ is
semi-conjugated to the product of two shift-flows on
$\overline{\R}$. More precisely, the situation is described by the
following theorem.

\begin{thm}
\label{conju} Assume that the gradient-like Morse vector field
$F$ satisfies (C1-2) with respect to a $(0)$-essential subbundle
$\mathcal{E}$ of $TM$. Assume that $F$ has the Morse-Smale
property up to order 2. Let $x,z$ be rest points such that
$m(x,\mathcal{E})-m(z,\mathcal{E})=2$, and let $W$ be a
connected component of $W^u(x) \cap W^s(z)$ such that
$\overline{W}$ is compact, and $W/\R$ is an open interval. Then
there exists a continuous surjective map
\[
h: \overline{\R} \times \overline{\R} \rightarrow \overline{W}
\]
with the following properties:
\begin{enumerate}
\item $\phi_t(h(u,v))=h(u+t,v+t)$, for every $(u,v)\in
  \overline{\R}\times \overline{\R}$, $t\in \R$;
\item $h(\R^2)=W$, and there exist rest points $y,y^{\prime}$ with
$m(y,\mathcal{E})=m(y^{\prime},\mathcal{E})=m(x,\mathcal{E})-1$,
and $W_1$, $W_2$, $W_1^{\prime}$, $W_2^{\prime}$ connected
components of $W^u(x)\cap W^s(y)$, $W^u(y)\cap W^s(z)$,
$W^u(x)\cap W^s(y^{\prime})$, $W^u(y^{\prime})\cap W^s(z)$,
respectively, such that $W_1\cup W_2 \neq W_1^{\prime} \cup
W_2^{\prime}$, and
\[
h(\R \times \{-\infty\})=W_1, \quad h(\{+\infty\} \times
    \R)=W_2, \quad
h(\{-\infty\}\times \R) = W_1^{\prime} , \quad h(\R \times
\{+\infty\}) = W_2^{\prime}.
\]
\item the restrictions of $h$ to $\R^2$, to $\{\pm \infty\} \times
  \R$, and to $\R \times \{\pm \infty\}$, are diffeomorphisms of class
  $C^1$;
\item if moreover the $(0)$-essential subbundle $\mathcal{E}$ can
be
  lifted to a subbundle $\mathcal{V}$, then
\[
\deg h = -\deg h|_{\{-\infty\}\times \R} \cdot \deg
  h|_{\R \times \{+\infty\}} = \deg h|_{\R \times \{-\infty\}} \cdot \deg
  h|_{\{+\infty\}\times \R},
\]
where $\deg$ denotes the $\Z$-topological degree, referred to the
orientations defined in section \ref{oriente}.
\end{enumerate}
\end{thm}

Concerning (ii), note that it may happen that $y=y^{\prime}$, and
in this case even that $W_1=W_1^{\prime}$ or $W_2=W_2^{\prime}$,
but the last two identities cannot hold simultaneously. When $y
\neq y^{\prime}$, $h$ is injective, so it is a conjugacy.
Statement (iv) expresses the coherence we need between the
orientations of the one-dimensional and two-dimensional
intersections of unstable and stable manifolds. The picture is
completed by the following proposition.

\begin{prop}
\label{c'e'} Assume that the gradient-like Morse vector field $F$
satisfies (C1-2) with respect to a $(0)$-essential subbundle
$\mathcal{E}$ of $TM$. Assume that $F$ has the Morse-Smale
property up to order 2. Let $x,y,z$ be rest points such that
$m(x)=m(y)+1=m(z)+2$, and let $W_1$, $W_2$ be connected
components of $W^u(x)\cap W^s(y)$, $W^u(y)\cap W^s(z)$,
respectively. Then there exists a unique connected component $W$
of $W^u(x)\cap W^s(z)$ such that $\overline{W_1 \cup W_2}$
belongs to the closure of
$\set{\overline{\phi(\R\times\{p\})}}{p\in W}$ with respect to
the Hausdorff distance.
\end{prop}

Both Theorem \ref{conju} and Proposition \ref{c'e'} will be proved
in section \ref{seconju}. The main tool in the proof is the graph
transform method, which allows us to study suitable portions of
$W^u(x)$ and $W^s(z)$ in a neighborhood of another rest point
$y\in \overline{W^u(x) \cap W^s(z)}$.

\section{The boundary homomorphism}
\label{sonasegaio}

Let $(M,\mathcal{E})$ be a pair consisting of a complete
Riemannian Hilbert manifold $M$ of class $C^2$, and of a $C^1$
$(0)$-essential subbundle of $TM$ having an admissible
presentation. Let $F$ be a $C^1$ Morse vector field on $M$,
admitting a non-degenerate Lyapunov function $f$. We shall assume
(PS), (C1-3), the Morse-Smale property up to order 2, and
\begin{description}
\item[(C4)] for every $q\in \Z$, $f$ is bounded below on
$\mathrm{rest}_q(F)=  \set{x\in \rest(F)}{m(x,\mathcal{E})=q}$.
\end{description}

\paragraph{Morse complex with coefficients in $\Z$.} We first consider
the situation in which $\mathcal{E}$ is the $(0)$-essential class
of a subbundle $\mathcal{V}$ of $TM$. In this case we can fix
arbitrary orientations of the Fredholm pairs
$(H^s_x,\mathcal{V}(x))$, for every $x\in \rest(F)$.

Let $x$ and $y$ be rest points with
$m(x,\mathcal{E})-m(y,\mathcal{E})=1$,  and let $W$ be a
connected component of $W^u(x)\cap W^s(y)$. Then $W$ is a flow
line, and it is endowed with the orientation described in section
\ref{oriente}. We can define the number
\[
\sigma(W):=\deg[\phi(\cdot,p):\R\rightarrow W],
\]
for $p\in W$. In other words $\sigma(W)$ equals $+1$ or $-1$
depending on whether $F(p)\in T_p W$ is positively or negatively
oriented. We define also
\[
\sigma(x,y):=\sum_W \sigma(W),
\]
where the sum ranges over all the connected components of
$W^u(x)\cap W^s(y)$.

Now let $x$ and $z$ be rest points with
$m(x,\mathcal{E})-m(z,\mathcal{E})=2$,  and let ${\cal
S}(x,z)$ be the set of broken flow lines from $x$ to $z$ with one
intermediate rest point (necessarily unique and of index
$m(z,\mathcal{E})+1$). By (PS) and by the Morse assumption
there are finitely many rest points $y$ with $f(y)\in ]f(z),
f(x)[$. By Theorem \ref{tcomp} and Proposition \ref{comp}, the
set ${\cal S}(x,z)$ is finite. By Theorem \ref{conju} and
Proposition \ref{c'e'}, there is an involution $\overline{W_1\cup
W_2}\mapsto \overline{W'_1\cup W'_2}$ without fixed points on the
set  ${\cal S}(x,z)$, and
\begin{equation}
\label{x} \sigma(W'_1)\sigma(W'_2)=-\sigma(W_1)\sigma(W_2).
\end{equation}
Let $q\in\Z$ and let $C_q(F)$ be the free Abelian group generated
by the rest points of index $q$:
\[
C_q(F)=\mathrm{span}_\Z \mathrm{rest}_q(F).
\]
Note that $C_q(F)$ may have infinite rank.

Assumption (C4) allows us to define the homomorphism
\[
\partial_q :\  C_q(F)\rightarrow C_{q-1}(F)
\]
by setting for every $x\in\mathrm{rest}_q(F)$
\begin{equation}
\label{bordo}
\partial_q x = \sum_{y\in\mathrm{rest}_{q-1}(F)}\sigma(x,y)y.
\end{equation}
The modules $C_q(F)$ together with the homomorphisms $\partial_q$
are the data of a chain complex. Indeed we have:

\begin{thm}
For every $q\in\Z$, $\partial_{q-1}\circ\partial_q=0$.
\end{thm}

\proof Let $x\in\mathrm{rest}_q(F)$ and $z\in\mathrm{rest}_{q-2}(F)$. The
coefficient of $z$ in $\partial_{q-1}\partial_q x$ is
\[
\sum_{y\in\mathrm{rest}_{q-1}(F)}\sigma(x,y)\sigma(y,z)=\sum_{\overline{W_1\cup
W_2}\in{\cal S}(x,z)} \sigma(W_1)\sigma(W_2),
\]
which is zero by (\ref{x}). \qed

We will call $\{ C_q(F),\partial_q\}_{q\in \Z}$ the Morse complex
of $F$. Clearly, the construction depends on the choice of the
subbundle $\mathcal{V}$ and on the orientation of each Fredholm
pair $(H^s_x,\mathcal{V}(x))$. Changing the subbundle
$\mathcal{V}$ by a compact perturbation changes the Morse complex
by a shift of the indices (when $M$ is connected). A change of
the orientation of $(H^s_x,\mathcal{V}(x))$ produces an isomorphic
Morse complex.

\paragraph{Morse complex with coefficients in $\Z_2$.} In the general
case of a $(0)$-essential subbundle $\mathcal{E}$, statement (iv)
of Theorem \ref{conju} is not available, but we can still define a
Morse complex with $\Z_2$ coefficients. Indeed, defining
$\sigma(x,y)\in \Z_2$ to be the number of connected components of
$W^u(x)\cap W^s(y)$ counted modulo 2, and $C_q(F)$ to be the
$\Z_2$-vector space generated by the rest points of index $q$,
(\ref{bordo}) defines a complex of $\Z_2$-vector spaces.

\section{Proof of the conjugacy theorem}

\paragraph{Construction of $h$ near a broken flow line.}
The main point in the proof of Theorem \ref{conju} and Proposition
\ref{c'e'} is to construct $h$ near a broken gradient flow line.

\begin{prop}
\label{locconju} Assume that the Morse vector field $F$ has a
non-degenerate Lyapunov function $f$, satisfies (C1-2) with
respect to a $(0)$-essential subbundle $\mathcal{E}$ of $TM$, and
satisfies the Morse-Smale condition up to order 2. Let $x,y,z$ be
rest points such that
$m(x,\mathcal{E})=m(y,\mathcal{E})+1=m(z, \mathcal{E})+2$.
Let $W_1$ and $W_2$ be connected components of $W^u(x)\cap W^s(y)$
and $W^u(y)\cap W^s(z)$, respectively. Then there exists a
continuous injective map
\[
h: \Delta := \set{(u,v)\in \overline{\R} \times \overline{\R}}{ v
\leq u} \rightarrow \overline{W^u(x)\cap W^s(z)}
\]
with the following properties:
\begin{enumerate}
\item $\phi_t(h(u,v))=h(u+t,v+t)$, for every $(u,v)\in \Delta$,
$t\in \R$;
\item $h(\Delta \cap \R^2)\subset W^u(x)\cap W^s(z)$,
$h(\R \times \{-\infty\})=W_1$, $h(\{+\infty\} \times \R)=W_2$,
and the restrictions of $h$ to $\Delta \cap \R^2$, to $\R \times
\{- \infty\}$, and to $\{+ \infty\}\times \R$, are diffeomorphisms
of class $C^1$;
\item there exists $\delta>0$ such that for any
$p\in W^u(x)\cap W^s(z)$, if $S=\overline{\phi(\R\times \{p\})}$
has Hausdorff distance less than $\delta$ from $\overline{W_1\cup
W_2}$, then $S\subset h(\Delta)$;
\item if moreover the
$(0)$-essential subbundle $\mathcal{E}$ can be
  lifted to a subbundle $\mathcal{V}$, then
\[
\deg h = -\deg h|_{\R \times \{-\infty\}} \cdot \deg
  h|_{\{+\infty\}\times \R},
\]
where $\deg$ denotes the $\Z$-topological degree, referred to the
orientations defined in section \ref{oriente}.
\end{enumerate}
\end{prop}

Let us identify a neighborhood of $y$ in $M$ with a neighborhood of
$0$ in the Hilbert space $H$, identifying $y$ with $0$.
We endow $H$ with an equivalent Hilbert product $\langle
\cdot, \cdot \rangle$ which is adapted to $\nabla F(y) = DF(0)$ 
(see Appendix C),
and we set $H^u := H^u_y$, $H^s := H^s_y$, so that $H^u \oplus
H^s$ is the splitting of $H$ given by the decomposition of the
spectrum of $\nabla F(y)$ into the subset with positive real part and the
one with negative real part. Let $P^u$ and $P^s$ denote the
corresponding projectors. We shall often identify $H=H^u \oplus
H^s$ with $H^u \times H^s$. By $H^u(r)$, resp.\ $H^s(r)$, we shall
denote the closed $r$-ball centered in $0$ of the linear subspace
$H^u$, resp.\ $H^s$. We set $Q(r) := H^u(r) \times H^s(r)$.
If $X$ and $Y$ are metric spaces and $\theta>0$,
$\mathrm{Lip}_{\theta}(X,Y)$ will denote the space of
$\theta$-Lipschitz maps from $X$ into $Y$, endowed with the $C^0$
topology.

Let $p_1\in W_1$ and $p_2\in W_2$. Choose $\rho_0>0$ so small
that the sets
\[
X:= W^u(x) \cap f^{-1}(f(p_1)) \cap B_{\rho_0}(p_1), \quad Z:=
W^s(z) \cap f^{-1}(f(p_2)) \cap B_{\rho_0}(p_2),
\]
do not contain rest points, and
\[
X \cap W^u(x) \cap W^s(y) = \{p_1\}, \quad Z \cap W^u(y) \cap
W^s(z) = \{p_2\}.
\]
Then $X$ and $Z$ are submanifolds of class $C^1$, and the
Morse-Smale condition implies that $X$ is transverse to $W^s(y)$,
and $Z$ is transverse to $W^u(y)$.

\begin{lem}
\label{osco} For any $\theta>0$ there exist $r_0>0$, $\rho>0$,
$t_0\geq 0$, and two continuous families
\[
\{\sigma_t\}_{t\in [0,+\infty]} \subset \mathrm{Lip}_{\theta}
(H^u(r_0),H^s(r_0)), \quad    \{\tau_t\}_{t\in [-\infty,0]}
\subset \mathrm{Lip}_{\theta} (H^s(r_0),H^u(r_0)),
\]
such that each $\sigma_t$ and each $\tau_t$ is $C^1$, and:
\begin{enumerate}
\item $\phi_{t_0 + t} (X\cap B_{\rho}(p_1)) \cap Q(r_0)=
  \graf \sigma_t$, for every $t\in [0,+\infty[$;
\item $W^u(y) \cap Q(r_0) = \graf \sigma_{+\infty}$; 
\item $\phi_{-t_0 +t} (Z\cap B_{\rho}(p_2)) \cap Q(r_0)=
\graf \tau_t$, for every $t\in ]-\infty,0]$;
\item $W^s(y) \cap Q(r_0) = \graf \tau_{-\infty}$; \item for
any $\theta_1>0$ there exist $r_1\in ]0,r_0]$ and $t_1\geq 0$ such
that
\[
\sigma_t|_{H^u(r_1)} \in \mathrm{Lip}_{\theta_1}
(H^u(r_1),H^s(r_1)), \quad \tau_{-t}|_{H^s(r_1)} \in
\mathrm{Lip}_{\theta_1} (H^s(r_1),H^s u(r_1)),
\]
for any $t\geq t_1$.
\end{enumerate}
\end{lem}

\proof Let $r$ be as small as required by Propositions
\ref{grafici} and \ref{conley}. Since $T_{p_1} X\oplus T_{p_1}
W^s(y)=T_{p_1} M$, the path of subspaces $D \phi_t(p_1) T_{p_1}
X$ converges to $T_y W^u(y)$ for $t\rightarrow +\infty$, by
Theorem \ref{abue} (iii). Therefore, we can find $s_1\geq 0$ such
that $\phi(s_1,p_1)$ is in the interior of $Q(r)$, and
$D \phi_{s_1}(p_1) T_{p_1} X \subset H^u \times H^s$
is the graph of a linear operator from $H^u$ to $H^s$ of norm
strictly less than 1. By the implicit function theorem, there
exists $\rho>0$ such that $\phi_{s_1} (X \cap B_{\rho}(p_1))$
is the graph of a 1-Lipschitz map $\sigma: U \rightarrow H^s(r)$,
where $U\subset H^u(r)$ is open. Moreover, $\graf \sigma
\cap W^s(y)=\{\phi(s_1,p_1)\}$, so by Proposition \ref{grafici}
(v), there exist $s_2\geq 0$ and $\sigma^{\prime}\in
\mathrm{Lip}_1 (H^u(r),H^s(r))$ such that
\[
\graf \sigma^{\prime} = \phi_{s_1+s_2} (X \cap
B_{\rho}(p_1)) \cap Q(r),
\]
where we have also used Proposition \ref{conley}. Let
\[
\Gamma : [0,+\infty] \times \mathrm{Lip}_1 (H^u(r) , H^s(r))
\rightarrow \mathrm{Lip}_1 (H^u(r), H^s(r))
\]
be the graph transform map provided by Proposition \ref{grafici}.
By Proposition \ref{grafici} (iv), there exist $r_0\in ]0,r]$ and
$s_3\geq 0$ such that $\Gamma(t,\sigma^{\prime})\in
\mathrm{Lip}_{\theta} (H^u(r_0), H^s(r_0))$ for any $t\geq s_3$.
Setting $t_0:=s_1 + s_2 + s_3$ and $\sigma_t =
\Gamma(t-s_3,\sigma^{\prime})$ for $t\in [0,+\infty]$, statements
(i), (ii), and the first part of (v) follow immediately from
Propositions \ref{grafici} and \ref{conley}.

Changing the sign of $t$ and considering the evolution of $Z$, we
obtain a family of maps $\{\tau_t\}$ satisfying (iii), (iv), and
the second part of (v).
\qed

\proof (of Proposition \ref{locconju}) Let $\theta$ be a positive
number strictly less than 1, and let
$r_0,\rho,t_0,\sigma_t,\tau_t$ be as in the lemma above. Since
$\theta<1$, the contracting mapping principle implies the
existence of a Lipschitz continuous map
\[
\Lambda : \mathrm{Lip}_{\theta} (H^u(r_0),H^s(r_0)) \times
\mathrm{Lip}_{\theta} (H^s(r_0),H^u(r_0)) \rightarrow Q(r_0),
\]
which associates to $(\sigma,\tau)$ the unique intersection of the
graphs of $\sigma$ and $\tau$, i.e.\ the unique fixed point of the
contraction
\[
H^u(r_0) \times H^s(r_0) \ni (\xi,\eta) \mapsto
(\tau(\eta),\sigma(\xi)) \in H^u(r_0) \times H^s(r_0).
\]
For $(u,v)\in [0,+\infty]\times [-\infty,0]$, set
\begin{eqnarray*}
X_u= \begin{cases} \phi_{u+t_0}(X\cap B_{\rho}(p_1))\cap Q(r_0), &
\mbox{if } u\in [0,+\infty[, \\ W^u(y) \cap Q(r_0), & \mbox{if }
u=+\infty, \end{cases} , \\ Z_v= \begin{cases} \phi_{v-t_0}(Z\cap
B_{\rho}(p_2))\cap Q(r_0), & \mbox{if } v\in ]-\infty,0], \\
W^s(y) \cap Q(r_0), & \mbox{if } v=-\infty.
\end{cases}
\end{eqnarray*}
Then we define $h(u,v)$ to be the unique point of the intersection
$X_u \cap Z_v$. The map $h$ is well defined and continuous on
$[0,+\infty]\times [-\infty,0]$ because it can be written as the
composition
\[
[0,+\infty]\times [-\infty,0] \ni (u,v) \mapsto (\sigma_u,\tau_v)
\stackrel{\Lambda}{\longrightarrow} Q(r_0)\hookrightarrow M.
\]
Since $X$ and $Z$ are contained in level sets of $f$, and they
contain no rest points, $X_u \cap X_{u^{\prime}} = \emptyset$ if
$u\neq u^{\prime}$, and  $Z_v \cap Z_{v^{\prime}} = \emptyset$ if
$v\neq v^{\prime}$. So $h$ is injective. By definition, for every
$(u,v)\in [0,+\infty]\times [-\infty,0]$, and $-u \leq t\leq -v$,
\begin{equation}
\label{canc} \phi_t(h(u,v))= h(u+t,v+t),
\end{equation}
and we can use formula (\ref{canc}) to extend $h$ to a continuous
injective map on
\[
\Delta= \set{(u,v)\in \overline{\R} \times \overline{\R}}{v\leq
u},
\]
still verifying (\ref{canc}), proving (i). Since $X\subset
W^u(x)$ and $Z\subset W^s(z)$, for every $(u,v)\in \Delta\cap
\R^2$ the point $h(u,v)$ belongs to $W^u(x)\cap W^s(z)$.  Since
\[
h(u,-\infty) = \phi(u,p_1), \quad \mbox{and} \quad h(+\infty,v) =
\phi(v,p_2),
\]
$h|_{\R\times \{-\infty\}}$ and $h|_{\{+\infty\}\times \R}$ are
diffeomorphisms of class $C^1$ onto $W_1$ and $W_2$. Since $X$
and $Z$ are of class $C^1$ and so is $\phi$, the implicit
function theorem implies that $h|_{\Delta\cap \R^2}$ is a
diffeomorphism of class $C^1$, proving (ii). Notice that,
differentiating the identities
\[
\phi_t(h(u,v)) = h(u+t,v+t), \quad \phi_t(h(u,-\infty)) =
h(u+t,-\infty), \quad \phi_t(h(+\infty,v)) = h(+\infty,v+t),
\]
with respect to $t$ in $t=0$, we obtain
\begin{equation}
\label{lc} F(h(u,v)) = \frac{\partial h}{\partial u} (u,v) +
\frac{\partial
  h}{\partial v} (u,v), \quad F(u,-\infty) = \frac{\partial
  h}{\partial u} (u,-\infty), \quad  F(+\infty,v) = \frac{\partial
  h}{\partial v} (+\infty,v),
\end{equation}
for every $u\in \R$, $v\in \R$.

Since $y$ is a rest point, and $p_1\in W^s(y)$, $p_2\in W^u(y)$,
we can find $\delta\in ]0,\rho]$ so small that, if $\phi(t_1,p)\in
B_{\delta}(p_1)$ and $\phi(t_2,p)\in B_{\delta}(p_2)$, we have
\begin{equation}
\label{pero} t_2-t_1\geq 2t_0, \quad \phi(t_1+t_0,p)\in Q(r_0),
\quad \phi(t_2-t_0,p)\in Q(r_0).
\end{equation}
If $p\in W^u(x)\cap W^s(z)$ and $S:=\overline{\phi(\R\times
\{p\})}$ has Hausdorff distance less than $\delta$ from
$\overline{W_1\cup W_2}$, we can find $t_1,t_2\in \R$ such that
$\phi(t_1,p)\in B_{\delta}(p_1)$ and $\phi(t_2,p)\in
B_{\delta}(p_2)$. By Proposition \ref{conley}, the set of $t\in
\R$ such that $\phi(t,p)\in Q(r_0)$ is connected, so by
(\ref{pero}),
\begin{equation}
\label{melo} \phi\left( \frac{t_1+t_2}{2},p\right) = \phi\left(
\frac{(t_1+t_0)+(t_2-t_0)}{2},p\right) \in Q(r_0).
\end{equation}
Then, setting $u:= (t_2-t_1)/2-t_0\geq 0$ and $v:=-u\leq 0$, we
obtain
\begin{equation}
\label{sero} \phi\left( \frac{t_1+t_2}{2},p\right) \in
\phi_{u+t_0} (X\cap B_{\rho}(p_1)) \cap \phi_{v-t_0} (Z\cap
B_{\rho}(p_2)).
\end{equation}
So, by (\ref{melo}) and (\ref{sero}), $\phi((t_1+t_2)/2,p)\in X_u
\cap Z_v$, that is $\phi((t_1+t_2)/2,p) = h(u,v)$. By (i) the
whole flow line through $p$ is in $h(\Delta)$ which,
being closed, must contain also $S$, proving conclusion (iii).

In order to prove (iv), we shall need the following:

\begin{lem}
\label{spazi} There exist $u_0\geq 0$, $v_0\leq 0$, and a
continuous map
\[
\mathcal{W}:[u_0,+\infty] \times [-\infty,v_0] \rightarrow \mathrm{Gr}(H),
\]
such that $(\mathcal{W}(u,v), \mathcal{V}(h(u,v)))$ is a Fredholm pair for
every $(u,v) \in [u_0,+\infty] \times [-\infty,v_0]$, and:
\begin{enumerate}

\item $T_{h(u,v_0)} W^s(z) = \R F(h(u,v_0)) \oplus
\mathcal{W}(u,v_0)$, for every $u\in [u_0,+\infty]$;

\item $\mathcal{W}(u,-\infty)=T_{h(u,-\infty)} W^s(y)$, for every
$u\in [u_0,+\infty]$;

\item $D \phi_{v-v_0}(h(+\infty,v_0)) \mathcal{W}(+\infty,v_0)
= \mathcal{W} (+\infty,v)$, for every $v\in ]-\infty,v_0]$;

\item $\mathcal{W}(u_0,v) + T_{h(u_0,v)} W^u(x) = H$
for every $v\in [-\infty,v_0]$, and there exists a non-vanishing
continuous vector field $G:[-\infty,v_0] \rightarrow H$ along
$v\mapsto h(u_0,v)$ such that
\[
\mathcal{W}(u_0,v) \cap T_{h(u_0,v)} W^u(x) = \R G(v), \quad
G(-\infty) = \frac{\partial h}{\partial u}(u_0,-\infty), \quad
G(v_0) = \frac{\partial h}{\partial u}(u_0,v_0),
\]
for every $v\in [-\infty,v_0]$.
\end{enumerate}
\end{lem}

\proof 
Recalling that by (C1) $(H^s,\mathcal{V}(y))$ is a Fredholm pair,
we can find $\theta_1\in ]0,\theta]$ so small that
\begin{equation}
\label{la1} \forall S\in \mathcal{L} (H^s, H^u) \; \mbox{such that
} \|S\| \leq 2 \theta_1, \;\; (\graf S,\mathcal{V}(y)) \mbox{ is a
Fredholm pair.}
\end{equation}
Moreover, we may assume that setting $L:= \nabla F(y) = DF(0)$,
\begin{equation}
\label{la2} \frac{4\theta_1^2}{1-\theta_1^2} + 2\theta_1 <
\frac{1}{\|L\|^2 \|L^{-1}\|^2}.
\end{equation}
Let $r_1$ and $t_1$ be the positive numbers given by Lemma
\ref{osco} (v): there holds
\begin{equation}
\label{la0} \|D\sigma_t(\xi)\| \leq \theta_1, \quad
\|D\tau_t(\eta)\|\leq \theta_1,
\end{equation}
for any $t\geq t_1$, $\xi\in H^u(r_1)$, $\eta \in H^s(r_1)$. Since
\begin{equation}
\label{la3} F(\xi) = L \xi + o(\xi) \quad \mbox{for }
\xi\rightarrow 0,
\end{equation}
the quadratic form
\[
g(\xi) = - \frac{1}{2} \langle L \xi, \xi\rangle
\]
is a Lyapunov function for $F$ in a neighborhood of $0$. Therefore
there exists $u_1\geq 0$ such that the function $u \mapsto
g(h(u,-\infty))$ is strictly decreasing in $[u_1,+\infty[$. So
when $v\rightarrow -\infty$ the functions
\[
[u_1,+\infty[ \ni u \mapsto g(h(u,v)) \in \R
\]
converge uniformly to a strictly decreasing function. Therefore,
using also (\ref{la1}), the fact that $W^s(y)$ is tangent to $H^s$
at $y=0$, (\ref{la2}), (\ref{la0}), and (\ref{la3}), it is easy to
check that there are $u_0\geq 0$ and $v_0\leq 0$ such that:
\renewcommand{\theenumi}{\alph{enumi}}
\renewcommand{\labelenumi}{(\theenumi)}
\begin{enumerate}

\item $Dg(h(h_0,v_0))\left[ \frac{\partial h}{\partial u} (u_0,
v_0)\right] = \frac{\partial}{\partial u} g(h(u,v_0))|_{u=u_0}<0$;

\item $g$ is a Lyapunov function for $F$ on $h([u_0,+\infty]
\times [-\infty,v_0])$;

\item if $S\in \mathcal{L} (H^s, H^u)$ has norm $\|S\| \leq
\theta_1$, and $(u,v) \in [u_0,+\infty] \times [-\infty,v_0]$,
then $(\graf S, \mathcal{V}(h(u,v)))$ is a Fredholm pair;

\item for any $v\in [-\infty,v_0]$ there holds $\|P^u h(u_0,v)\|
\leq \epsilon \|P^s h( u_0,v)\|$, and
\[
\|F(h(u_0,v))\| \leq 2 \|L\| \, \|P^s h(u_0,v)\|, \quad
\|F(h(u_0,v)) - L h(u_0,v)\| \leq \epsilon \|P^s h(u_0,v)\|,
\]
where $\epsilon>0$ is so small that
\begin{equation}
\label{la4} \frac{4\theta_1}{1-\theta_1^2} (\theta_1 + \epsilon) +
2\theta_1 + \epsilon < \frac{1}{\|L\|^2 \|L^{-1}\|^2};
\end{equation}

\item for any $(u,v)\in [u_0,+\infty] \times [-\infty,v_0]$ we
have
\[
\|D\sigma_u (P^u h(u,v))\| \leq \theta_1, \quad \|D\tau_v (P^s
h(u,v))\| \leq \theta_1.
\]

\end{enumerate}
\renewcommand{\theenumi}{\roman{enumi}}
\renewcommand{\labelenumi}{(\theenumi)}

We will define $\mathcal{W}(u,v)$ to be the graph of suitable
linear maps $S(u,v) \in \mathcal{L}(H^s,H^u)$. We start by
defining $S$ on three edges of the square $[u_0,+\infty] \times
[-\infty,v_0]$. For $u\in [u_0,+\infty]$ we set
\[
S(u,v_0) = D\tau_{v_0} (P^s h(u,v_0)), \quad S(u,-\infty) =
D\tau_{-\infty} (P^s h(u,-\infty)),
\]
and for $v\in ]-\infty,v_0]$ we set
\[
S(+\infty,v) = D\tau_{v} (P^s h(+\infty,v)).
\]
By (e), $\|S\|\leq \theta_1$. The map $S$ is clearly continuous on
$[u_0,+\infty] \times \{v_0\}$ and on $[u_0,+\infty] \times
\{-\infty\}$. If $v\in ]-\infty,v_0]$,
\[
\graf S(+\infty,v) = D\phi_{v-v_0} (h(+\infty,v_0)) [\graf
S(+\infty,v_0)],
\]
so $S$ is continuous on $\{+\infty\} \times ]-\infty,v_0]$. By
Theorem \ref{abue} (iii), $S(+\infty,v)$ converges to
$S(+\infty,-\infty)$ for $v\rightarrow -\infty$, so
\[
S : ([u_0,+\infty] \times \{-\infty,v_0\}) \cup (\{+\infty\}
\times [-\infty,v_0]) \rightarrow \mathcal{L}(H^s,H^u)
\]
is a continuous map. We can extend the map $S$ by convexity to a
continuous map on $[u_0,+\infty] \times [-\infty,v_0]$ in such a
way that $\|S(u,v)\|\leq \theta_1$ everywhere, and we set
\[
\mathcal{W}(u,v) = \graf S(u,v).
\]
By (c), $(\mathcal{W}(u,v), \mathcal{V}(h(u,v)))$ is always a
Fredholm pair, and by construction $\mathcal{W}$ satisfies the
requirements (i), (ii), and (iii).

There remains to check (iv). Let $v\in [-\infty,v_0]$. The tangent
space to the unstable manifold of $x$ at $h(u_0,v)$ is
\[
T_{h(u_0,v)} W^u(x) = \R F(h(u_0,v)) \oplus \graf D\sigma_{u_0}
(P^u h(u_0,v)).
\]
Since $S(u_0,v) \in \mathcal{L} (H^s,H^u)$ and $D\sigma_{u_0} (P^u
h(u_0,v)) \in \mathcal{L} (H^u,H^s)$ have norm not exceeding
$\theta_1<1$, we have
\[
\mathcal{W} (u_0,v) + T_{h(u_0,v)} W^u(x) = H.
\]
Moreover, a simple computation shows that the intersection
\[
\mathcal{W} (u_0,v) \cap T_{h(u_0,v)} W^u(x) = \graf S(u_0, v)
\cap \Bigl(\R F(h(u_0,v)) \oplus \graf D\sigma_{u_0} (P^u
h(u_0,v))\Bigr)
\]
is a one-dimensional space spanned by the vector
\[
\tilde{G}(v) = ( \tilde{G}^u(v), \tilde{G}^s (v)) \in H^u \times
H^s,
\]
where
\[
\tilde{G}^u(v) = S(u_0,v) \tilde{G}^s(v), \quad \tilde{G}^s(v) =
(I-T(v)S(u_0,v))^{-1} (P^s F(h(u_0,v)) - T(v) P^u F(h(u_0,v))),
\]
and $T(v) = D\sigma_{u_0} (P^u h(u_0,v))$. Indeed, $\|TS\|\leq
\theta_1^2 <1$, so $I-T S$ is invertible, and
\[
\|(I - T S)^{-1} \| \leq \frac{1}{1-\theta_1^2}, \quad \|(I - T
S)^{-1} -I\| \leq \frac{\theta_1^2}{1-\theta_1^2}.
\]
By (d) we have the estimates
\begin{eqnarray*}
\|\tilde{G}^s\| \leq \frac{1}{1-\theta_1^2} ( \|P^s F(h)\| +
\theta_1 \|P^u F(h)\| ) \leq \frac{2}{1-\theta_1^2} \|F(h)\|
\leq \frac{4}{1-\theta_1^2} \|L\|\, \|P^s h\|, \\
|\langle LP^u h,\tilde{G}^u\rangle| \leq \epsilon \|L\|\,
\|\tilde{G}^u\| \,\|P^s h\| \leq \epsilon \theta_1 \|L\|\,
\|\tilde{G}^s\| \,\|P^s h\| \leq \frac{4\epsilon
\theta_1}{1-\theta_1^2} \|L\|^2 \|P^s h\|^2, \\
\| \tilde{G}^s -P^s F(h) \| \leq \| \tilde{G}^s - (P^s F(h) - TP^u
F(h))\| + \|TP^u F(h)\| \leq \frac{\theta_1^2}{1-\theta_1^2} \|P^s
F(h) - TP^u F(h)\| \\ + 2\theta_1 \|L\| \, \| P^s h\| \leq \left(
\frac{4\theta_1^2}{1-\theta_1^2} + 2\theta_1 \right) \|L\|\,
\|P^s h\|, \\
\| \tilde{G}^s - LP^s h\| \leq \|\tilde{G}^s - P^s F(h)\| + \|P^s
F(h) - P^s L h\| \leq \left( \Bigl(
\frac{4\theta_1^2}{1-\theta_1^2} + 2\theta_1 \Bigr) \|L\|
+\epsilon \right) \|P^s h\|.
\end{eqnarray*}
Thus
\begin{eqnarray*}
Dg(h) [\tilde{G} ] = - \langle L h , \tilde{G} \rangle = - \langle
L P^u h, \tilde{G}^u \rangle - \langle L P^s h, \tilde{G}^s
\rangle \leq \frac{4 \epsilon \theta_1}{1-\theta_1^2} \|L\|^2
\|P^s h\|^2 - \langle L P^s h, L P^s h\rangle \\ + \| L \| \,
\|P^s h\| \, \| \tilde{G}^s - L P^s h\| \leq \left( \Bigl( \frac{4
\theta_1}{1-\theta_1^2} (\theta_1 + \epsilon) + 2\theta_1 +
\epsilon \Bigr) \|L\|^2 - \frac{1}{\|L^{-1}\|^2} \right) \|P^s
h\|^2,
\end{eqnarray*}
and by (\ref{la4}) we get
\begin{equation}
\label{la5} Dg(h(u_0,v))[\tilde{G} (v) ]<0 \quad \forall v\in
[-\infty,v_0].
\end{equation}
On the other hand, by (a),
\begin{equation}
\label{la6} Dg(h(u_0,v_0)) \left[ \frac{\partial h}{\partial u}
(u_0,v_0) \right] <0,
\end{equation}
and by (\ref{lc}) and (b),
\begin{equation}
\label{la7} Dg(h(u_0,-\infty)) \left[ \frac{\partial h}{\partial
u} (u_0, -\infty) \right] = Dg (h(u_0,-\infty))
[F(h(u_0,-\infty))] <0.
\end{equation}
By construction, $h(u,v_0) \in \graf \tau_{v_0} \cap W^u(x)$, for
every $u\geq 0$, so
\[
\frac{\partial h}{\partial u} (u_0,v_0) \in \graf D\tau_{v_0} (P^s
h(u_0,v_0)) \cap T_{h(u_0,v_0)} W^u(x) = \mathcal{W} (u_0,v_0)
\cap T_{h(u_0,v_0)} W^u(x).
\]
Moreover,
\[
\frac{\partial h}{\partial u} (u_0,-\infty) \in T_{h(u_0,-\infty)}
W^s(y) \cap T_{h(u_0,-\infty)} W^u(x) = \mathcal{W} (u_0,-\infty)
\cap T_{h(u_0,-\infty)} W^u(x).
\]
Then (\ref{la5}), (\ref{la6}), and (\ref{la7}) imply that a vector
field $G$ satisfying the requirements of (iv) can be defined by
multiplying $\tilde{G}$ by a suitable positive function. \qed

\begin{rem}
In some particular cases, such as when $F$ is linear in a
neighborhood of $y=0$, a map $\mathcal{W}$ satisfying the
requirements of the above lemma can be defined simply as
\[
\mathcal{W}(u,v) = \graf D\tau_v (P^s h(u,v)),
\]
providing us with a drastic simplification of the proof. However
in general, the above expression does not define a continuous map
$\mathcal{W}$, the reason being that the graph transform $\Gamma$
of Proposition \ref{grafici} needs not be continuous with respect
to the $C^1$ topology.
\end{rem}

We are now ready to prove assertion (iv) of Proposition
\ref{locconju}. The continuous maps
\begin{eqnarray*}
\tau_{W_1} :[u_0,+\infty[ \rightarrow \mathrm{Gr}_{1,\infty}(TM), \quad
u\mapsto T_{h(u,-\infty)} W_1 = \R F(h(u,-\infty)), \\
\tau_{W_2} :]-\infty,v_0] \rightarrow \mathrm{Gr}_{1,\infty}(TM), \quad
v\mapsto T_{h(+\infty,v)} W_2 = \R F(h(+\infty,v)), \\
\tau_W : [u_0,+\infty[ \times ]-\infty,v_0] \rightarrow
\mathrm{Gr}_{2,\infty} (TM), \quad (u,v) \mapsto T_{h(u,v)} W,
\end{eqnarray*}
have continuous liftings $\widehat{\tau_{W_1}}$,
$\widehat{\tau_{W_2}}$ to $\mathrm{Or}(\mathrm{Gr}_{1,\infty}(TM))$, and
$\widehat{\tau_W}$ to $\mathrm{Or}(\mathrm{Gr}_{2,\infty}(TM))$, corresponding to
the orientations of $W_1$, $W_2$, and $W$, defined in section
\ref{oriente}. Moreover, the continuous map
\[
\omega:[u_0,+\infty] \times ]-\infty,v_0] \rightarrow \mathrm{Fp}(TM),
\quad (u,v) \mapsto (T_{h(u,v)} W^s(z),\mathcal{V}(h(u,v))),
\]
has a continuous lifting $\hat{\omega}$ to
$\mathrm{Or}(\mathrm{Fp}(TM))$, 
which is determined by the orientation $o_z$ of $(H^s_z,\mathcal{V}(z))$.

By Lemma \ref{spazi} (ii), $\mathcal{W}(+\infty,-\infty)=H^s_y$,
so the continuous map
\[
\alpha: [u_0,+\infty] \times [-\infty,v_0] \rightarrow \mathrm{Fp}(TM),
\quad (u,v) \mapsto (\mathcal{W}(u,v), \mathcal{V}(h(u,v))),
\]
has a unique continuous lifting $\hat{\alpha}$ to 
$\mathrm{Or}(\mathrm{Fp}(TM))$
such that $\hat{\alpha}(+\infty,-\infty)=o_y$. By Lemma
\ref{spazi} (iv), there is a continuous curve $\mathcal{X}:
[-\infty,v_0] \rightarrow \mathrm{Gr}(TM)$ such that $\mathcal{W}(u_0,v) =
\R G(v) \oplus \mathcal{X}(v)$ for every $v\in [-\infty,v_0]$,
and we can define the continuous map
\[
\beta: [-\infty,v_0] \rightarrow \mathrm{Fp}(TM), \quad v \mapsto
(\mathcal{X}(v), \mathcal{V}(h(u_0,v))).
\]
Then $\mathcal{X}(v)$ is a linear supplement of $T_{h(u_0,v)}
W^u(x)$ in $T_{h(u_0,v)} M$, so Theorem \ref{abue} (iii) implies
that
\[
\lim_{t\rightarrow -\infty} D \phi_t (h(u_0,v)) \mathcal{X}(v) =
H^s_x,
\]
uniformly in $v\in [-\infty,v_0]$. Therefore, the orientation
$o_x$ of $(H^s_x,\mathcal{V}(x))$ determines a continuous lifting
$\hat{\beta} : [-\infty,v_0] \rightarrow \mathrm{Or}(\mathrm{Fp}(TM))$ of $\beta$.

Denote by $h_1$ and by $h_2$ the restrictions of $h$ to $\R \times
\{-\infty\}$ and to $\{+\infty\} \times \R$. If $X$ is an
$n$-dimensional real vector space and $\xi$ is a non-zero element
of $\Lambda^n(X)$, the same symbol $\xi$ will also denote the
orientation of $X$ induced by $\xi$. When $n=1$, we shall
identify $\Lambda^1(X)$ with $X$. If $o$ is an orientation of
$X$, $-o$ will denote the other orientation.

By Lemma \ref{spazi} (i), (iv), and by (\ref{lc}),
\[
T_{h(u_0,v_0)} W^s(z)= \R F(h(u_0,v_0)) \oplus \R G(v_0) \oplus
\mathcal{X}(v_0)  = T_{h(u_0,v_0)} W \oplus \mathcal{X}(v_0),
\]
so by the definition of the orientation of $W\subset W^u(x) \cap
W^s(z)$,
\begin{equation}
\label{struno} \hat{\omega}(u_0,v_0) = \widehat{\tau_W} (u_0,v_0)
{\textstyle \bigwedge}  \hat{\beta}(v_0) = (\deg h)\,  \left(
\frac{\partial h}{\partial u}
  (u_0,v_0) \wedge \frac{\partial h}{\partial v}(u_0,v_0)\right)
  {\textstyle \bigwedge}  \hat{\beta}(v_0).
\end{equation}
Moreover, $T_{h(+\infty,v_0)} W^s(z)= T_{h(+\infty,v_0)} W_2
\oplus \mathcal{W}(+\infty,v_0)$, so by Lemma \ref{spazi} (iii)
and by the definition of the orientation of $W_2\subset W^u(y)
\cap W^s(z)$,
\begin{eqnarray*}
\hat{\omega} (+\infty,v_0) = \widehat{\tau_{W_2}} (v_0){\textstyle
\bigwedge}
  \hat{\alpha} (+\infty,v_0) = (\deg h_2) \, \frac{\partial
  h}{\partial v} (+\infty,v_0) {\textstyle \bigwedge}  \hat{\alpha} (+\infty,v_0) \\=
(\deg h_2) \,S F(h(+\infty,v_0) {\textstyle \bigwedge}
\hat{\alpha} (+\infty,v_0),
\end{eqnarray*}
where we have taken (\ref{lc}) into account. By Lemma \ref{spazi}
(i),
\[
\omega(u,v_0) = (\R F(h(u,v_0)) \oplus \alpha_1(u,v_0), \alpha_2
(u,v_0)),
\]
for every $u\in [u_0,+\infty]$, so by the continuity of the
product on the orientation bundle we obtain
\[
\hat{\omega} (u_0,v_0) = (\deg h_2)\, F(h(u_0,v_0)){\textstyle
\bigwedge} \hat{\alpha} (u_0,v_0).
\]
Hence by (\ref{lc}),
\begin{equation}
\label{strue} \hat{\omega} (u_0,v_0) = (\deg h_2)\, \left(
\frac{\partial h}{\partial u}
  (u_0,v_0) + \frac{\partial h}{\partial u} (u_0,v_0) \right) {\textstyle \bigwedge}
  \hat{\alpha}(u_0,v_0).
\end{equation}
By Lemma \ref{spazi} (ii), (iv), and by (\ref{lc}),
\[
T_{h(u_0,-\infty)} W^s(y) = \mathcal{W}(u_0,-\infty) = \R
G(-\infty) \oplus \mathcal{X}(-\infty) = T_{h(u_0,-\infty)} W_1
\oplus \mathcal{X}(-\infty),
\]
so by the definition of the orientation of $W_1 \subset W^u(x)
\cap W^s(y)$,
\begin{eqnarray*}
\hat{\alpha} (u_0,-\infty) = \widehat{\tau_{W_1}} (u_0)
{\textstyle \bigwedge} \hat{\beta}(-\infty)) = (\deg h_1)\,
\frac{\partial h}{\partial
    u} (u_0,-\infty) {\textstyle \bigwedge}  \hat{\beta}(-\infty) \\ = (\deg h_1)\,
G(-\infty) {\textstyle \bigwedge}  \hat{\beta}(-\infty).
\end{eqnarray*}
Then by the identity
\[
\alpha(u_0,v) = (\R G(v) \oplus \beta_1(v) , \beta_2(v)), \quad
\forall v\in [-\infty,v_0],
\]
and by the continuity of the product on the orientation bundle we
obtain
\begin{equation}
\label{stre} \hat{\alpha}(u_0,v_0)= (\deg h_1)\, G(v_0)
{\textstyle \bigwedge} \hat{\beta}(v_0) = (\deg h_1)\,
\frac{\partial h}{\partial u} (u_0,v_0) {\textstyle \bigwedge}
\hat{\beta}
  (v_0).
\end{equation}
Identities (\ref{strue}) and (\ref{stre}), together with the
associativity of the product of orientation, imply that
\begin{eqnarray*}
\hat{\omega} (u_0,v_0) = (\deg h_2) (\deg h_1) \, \left(
  \frac{\partial h}{\partial u} (u_0,v_0) + \frac{\partial h}{\partial
  v} (u_0,v_0)\right) {\textstyle \bigwedge}  \left( \frac{\partial h}{\partial u}
  (u_0,v_0) {\textstyle \bigwedge}
  \hat{\beta}(v_0) \right) \\
= (\deg h_2) (\deg h_1)\, \left( \Bigl( \frac{\partial h}{\partial
u} (u_0,v_0) + \frac{\partial h}{\partial v} (u_0,v_0) \Bigr)
\wedge \frac{\partial h}{\partial u} (u_0,v_0) \right) {\textstyle
\bigwedge} \hat{\beta}(v_0)
\\ = - (\deg h_2) (\deg h_1)\, \left(\frac{\partial h}{\partial
u} (u_0,v_0) \wedge
  \frac{\partial h}{\partial v} (u_0,v_0)\right) {\textstyle \bigwedge}
  \hat{\beta}(v_0),
\end{eqnarray*}
and comparing the above identity with (\ref{struno}) we obtain
\[
\deg h = - (\deg h_1) (\deg h_2),
\]
proving (iv). \qed

\paragraph{Conclusion.}
\label{seconju}

\proof (of Theorem \ref{conju}) Fix a value $c\in ]f(z), f(x)[$.
By assumption, $W\cap \{f=c\}\cong W/\R$ is an open interval, so
it is parameterized by a $C^1$ diffeomorphism $\gamma:\R
\rightarrow M$. By Proposition \ref{comp}, there exist an
increasing, unbounded sequence $(s_n)$ and a broken gradient flow
line $S^+$ from $x$ to $z$ such that
\begin{equation}
\label{z} \lim_{n\rightarrow \infty}
\overline{\phi(\R\times\{\gamma(s_n)\})} =S^+
\end{equation}
in the Hausdorff distance, and also such that $\gamma(s_n)$
converges to a point $p\in S^+$. Since $\gamma$ is a
homeomorphism, $p$ is not in $W$, and since $W$ is closed in
$W^u(x)\cap W^s(z)$, $p$ is not in  $W^u(x)\cap W^s(z)$ either.
So $S^+$ contains a rest point $y$ of  intermediate level. As
already noticed, the Morse-Smale property and the fact that
$m(x,\mathcal{E})=m(z,\mathcal{E})+2$ imply that
$m(y,\mathcal{E})=m(z,\mathcal{E})+1$, and that there exist
$W^+_1$ and $W^+_2$, connected components of $W^u(x)\cap W^s(y)$
and of $W^u(y)\cap W^s(z)$, such that $S^+=\overline{W^+_1\cup
W^+_2}$. Proposition \ref{locconju} provides us with a map
\[
h^+:\Delta^+:=
\set{(u,v)\in\overline{\R}\times\overline{\R}}{v\leq
  u} \rightarrow\overline{W^u(x)\cap W^s(z)},
\]
verifying properties (i) to (iv). In particular by (iii) and
(\ref{z}), $h^+(\Delta^+\cap\R^2)\subset W$, and by (iv),
\begin{equation}
\label{zz} \deg(h^+) = -\deg(h^+|_{\R \times\{-\infty\}})\cdot
\deg(h^+|_{\{+\infty\}\times\R}).
\end{equation}
By Proposition \ref{locconju} (i), for any $t\geq 0$,
\begin{eqnarray*}
\lim_{s\rightarrow -\infty} f(h^+(s+t,s)) = \lim_{s\rightarrow
  -\infty} f(\phi(s,h^+(t,0))) = f(x) >c, \\
\lim_{s\rightarrow +\infty} f(h^+(s+t,s)) = \lim_{s\rightarrow
  +\infty} f(\phi(s,h^+(t,0))) = f(z) <c, \\
\frac{\partial}{\partial s}\left[f(h^+(s+t,s))\right]=
\frac{\partial}{\partial s}\left[f(\phi(s, h^+(t,0))\right]=
Df(\phi_s(h^+(t,0))) [ F(\phi_s(h^+(t,0))] <0,
\end{eqnarray*}
so by the implicit function theorem there exists a function
$\eta_+\in C^1([0,+\infty[,\R)$ such that
\[
f(h^+(\eta_+(t)+t,\ \eta_+(t))) = c,\quad\forall t\geq0.
\]
Then $h^+(\eta_+(t)+t,\ \eta_+(t))\in W\cap\{f=c\}$, and
\[
\theta_+(t):=\gamma^{-1}(h^+(\eta_+(t)+t,\ \eta_+(t)))
\]
defines a $C^1$ function $\theta_+:[0,+\infty[ \rightarrow\R$. An
application of $\phi_{v-\eta_+(u-v)}$ to $h^+(\eta_+(u-v)+t,\
\eta_+(u-v))=\gamma(\theta_+(u-v))$ yields to the the
representation
\begin{equation}
\label{zzz} h^+(u,v)= \phi(v-\eta_+(u-v),
\gamma(\theta_+(u-v))),\quad (u,v)\in\Delta^+\cap\R^2.
\end{equation}
Since $h^+$ is a diffeomorphism, the vectors $\partial
h^+/\partial u$ and $\partial h^+/\partial v$ are linearly
independent, so
\[
(\gamma\circ\theta_+)'= (1+\eta_+^{\prime}) \frac{\partial
  h^+}{\partial u}  + \eta_+^{\prime} \frac{\partial h^+}{\partial v},
\]
never vanishes, and from the fact that $\gamma$ is a
diffeomorphism we deduce that $\theta_+^{\prime}(t)\neq 0$ for
every $t\geq0$. Moreover, from Proposition \ref{locconju} (iii),
$\phi(\R\times\{\gamma(s_n)\})\subset h_+(\Delta^+)$ for $n$
large, which implies that $\phi(\R\times\{\gamma(s_n)\})=
\phi(\R\times\{h^+(t_n,0)\})=
\phi(\R\times\{\gamma(\theta_+(t_n))\})$ for some $t_n\geq0$.
Since $\gamma$ is injective and meets any flow line at most once,
the last equality implies that $\theta_+(t_n)=s_n\rightarrow
+\infty$. Therefore
\begin{equation}
\label{zzzz} \theta_+'>0,\quad\quad
\lim_{t\to+\infty}\theta_+(t)=+\infty.
\end{equation}
The same construction, starting with a sequence $s'_n\rightarrow
-\infty$, yields to a rest point $y'$ with $m(y',\mathcal{E})=
m(z,\mathcal{E})+1$ such that
$\overline{\phi(\R\times\{\gamma(s'_n)\}}$ converges to $S^-=
\overline{W^-_1\cup W^-_2}$, for some connected components
$W^-_1$ and  $W^-_2$ of $W^u(x)\cap W^s(y')$ and $W^u(y')\cap
W^s(z)$, respectively. As before we obtain a map
\[
h^-:\Delta^-:=\{(u,v)\in\overline{\R}\times\overline{\R}\ :\ v\geq
  u\}\rightarrow\overline{W}
\]
where we also used the orientation reversing change of variables
$\Delta^-\ni(u,v)\mapsto (v,u)\in\Delta^+$. Hence
\begin{equation}
\label{zzzzz} \deg(h^-) = \deg(h^-|_{\{-\infty\}\times\R})\cdot
 \deg(h^-|_{\R \times\{+\infty\}})
\end{equation}
and we have the representation
\begin{equation}
\label{zzzzzz} h^-(u,v)= \phi(v-\eta_-(u-v),
\gamma(\theta_-(u-v))),\quad (u,v)\in\Delta^-\cap\R^2.
\end{equation}
for suitable $C^1$ functions $\eta_-$ and $\theta_-$ on
$]-\infty,0]$, with
\begin{equation}
\label{zzzzzzz} \theta_-'>0\quad\quad
\lim_{t\to-\infty}\theta_-(t)=-\infty.
\end{equation}

Proposition \ref{locconju} (iii) together with (\ref{zzz}) and
(\ref{zzzzzz}), implies that $S^-\neq S^+$, as claimed in (ii).
Now we can choose two $C^1$ functions $\eta,\
\theta:\R\rightarrow\R$, with $\theta'>0$, coinciding with
$\eta_-$, $\theta_-$ in a neighborhood of $-\infty$ and with
$\eta_+$, $\theta_+$ in a neighborhood of $+\infty$. The map
\[
h(u,v):= \phi(v-\eta(u-v), \gamma(\theta(u-v))),\quad (u,v)\in
\R^2.
\]
has a continuous extension to $\overline{\R}\times\overline{\R}$
and clearly satisfies all the requirements (i) to (iv).\qed

\medskip

\noindent \proof (of Proposition \ref{c'e'}) The conclusion
follows immediately from Proposition \ref{locconju} (iii). \qed

\renewcommand{\thesection}{\Alph{section}}
\setcounter{section}{0}

\section{Appendix - Infinite dimensional Grassmannians}

The aim of this appendix is to gather the definitions and the
relevant properties of some infinite dimensional Grassmannians.
Unless otherwise stated, detailed proofs can be found in
\cite{ama03b} (but see also
\cite{pal65,luf67,qui85,sw85,ps86,cjs95,shu96}).

\paragraph{The Hilbert Grassmannian and the space of Fredholm pairs.}
By $\mathcal{L}(E,F)$, respectively $\mathcal{L}_c(E,F)$, we will
denote the space of continuous linear, respectively compact
linear, maps from the Banach space $E$ to the Banach space $F$.
If $F=E$ we will use the abbreviations $\mathcal{L}(E)$ and
$\mathcal{L}_c(E)$. The norm of the operator $T\in
\mathcal{L}(E,F)$ will be denoted by $\|T\|$. By $\sigma(L)$ and
by $\sigma_{\mathrm{ess}}(L)$ we will denote the spectrum and the
essential spectrum of the operator $L\in \mathcal{L}(E)$, that is the
spectrum of $[L]$ in the Calkin algebra
$\mathcal{L}(E)/\mathcal{L}_c(E)$.

Let $H$ be a real infinite dimensional separable Hilbert space.
The orthogonal projection onto a closed subspace $V\subset H$
will be denoted by $P_V$, while $V^{\perp}$ will denote the
orthogonal complement of $V$ in $H$.

Let $\mathrm{Gr}(H)$ be the {\em Grassmannian of} $H$, i.e.\ the set of
closed linear subspaces of $H$. The assignment $V\mapsto P_V$ is
an inclusion of $\mathrm{Gr}(H)$ into $\mathcal{L}(H)$, onto the closed
subset of the orthogonal projectors of $H$. We can therefore
define, for any $V,W\in \mathrm{Gr}(H)$ the distance
\[
\dist(W_1,W_2) := \|P_{W_1} - P_{W_2}\|,
\]
which makes $\mathrm{Gr}(H)$ a complete metric space. It can be proved
that $\mathrm{Gr}(H)$ is an analytic Banach submanifold of the Banach space
$\mathcal{L}(H)$: indeed, the subspace of symmetric idempotent
elements of a $C^*$-algebra is always an analytic Banach
submanifold.

The connected components of $\mathrm{Gr}(H)$ are the subsets
\[
\mathrm{Gr}_{n,k}(H):= \set{V\in \mathrm{Gr}(H)}{\dim V=n, \; \codim V=k}, \quad
n,k\in \N\cup \{\infty\}, \; n+k = \infty.
\]
The orthogonal group $\mathrm{O}(H)$ is contractible, by a well known
result by Kuiper \cite{kui65}, and it acts transitively on each of
these components. These facts imply that $\mathrm{Gr}_{\infty,\infty}(H)$
is contractible, while $\mathrm{Gr}_{n,\infty}(H)$ and $\mathrm{Gr}_{\infty,n}(H)$
have the homotopy type of $\mathrm{BO}(n)$, the classifying space of the
orthogonal group of $\R^n$.

A pair $(V,W)$ of closed subspaces of $H$ is said a {\em Fredholm
  pair} if $V\cap W$ is finite dimensional, and $V+W$ is
finite codimensional (see also \cite{kat80}, section IV \S 4). In
this situation, the {\em index} of $(V,W)$ is the number
\[
\ind (V,W) = \dim V\cap W - \codim (V+W).
\]
The set of Fredholm pairs in $H$ will be denoted by $\mathrm{Fp}(H)$: it is
open in $\mathrm{Gr}(H)\times \mathrm{Gr}(H)$, and the index is a continuous
function on $\mathrm{Fp}(H)$. The connected components of $\mathrm{Fp}(H)$ are the
subsets
\begin{eqnarray*}
\mathrm{Gr}_{n,\infty}(H)\times \mathrm{Gr}_{\infty,m}(H), \quad
\mathrm{Gr}_{\infty,n}(H)\times \mathrm{Gr}_{m,\infty}(H), \quad n,m\in \N, \\
\mathrm{Fp}_k^*(H) := \set{(V,W)\in \mathrm{Fp}(H)}{V,W\in \mathrm{Gr}_{\infty,\infty}(H),
\; \ind
  (V,W)=k}, \quad k\in \Z.
\end{eqnarray*}
The space of Fredholm pairs consisting of infinite dimensional
spaces will be denoted by
\[
\mathrm{Fp}^*(H) := \bigcup_{k\in \Z} \mathrm{Fp}_k^*(H).
\]

It can be proved that $\mathrm{Fp}_k^*(H)$ has the homotopy type of
$\mathrm{BO}(\infty)$, the classifying space of the infinite real
orthogonal group $\mathrm{O}(\infty) = \lim_{n\rightarrow \infty} 
\mathrm{O}(n)$. So
$\mathrm{Fp}^*(H)$ is homotopically equivalent to $\Z \times \mathrm{BO}(\infty)$,
and by the Bott periodicity theorem we get
\begin{equation}
\label{bott} \pi_i(\mathrm{Fp}^*(H)) = \begin{cases} \Z & \mbox{for }
i\equiv 0,4 \mod
  8, \\ \Z_2 & \mbox{for } i\equiv 1,2 \mod
  8, \\ 0 & \mbox{for } i\equiv 3,5,6,7 \mod
  8. \end{cases}
\end{equation}

We conclude this section with a result about the existence of
hyperbolic rotations, which will be useful in Appendix B.

\begin{prop}
\label{hr} Let $V,W\in \mathrm{Gr}_{\infty,\infty}(H)$ be such that $\dist
(V,W)<1$. Then there exists $A\in \mathcal{L}(H)$ self-adjoint,
invertible, with $\sigma_{\mathrm{ess}} (A) \cap \R^- \neq
\emptyset$, $\sigma_{\mathrm{ess}} (A) \cap \R^+ \neq \emptyset$,
such that $e^A V=W$.
\end{prop}

\proof Since $\dist(V,W)<1$, $W=\graf L$, with $L=P_{V^{\perp}}
(P_V|_W)^{-1} \in \mathcal{L}(V,V^{\perp})$. Consider the
self-adjoint bounded operator
\[
S = \left( \begin{array}{cc} \theta & \eta L^* \\ \eta L &
1/\theta
\end{array} \right), \quad 0<\theta <1, \quad \eta\in \R,
\]
in the splitting $H=V\oplus V^{\perp}$. Then
\[
(S-\theta)(S-1/\theta) = \eta^2 \left( \begin{array}{cc} L^*L & 0
\\ 0
  & LL^* \end{array} \right).
\]
We fix a $\theta< 1/\|L\|$, so the positive self-adjoint operator on
the right-hand side has its spectrum in $[0,1[$, for every $\eta\in
[0,\theta]$. The spectral mapping theorem implies that
\[
\set{(s-\theta)(s-1/\theta)}{s\in \sigma(S)} \subset [0,1[,
\]
so we have
\[
\sigma(S) \subset ]0,\theta] \; \cup \; [1/\theta,1/\theta+\theta[ \;
\subset \; ]0,1[ \; \cup \; ]1,+\infty[,
\]
for any $\eta\in [0,\theta]$. For $\eta=0$,
$\sigma_{\mathrm{ess}} (S) = \{\theta,1/\theta\}$, so by the
semi-continuity of the essential spectrum
\[
\sigma_{\mathrm{ess}} (S) \cap ]0,1[ \neq \emptyset, \quad
\sigma_{\mathrm{ess}} (S) \cap ]1,+\infty[ \neq \emptyset,
\]
for any $\eta\in [0,\theta]$. In particular for $\eta=\theta$,
$A=\log S$ is a well-defined operator satisfying the requirements.
\qed

\paragraph{The determinant and the orientation of Fredholm pairs.}
Let $n\in \N$. The Grassmannian of $n$-dimensional linear
subspaces $\mathrm{Gr}_{n,\infty}(H)$ is the base space of a non-trivial
real line bundle, the determinant bundle
\[
\Det(\mathrm{Gr}_{n,\infty}(H)) \rightarrow \mathrm{Gr}_{n,\infty}(H),
\]
whose fiber at $X\in \mathrm{Gr}_{n,\infty}(H)$ is the line $\Det(X):=
\Lambda^{\dim
  X}(X)$, the component of the exterior algebra of $X$ consisting of
  tensors of top degree. Such a line bundle has a natural analytic
  structure. Its $\Z_2$ reduction is the non-trivial double covering
\[
\mathrm{Or}(\mathrm{Gr}_{n,\infty}(H)) \rightarrow  \mathrm{Gr}_{n,\infty}(H),
\]
called the orientation bundle, whose fiber at $X$ is the set
$\mathrm{Or}(X)$ consisting of the two elements of $\Det(X)\setminus
\{0\}/\R^+$. If $o_X$ is an element of $\mathrm{Or}(X)$, the other element
will be denoted by $-o_X$. If $n,m\in \N$, the space
\[
\mathcal{S}(n,m) = \set{(X,Y)\in \mathrm{Gr}_{n,\infty}(H) \times
  \mathrm{Gr}_{m,\infty}(H)}{X\cap Y = (0)}
\]
is the base space of the line bundle
\[
\Det (\mathcal{S}(n,m)) \rightarrow \mathcal{S}(n,m),
\]
whose fiber at $(X,Y)$ is the line $\Det(X) \otimes \Det(Y)$, and
the exterior product of tensors of top degree defines an analytic
morphism
\[
\wedge : \Det (\mathcal{S}(n,m)) \rightarrow
\Det(\mathrm{Gr}_{n+m,\infty}(H)), \quad \omega_X \otimes \omega_Y \mapsto
\omega_X \wedge \omega_Y,
\]
which lifts the analytic map $(X,Y) \rightarrow X+Y$. This
operation is associative. The morphism of line bundles $\wedge$
induces a morphism of coverings, denoted by the same symbol,
between the orientation bundles:
\[
\wedge:\mathrm{Or}(\mathcal{S}(n,m)) \rightarrow \mathrm{Or}(\mathrm{Gr}_{n+m,\infty}(H)),
\]
where the first space is the total space of the covering over
$\mathcal{S}(n,m)$ whose fiber at $(X,Y)$ is $\mathrm{Or}(X) \times \mathrm{Or}(Y)$.
The product of orientations satisfies the identity
\[
o_X \wedge o_Y = (-o_X) \wedge (-o_Y) = -(-o_X) \wedge o_Y = - o_X
\wedge (-o_Y),
\]
and it is associative.

\medskip

These constructions have a natural extension to the space of
Fredholm pairs. The {\em determinant bundle over $\mathrm{Fp}(H)$} is the
line bundle
\[
\Det(\mathrm{Fp}(H)) \rightarrow \mathrm{Fp}(H),
\]
whose fiber at $(V,W)$ is the line
\[
\Det(V,W) := \Det (V\cap W) \otimes \Det \left( \Bigl(\frac{H}{V+W}\Bigr) 
\right)^*.
\]
Although the intersection $V\cap W$ and the sum $V+W$ do not
depend even continuously on $(V,W)$, it can be shown that the
above bundle has an analytic structure. This line bundle is also
non-trivial, and its $\Z_2$ reduction is the non-trivial double
covering
\[
\mathrm{Or}(\mathrm{Fp}(H)) \rightarrow  \mathrm{Fp}(H),
\]
called the {\em orientation bundle over $\mathrm{Fp}(H)$}, whose fiber at
$(V,W)$ is the set $\mathrm{Or}(V,W)$ consisting of the two elements of
$\Det(V,W)\setminus \{0\}/\R^+$. If $o_{(V,W)}$ is an element of
$\mathrm{Or}(V,W)$, the other element will be denoted by $-o_{(V,W)}$.
Note that the fundamental group of each component of $\mathrm{Fp}^*(H)$ is
$\Z_2$, so the restriction of the orientation bundle to $\mathrm{Fp}^*(H)$
is the universal covering of $\mathrm{Fp}^*(H)$.

If $n\in \N$, the space
\[
\mathcal{S}(n,\mathrm{Fp}) = \set{(X,(V,W))\in \mathrm{Gr}_{n,\infty}(H) \times
  \mathrm{Fp}(H)}{X\cap V = (0)}
\]
is the base space of the line bundle
\[
\Det (\mathcal{S}(n,\mathrm{Fp})) \rightarrow \mathcal{S}(n,\mathrm{Fp}),
\]
whose fiber at $(X,(V,W))$ is the line $\Det(X) \otimes
\Det(V,W)$, and there is an analytic morphism
\[
{\textstyle \bigwedge}  : \Det (\mathcal{S}(n,\mathrm{Fp})) \rightarrow
\Det(\mathrm{Fp}(H)), \quad \omega_X \otimes \omega_{(V,W)} \mapsto
\omega_X {\textstyle \bigwedge} \omega_{(V,W)}
\]
which lifts the analytic map $(X,(V,W)) \rightarrow (X+V,W)$. Also
this operation is associative, meaning that
\[
\omega_X {\textstyle \bigwedge}  (\omega_Y {\textstyle \bigwedge}
\omega_{(V,W)}) = (\omega_X \wedge \omega_Y) {\textstyle
\bigwedge} \omega_{(V,W)},
\]
for any $\omega_X\in \Det(X)$, $\omega_Y\in \Det(Y)$,
$\omega_{(V,W)} \in \Det(V,W)$, where $X$, $Y$ are finite
dimensional subspaces of $H$, and $(V,W)$ is a Fredholm pair such
that $X\cap Y=(0)$, $(X+Y)\cap V=0$. The morphism of line bundles
${\textstyle \bigwedge} $ induces  a morphism of coverings,
denoted by the same symbol, between the orientation bundles:
\[
{\textstyle \bigwedge}  :\mathrm{Or}(\mathcal{S}(n,\mathrm{Fp})) \rightarrow
\mathrm{Or}(\mathrm{Fp}(H)),
\]
where the first space is the total space of the covering over
$\mathcal{S}(n,\mathrm{Fp})$ whose fiber at $(X,(V,W))$ is $\mathrm{Or}(X) \times
\mathrm{Or}(V,W)$. This map satisfies the identity
\[
o_X {\textstyle \bigwedge}  o_{(V,W)} = (-o_X) {\textstyle
\bigwedge} (-o_{(V,W)}) = - (-o_X) {\textstyle \bigwedge}
o_{(V,W)} = - o_X {\textstyle \bigwedge} (-o_{(V,W)}),
\]
and it is associative, meaning that
\[
o_X {\textstyle \bigwedge}  (o_Y {\textstyle \bigwedge} o_{(V,W)})
= (o_X \wedge o_Y) {\textstyle \bigwedge}  o_{(V,W)},
\]
for any $o_X\in \mathrm{Or}(X)$, $o_Y\in \mathrm{Or}(Y)$, $o_{(V,W)} \in \mathrm{Or}(V,W)$,
where $X$, $Y$ are finite dimensional subspaces of $H$, and
$(V,W)$ is a Fredholm pair such that $X\cap Y=(0)$, $(X+Y)\cap
V=0$.

\paragraph{The Grassmannian of compact perturbations.}
We shall say that the closed linear subspace $W$ is a {\em
compact perturbation} of $V$ if its orthogonal projector $P_W$ is a compact
perturbation of $P_V$. The subspace $W$ is a compact perturbation of
$W$ if and
only if the operators $P_{V^{\perp}} P_W$ and $P_{W^{\perp}} P_V$
are compact. The notion of compact perturbation produces an
equivalence  relation, and the
{\em Grassmannian of compact perturbations of} $V$,
\[
\mathrm{Gr}(V,H) := \set{W\in \mathrm{Gr}(H)}{W \mbox{ is a compact 
perturbation of }V}
\]
is a closed subspace of $\mathrm{Gr}(H)$. If $V$ has finite dimension
(respectively finite codimension), then
\[
\mathrm{Gr}(V,H) = \bigcup_{n\in \N} \mathrm{Gr}_{n,\infty}(H), \quad \Bigl(\mbox{
    resp. } = \bigcup_{n\in \N} \mathrm{Gr}_{\infty,n}(H) \Bigr).
\]
In the more interesting case, $V$ has both infinite dimension and
infinite codimension. In such a situation, $\mathrm{Gr}(V,H)$ is a closed
proper subset of $\mathrm{Gr}_{\infty,\infty}(H)$. It is an analytic Banach
manifold, although just a $C^0$ Banach submanifold of $\mathrm{Gr}(H)$.
This space is also called {\em restricted Grassmannian} by some
authors (see \cite{sw85,ps86,cjs95}).

If $W$ is a compact perturbation of $W$, then $(V,W^{\perp})$ is a
Fredholm pair, and the {\em relative dimension of $V$ with
respect to
 $W$} is the integer
\[
\dim(V,W) := \ind (V,W^{\perp}) = \dim V\cap W^{\perp} - \dim
V^{\perp} \cap W.
\]
When $V$ and $W$ are finite dimensional (resp.\ finite
codimensional), we have $\dim (V,W)=\dim V-\dim W$ (resp.\ $\dim
(V,W)=\codim W - \codim V$).

\begin{prop}{\em (\cite{ama03b}, Proposition 5.1)}
\label{indice} If $(V,Z)$ is a Fredholm pair and $W$ is a compact
perturbation of $V$, then $(W,Z)$ is a Fredholm pair, with
\[
\ind (W,Z) = \ind (V,Z) + \dim (W,V).
\]
\end{prop}

In particular, if $V,W,Y$ are compact perturbations of the same subspace,
\begin{equation}
\label{astar}
\dim(Y,V) = \dim (Y,W) + \dim (W,V).
\end{equation}
Nor the notion of compact perturbation, neither the relative dimension
depend on the choice of an equivalent inner product in $H$.

\begin{prop}\label{trasf}{\em (\cite{ama01}, Proposition 2.3)}
Let $H_1,H_2$ be Hilbert spaces and let $T,S\in
\mathcal{L}(H_1,H_2)$ be operators with closed range and compact
difference. Then $\ker T$ is a compact perturbation of $\ker S$, $\ran T$
is a compact perturbation of $\ran S$, and
\[
\dim(\ran T,\ran S) = - \dim (\ker T,\ker S).
\]
\end{prop}

\begin{prop}
\label{appu}
Let $T\in \mathrm{GL}(H)$, $V\in \mathrm{Gr}(H)$, and let $P$ be a 
projector onto $V$. 
Then $TV$ is a compact perturbation of $V$ if and only if the
operator $(I-P)TP$ is compact.
\end{prop} 

\proof
By choosing a suitable inner product on $H$, we may assume that
$P=P_V$ is an orthogonal projector. The operator $L:= TP +
{T^*}^{-1} (I-P)$ is invertible, and $P_{TV} = T P
L^{-1}$. Therefore, $TV$ is a compact perturbation of $V$ if and only
if the operator
\[
(P_{TV} - P_V) L = (I-P) T P - P {T^*}^{-1} (I-P) =: S
\]
is compact. 

Now, if $S$ is compact, so is $(I-P)TP = SP$. On the other hand,
since the set
\[
\set{X\in \mathrm{GL}(H)}{(I-P)XP \in \mathcal{L}_c(H)}
\]
is a subgroup of $\mathrm{GL}(H)$, if $(I-P)TP$ is compact so is
$(I-P)T^{-1}P$. Therefore,
\[
S = (I-P)T P - \left( (I-P)T^{-1} P \right)^*
\]
is compact.
\qed

Let $V\in \mathrm{Gr}_{\infty,\infty}(H)$. The connected components of
$\mathrm{Gr}(V,H)$ are the subsets
\[
\mathrm{Gr}_n(V,H) := \set{W\in \mathrm{Gr}(V,H)}{\dim (W,V)=n}, \quad n\in \Z.
\]
These components are pairwise diffeomorphic. Each of these
components has the homotopy type of $\mathrm{BO}(\infty)$, so the homotopy
groups of $\mathrm{Gr}(V,H)$ are those listed in (\ref{bott}).

We conclude this section with a result about the kernel of
semi-Fredholm operators.

\begin{prop}
\label{kersf} Let $A,B\in \mathcal{L}(H_1,H_2)$ be continuous
linear operators between Hilbert spaces, with finite-codimensional
range. Assume that the restrictions $A|_{\ker B}$ and $B|_{\ker
  A}$ are compact. Then $\ker A$ is a compact perturbation of $\ker
B$, the operator $AB^*\in \mathcal{L}(H_2)$ is Fredholm, and
\begin{equation}
\label{formu} \ind (AB^*) = \dim \coker B - \dim \coker A + \dim
(\ker A,\ker B).
\end{equation}
\end{prop}

\proof Since $A$ has closed range, there exists $S\in
\mathcal{L}(H_2,H_1)$ such that $SA=P_{(\ker A)^{\perp}}$. Then
$P_{(\ker A)^{\perp}} P_{\ker B}= SAP_{\ker B}$ is compact, and
symmetrically so is $P_{(\ker B)^{\perp}} P_{\ker A}$. Therefore
$\ker A$ is a compact perturbation of $\ker B$. Moreover,
$AP_{(\ker B)^{\perp}}=A-AP_{\ker
  B}$ is a compact perturbation of $A$, so it
has closed range $\ran (AP_{(\ker B)^{\perp}}) = \ran (AB^*)$.
Since $(AB^*)^*=BA^*$, the exactness of the sequences
\begin{eqnarray*}
0 \rightarrow \ker B^* \hookrightarrow \ker (AB^*)
\stackrel{B^*}{\longrightarrow} \ker A \cap (\ker B)^{\perp}
\rightarrow 0, \\
0 \rightarrow \ker A^* \hookrightarrow \ker (BA^*)
\stackrel{A^*}{\longrightarrow} \ker B \cap (\ker A)^{\perp}
\rightarrow 0,
\end{eqnarray*}
implies that $AB^*$ is Fredholm and that (\ref{formu}) holds. \qed

\paragraph{Essential Grassmannians.}
If $m\in \N$, we define the {\em $(m)$-essential Grassmannian}
$\mathrm{Gr}_{(m)}(H)$ to be the quotient space of $\mathrm{Gr}(H)$ by the
equivalence relation
\[
\set{(V,W)\in \mathrm{Gr}(H) \times \mathrm{Gr}(H)}{V \mbox{ is a compact
perturbation
    of } W \mbox{ and } \dim(V,W) \in m\Z}.
\]
The $(1)$-essential Grassmannian is simply called {\em essential
Grassmannian}. If $E\in \mathrm{Gr}_{(m)}(H)$ and $V\in \mathrm{Gr}(H)$ is
commensurable to the subspaces belonging to the equivalence class
$E$,
\[
\dim (V,E) := \dim (V,W), \quad W\in E,
\]
defines an element of $\Z/m\Z$.

The essential Grassmannian $\mathrm{Gr}_{(1)}(H)$ is homeomorphic to the
space of symmetric idempotent elements of the Calkin algebra
$\mathcal{L}(H)/\mathcal{L}_c(H)$, hence it inherits the structure
of a complete metric space, and of an analytic submanifold of the
Calkin algebra.

Every set $\mathrm{Gr}_{n,\infty}(H)$ or $\mathrm{Gr}_{\infty,n}(H)$, $n\in \N$,
represents an isolated point in $\mathrm{Gr}_{(0)}(H)$, which has thus
countably many isolated points. If $m\geq 1$, the sets
\[
\bigcup_{n\in m\Z+k} \mathrm{Gr}_{n,\infty}(H) \quad \mbox{and} \quad
\bigcup_{n\in m\Z+k} \mathrm{Gr}_{\infty,n}(H), \quad k=0,1,\dots,m-1,
\]
represent $2m$ isolated points in $\mathrm{Gr}_{(m)}(H)$. The remaining
part of $\mathrm{Gr}_{(m)}(H)$ is connected, being the quotient space of
$\mathrm{Gr}_{\infty,\infty}(H)$, and it is denoted by $\mathrm{Gr}_{(m)}^*(H)$.

The space $\mathrm{Gr}_{(0)}^*(H)$ is simply connected, while the
fundamental group of $\mathrm{Gr}_{(m)}^*(H)$ for $m\geq 1$ is infinite
cyclic. If $m\geq 1$ divides $k\in \N$, the natural projection
\[
\mathrm{Gr}_{(k)}^*(H) \rightarrow \mathrm{Gr}^*_{(m)}(H)
\]
is a covering map. It is the universal covering of $\mathrm{Gr}_{(m)}^*(H)$
if $k=0$, it induces the homomorphism $q\mapsto (k/m) q$ between
fundamental groups if $k\neq 0$. For $m=1$ we obtain a covering
map with a basis having the structure of an analytic Banach
manifold and of a complete metric space, hence the same structures
can be lifted to $\mathrm{Gr}_{(k)}(H)$, for any $k\neq 1$.

Finally, the natural projection
\begin{equation}
\label{fb} \mathrm{Gr}_{\infty,\infty}(H) \rightarrow \mathrm{Gr}_{(m)}^*(H)
\end{equation}
is a $C^0$ fiber bundle\footnote{Although the map (\ref{fb}) is
analytic, it has no differentiable trivializations.}. Its total
space is contractible, and its typical fiber is 
\[
\bigcup_{[n]\in \Z/m\Z} \mathrm{Gr}_n(V,H), \quad \mbox{where } V\in
\mathrm{Gr}_{\infty,\infty} (H),
\]
a disjoint union all of whose components have the homotopy type of
$\mathrm{BO}(\infty)$. 
Therefore, the exact homotopy sequence of a fibration
yields to the isomorphisms
\[
\pi_i(\mathrm{Gr}^*_{(m)}(H)) \cong \pi_{i-1}(\mathrm{Gr}(V,H)) = \begin{cases}
 \Z & \mbox{for } i\equiv 1,5 \mod
  8, \\ \Z_2 & \mbox{for } i\equiv 2,3 \mod
  8, \\ 0 & \mbox{for } i\equiv 0,4,6,7 \mod
  8, \end{cases}
\]
for $i\geq 2$.

\section{Appendix - Linear ordinary differential operators in
  Hilbert spaces}

This appendix summarizes some results about linear ordinary
differential operators in Hilbert spaces. See \cite{ama03} for a
detailed exposition (see also \cite{rs95} and \cite{lt03} for
related results in the framework of unbounded operators, and for
an extensive bibliography).

Let $H$ be a real Hilbert space. A bounded linear operator $L\in
\mathcal{L}(H)$ is said {\em hyperbolic} if its spectrum does not
meet the imaginary axis. In such a case, let $H = V^+(L) \oplus
V^-(L)$ be the $L$-invariant splitting of $L$ into closed
subspaces, corresponding to the decomposition of the spectrum of
$L$ into positive and negative real part.

\begin{prop} \label{abuno}
{\em (\cite{ama01}, Proposition 2.2)} Let
  $L,L^{\prime}\in \mathcal{L}(H)$ be hyperbolic operators. If
  $L^{\prime}$ is a compact perturbation of $L$, then
  $V^+(L^{\prime})$ is a compact perturbation of $V^+(L)$, and
  $V^-(L^{\prime})$ is a compact perturbation of $V^-(L)$.
\end{prop}

Let $A:[0,+\infty] \rightarrow \mathcal{L}(H)$ (resp.\
$A:[-\infty,0] \rightarrow \mathcal{L}(H)$) be a piecewise
continuous path such that $A(+\infty)$ (resp.\ $A(-\infty)$) is
hyperbolic. We shall denote by $X_A:[0,+\infty[\rightarrow \mathrm{GL}(H)$
(resp.\ $X_A:]-\infty,0] \rightarrow \mathrm{GL}(H)$) the solution of the
linear Cauchy problem $X_A^{\prime}(t)=A(t)X_A(t)$, $X_A(0)=I$.
The {\em
  linear stable space} of $A$ (resp. the {\em linear unstable space }
  of $A$) is the linear subspace of $H$
\[
 W^s_A = \set{\xi\in H}{\lim_{t\rightarrow +\infty} X_A(t)\xi=0},
 \quad \Bigl(\mathrm{resp}.\ W^u_A = \set{\xi\in H}{\lim_{t\rightarrow -\infty}
 X_A(t)\xi=0}\Bigr).
\]
The main properties of the linear stable space are listed in the
following:

\begin{thm} \label{abue}
{\em (\cite{ama03}, Proposition 1.2 and Theorems 2.1, 3.1)} Let
$A:[0,+\infty] \rightarrow \mathcal{L}(H)$ be a piecewise
continuous path such that $A(+\infty)$ is hyperbolic. Then
$W^s_A$ is a closed subspace of $H$, which depends continuously
on $A$ in the $L^{\infty}([0,+\infty[,\mathcal{L}(H))$ topology.
The following convergence results for $t\rightarrow +\infty$ hold:
\begin{enumerate}
\item $W^s_A$ is the only closed subspace $W$ such that $X_A(t)W$
  converges to $V^-(A(+\infty))$;
\item $\|X_A(t)|_{W^s_A} \|$ converges to 0 exponentially fast.
\end{enumerate}
The above limits are locally uniform in $A$, with respect to the
$L^{\infty}$ topology. Moreover, if $V\in \mathrm{Gr}(H)$ is a linear
supplement of $W^s_A$,
\begin{enumerate}
\setcounter{enumi}{2} \item $X_A(t)V$ converges to
$V^+(A(+\infty))$;

\item $\displaystyle{\inf_{\substack{\xi \in V\\ |\xi|=1}}}
|X_A(t)\xi|$ diverges exponentially fast.
\end{enumerate}
The above limits are locally uniform in $V\in \mathrm{Gr}(H)$, and in $A$,
with respect to the $L^{\infty}$ topology. Finally:
\begin{enumerate}
\setcounter{enumi}{4} \item $W^s_{-A^*} = (W^s_A)^{\perp}$.
\end{enumerate}
\end{thm}

The analogous statements for the linear unstable space can be
deduced from the above theorem, taking into account the identity
$X_A(t)=X_B(-t)$ for $B(t)=-A(-t)$. The following proposition
characterizes those paths $A$ for which the evolution of the
linear stable space remains in a fixed essential class:

\begin{prop}{\em (\cite{ama03}, Proposition 3.8)}
\label{utl} Let $A:[0,+\infty] \rightarrow \mathcal{L}(H)$ be a
piecewise continuous path such that $A(+\infty)$ is hyperbolic.
Let $V$ be a closed subspace of $H$, and let $P$ be a projector
onto $V$. Then the following statements are equivalent:
\begin{enumerate}
\item $X_A(t)W^s_A$ is a compact perturbation of $V$ for any
$t\geq 0$; \item $V^-(A(+\infty))$ is a compact perturbation of
$V$ and $[A(t),P]P$ is
  compact for any $t\geq 0$.
\end{enumerate}
\end{prop}

The proof of the above proposition is based on the following
fact: if $V$ is a closed linear subspace of $H$, then the
orthogonal projector $P(t)$ onto $X_A(t)V$ solves the Riccati
equation
\begin{equation}
\label{ricceq} P^{\prime}(t) = (I-P(t)) A(t) P(t) + P(t) A(t)^*
(I-P(t)),
\end{equation}
as shown in \cite{ama03}, formula (35).

Now let $A:[-\infty,+\infty] \rightarrow \mathcal{L}(H)$ be a
continuous path such that $A(-\infty)$ and $A(+\infty)$ are
hyperbolic. If $C^0_0(\R,H)$ (resp.\ $C^1_0(\R,H)$) denotes the
Banach space of continuous curves $u:\R\rightarrow H$ such that
$u(t)$ is infinitesimal (resp.\ $u(t)$ and $u^{\prime}(t)$ are
infinitesimal) for $t\rightarrow \pm \infty$, we can consider the
bounded linear operator
\[
F_A : C^1_0(\R,H) \rightarrow C^0_0(\R,H), \quad (F_A u)(t) =
u^{\prime}(t) - A(t)u(t).
\]
Its main properties are listed in the following:

\begin{thm} \label{abre}
{\em (\cite{ama03}, Theorem 5.1 and Remark 5.1)} Let
$A:[-\infty,+\infty] \rightarrow \mathcal{L}(H)$ be a continuous
path such that $A(-\infty)$ and $A(+\infty)$ are hyperbolic. Then:
\begin{enumerate}
\item $F_A$ has closed range if and only if the linear subspace
  $W^s_A+W^u_A$ is closed;
\item $F_A$ is surjective if and only if $W^s_A+W^u_A=H$; \item
$F_A$ is injective if and only if $W^s_A \cap W^u_A=(0)$; \item
$F_A$ is a Fredholm operator if and only if $(W^s_A,W^u_A)$ is a
  Fredholm pair, and in this case $\ind F_A = \ind (W^s_A,W^u_A)$.
\end{enumerate}
\end{thm}

It is easy to build examples of paths $A$ having two arbitrary
closed linear subspaces as linear stable space and linear
unstable space, so the above theorem shows that in general $F_A$
may not have closed range, and its kernel and cokernel may be
infinite dimensional. Even $F_A$ is Fredholm, $A(t)$ is
self-adjoint and invertible for any $t$, and
$A(-\infty)=A(+\infty)$, the operator $F_A$ may still have any
index. In the following proposition we exhibit such an example
with positive index. By Theorem \ref{abue} (v), we obtain an
example with negative index by considering the path
$B(t)=-A(t)$.

\begin{prop}
\label{costra} Let $H$ be a separable infinite dimensional real
Hilbert space. Let $H=H^- \oplus H^+$ be an orthogonal splitting,
with $H^-,H^+\in \mathrm{Gr}_{\infty,\infty}(H)$. For any $k\in \N$ there
exists $A\in C^{\infty}(\R, \mathrm{GL}(H) \cap \mathrm{Sym}(H))$ such that
$A(t)=P_{H^+} - P_{H^-}$ for $t\notin (0,1)$, and $W^s_A +
W^u_A=H$, $\dim W^s_A \cap W^u_A = k$. In particular, $F_A$ is
a surjective Fredholm operator of index $k$.
\end{prop}

\proof Let $W\subset H^+$ be a linear subspace of dimension $k$.
Since $\mathrm{Gr}_{\infty,\infty}(H)$ is connected, there exist closed
subspaces $V_0 = H^- \oplus W, V_1, \dots, V_{m-1},V_m=H^-$ in
$\mathrm{Gr}_{\infty,\infty}(H)$ with $\dist (V_{j-1},V_j) <1$ for any
$j\in \{1,\dots,m\}$ (in fact it is possible to take $m=4$).
Denote by $\mathcal{S}$ the open subset of $\mathrm{Sym}(H)$ consisting of
the invertible operators $A$ with $\sigma_{\mathrm{ess}} (A) \cap
\R^{\pm} \neq \emptyset$. By Proposition \ref{hr} we can find
operators $A_1,\dots,A_m$ in $\mathcal{S}$ such that $e^{A_j/m}
V_{j-1} = V_j$. Define the piecewise constant path $B:\R
\rightarrow \mathcal{S}$ as
\[
B(t) = \begin{cases} P_{H^+} - P_{H^-} & \mbox{for } t<0 \mbox{
or }
  t\leq 1, \\ A_j & \mbox{for } \frac{j-1}{m} \leq t < \frac{j}{m},
  \quad j\in \{1,\dots,m\}. \end{cases}
\]
Since $X_B(t) = e^{tA_j/m} e^{A_{j-1}/m} \dots e^{A_1/m}$ for
$(j-1)/m \leq t \leq j/m$, there holds
\[
X_B(1) (H^-\oplus W) = e^{A_m/m} \dots e^{A_1/m} V_0 = V_m =H^-.
\]
Since $\mathcal{S}$ is connected, there is a sequence $(B_n)
\subset C^{\infty}(\R,\mathcal{S})$ with $B_n(t)=P_{H^+} -
P_{H^-}$ for $t\notin (0,1)$, $(B_n)$ bounded in
$L^{\infty}(\R,\mathcal{L}(H))$, and $B_n \rightarrow B$ in
$L^1(\R,\mathcal{L}(H))$. By the identity
\[
X_A(t) = X_B(t) + \int_0^t X_B(t) X_B(\tau)^{-1} (A-B)(\tau)
X_A(\tau)\, d\tau,
\]
the sequence $(X_{B_n}(1))$ converges to $X_B(1)$, hence
\[
W^s_{B_n} = X_{B_n} (1)^{-1} W^s_{B_n(\cdot+1)} = X_{B_n} (1)^{-1}
W^s_{P_{H^+} - P_{H^-}} = X_{B_n}(1)^{-1} H^- \rightarrow
X_{B}(1)^{-1} H^- = H^- \oplus W.
\]
Moreover, $W^u_{B_n} = W^u_{P_{H^+} - P_{H^-}} = H^+$, so for $n$
large enough $A=B_n$ satisfies $W^u_A + W^s_A=H$ and $\dim W^s_A
\cap W^u_A=k$. \qed

\section{Appendix - Hyperbolic rest points}

This appendix summarizes some well known results about hyperbolic
dynamics. See \cite{shu87}.

\paragraph{Local statements.}
Let $F$ be a vector field of class $C^1$ defined on a neighborhood
$U$ of $0$ in the real Hilbert space $H$. We denote by $\Omega(F)$
the maximal subset of $\R\times U$ where the local flow of $F$,
i.e. the solution of
\[
\partial_t \phi (t,p) = F(\phi(t,p)), \quad \phi(0,p)=p,
\]
is defined. We assume that $0$ is a hyperbolic rest point for $F$,
meaning that $F(0)=0$ and $L:= DF(0)$ is a hyperbolic operator,
that is $\sigma(L) \cap i\R = \emptyset$. Let $H^u \oplus H^s$ be
the splitting of $H$ corresponding to the partition of the
spectrum of $L$ into the closed subsets $\sigma(L) \cap \set{z\in
\C}{\re z>0}$ and $\sigma(L) \cap \set{z\in \C}{\re z<0}$. By
$P^u$ and $P^s=I-P^u$ we shall denote the projections onto $H^u$
and $H^s$, and we shall often identify $H=H^u \oplus H^s$ with
$H^u \times H^s$.

There exists an equivalent inner product $\langle \cdot, \cdot
\rangle$ on $H$ with associated norm $\|\cdot\|$ which is {\em
  adapted} to $L$, meaning that $H^u$ and $H^s$ are orthogonal, and
\begin{equation}
\label{adapt} \langle L \xi,\xi \rangle \geq \lambda \|\xi\|^2
\;\;\; \forall \xi \in H^u, \quad \langle L \xi,\xi \rangle \leq
-\lambda \|\xi\|^2 \;\;\; \forall \xi \in H^s,
\end{equation}
for some $\lambda>0$. Indeed, we may choose any positive $\lambda$
which is strictly less than $\min |\re \sigma(L)|$, as shown by
the following lemma, applied to $L|_{H^s}$ and to $-L|_{H^u}$.

\begin{lem}
Let $L$ be a bounded linear operator on $H$ and let $\lambda$ be a
real number such that $\lambda> \max \re \sigma(L)$. Then there
exists an equivalent inner product $\langle \cdot, \cdot \rangle$
on $H$ such that
\[
\langle L \xi,\xi \rangle \leq \lambda \langle \xi,\xi \rangle
\quad \forall \xi \in H.
\]
\end{lem}

\proof Up to replacing $L$ by $L-\lambda I$, we may assume that
$\lambda=0$. Let $\langle \cdot,\cdot \rangle_*$ be any Hilbert
product on $H$, and denote by $\|\cdot\|_*$ both the associated
norm on $H$ and the induced norm on $\mathcal{L}(H)$. By the
spectral radius formula and by the spectral mapping theorem,
\[
\lim_{n\rightarrow \infty} \|e^{nL}\|_*^{1/n} = \max
|\sigma(e^L)| = \max | e^{\sigma(L)} | < 1.
\]
Let $k\in \N$ be so large that $\|e^{kL}\|_*\leq 1$, and set
\[
\langle \xi, \eta \rangle := \int_0^k \langle e^{tL} \xi, e^{tL}
\eta \rangle_* \, dt, \quad \forall \xi,\eta\in H.
\]
Then $\langle \cdot,\cdot \rangle$ is an equivalent inner product
on $H$, and for any $\xi\in H$
\begin{eqnarray*}
\langle L\xi,\xi \rangle = \int_0^k \langle e^{tL} L \xi, e^{tL}
\xi \rangle_* \, dt = \frac{1}{2} \int_0^k \frac{d}{dt} \|e^{tL}
\xi\|_*^2\, dt \\
= \frac{1}{2} \left( \|e^{kL} \xi\|^2_* - \|\xi\|_*^2 \right) \leq
\frac{1}{2} \left( \|e^{kL}\|_*^2 - 1\right) \|\xi\|_*^2 \leq 0,
\end{eqnarray*}
concluding the proof. \qed

If $V$ is a closed linear subspace of $H$ and $r>0$, $V(r)$ will
denote the closed ball of $V$ centered in 0 with radius $r$.
Moreover, we set
\[
Q(r) := \set{\xi \in H}{\|P^u \xi\|\leq r, \; \|P^s \xi\|\leq r}.
\]
If $A\subset X\subset H$, the set $A$ is said {\em positively
  (negatively) invariant with respect to $X$} if for every $\xi\in A$
  and for every $t>0$,
  $\phi([0,t] \times \{\xi\})\subset X$ implies $\phi([0,t] \times
  \{\xi\})\subset A$ (resp.\ for every $\xi\in A$
  and for every $t<0$, $\phi([t,0] \times \{\xi\})\subset X$
  implies $\phi([t,0] \times \{\xi\})\subset A$).

\begin{lem}
\label{llmm} For any $r>0$ small enough, the set
\[
\set{\xi\in Q(r)}{\|P^s \xi\| \leq \|P^u \xi\|} \quad
\mbox{(resp.\ } \set{\xi\in Q(r)}{\|P^u \xi\| \leq \|P^s
\xi\|}\mbox{)}
\]
is positively (resp.\ negatively) invariant with respect to
$Q(r)$. Moreover, if $\xi$ belongs to the set
\[
\set{\xi \in Q(r)}{\|P^u \xi\| =r}\quad
\mbox{(resp.\ } \set{\xi \in Q(r)}{\|P^s \xi\|
=r}\mbox{)},
\]
then $\phi(t,\xi)\notin Q(r)$ (resp.\ $\phi(-t,\xi)\notin Q(r)$)
for every $t>0$ small enough.
\end{lem}

\proof By a first order expansion of $F$ at $0$ and by
(\ref{adapt}),
\begin{eqnarray*}
\frac{d}{dt} \left\| P^u \phi(t,\xi) \right\|^2 \Bigr|_{t=0} \geq
2 \lambda \|P^u \xi\|^2 + o(\|P^u \xi\|^2) \quad \mbox{if } \|P^s
\xi\| \leq \|P^u \xi\|, \\ \frac{d}{dt} \left\| P^s \phi(t,\xi)
\right\|^2 \Bigr|_{t=0} \leq -2 \lambda \|P^s \xi\|^2 + o(\|P^s
\xi\|^2) \quad \mbox{if } \|P^u \xi\| \leq \|P^s \xi\|.
\end{eqnarray*}
All the statements follow from the above inequalities. \qed

Given $r>0$, the {\em local unstable manifold} and the {\em local
  stable manifold} of $0$ are the sets
\begin{eqnarray*}
W^u_{\mathrm{loc},r}(0) = \set{\xi \in Q(r)}{]-\infty,0]\times
\{\xi\}
      \subset \Omega(F) \mbox{ and }
      \phi(]-\infty,0] \times \{\xi\})\subset Q(r)}, \\
W^s_{\mathrm{loc},r}(0) = \set{\xi \in Q(r)}{[0,+\infty[\times
\{\xi\}
      \subset \Omega(F) \mbox{ and } \phi([0,+\infty[ \times \{\xi\})
      \subset Q(r)}.
\end{eqnarray*}
Then the local stable manifold theorem (see \cite{shu87}, chapter
5) states that:
\begin{thm}
\label{stabile} For any $r>0$ small enough,
$W^u_{\mathrm{loc},r}(0)$ (respectively $W^s_{\mathrm{loc},r}(0)$)
is the graph of a $C^1$ map $\sigma^u: H^u(r) \rightarrow H^s(r)$
such that $\sigma^u(0)=0$ and $D\sigma^u(0)=0$ (resp.\ of a $C^1$
map $\sigma^s: H^s(r) \rightarrow H^u(r)$ such that
$\sigma^s(0)=0$ and $D\sigma^s(0)=0$). Moreover, for any $\xi\in
W^u_{\mathrm{loc},r}(0)$ (resp.\ for any $\xi\in
W^s_{\mathrm{loc},r}(0)$), there holds
\[
\lim_{t\rightarrow -\infty} \phi(t,\xi) = 0 \quad ( \mbox{resp. }
\lim_{t\rightarrow +\infty} \phi(t,\xi) = 0 ).
\]
\end{thm}

We recall that a non-degenerate local Lyapunov function for the
vector field $F$ at the rest point $0$ is a $C^1$ real function
defined on a neighborhood of 0 in $H$, such that $Df(\xi) [F(\xi)]
<0$ for $\xi\neq 0$, and which is twice differentiable at $0$,
with the quadratic form $D^2 f(0)$ coercive on $H^u$, and the
quadratic form $-D^2 f(0)$ coercive on $H^s$ (necessarily,
$Df(0)=0$). A first order expansion of $F$ at $0$ shows that the
restriction of the function
\[
f(\xi) = -\frac{1}{2} \langle L\xi,\xi \rangle
\]
to a suitably small neighborhood of $0$ is a non-degenerate local
Lyapunov function for $F$ at $0$.

\begin{lem}
\label{hg} For any $r>0$ small enough, for every sequence
$(\xi_n)\subset H$ converging to 0 and for every sequence
$(t_n)\subset [0,+\infty[$ such that $\phi([0,t_n] \times
\{\xi_n\}) \subset Q(r)$ and $\phi(t_n,\xi_n)\in
\partial Q(r)$, there holds
\[
\dist\left(\phi(t_n,\xi_n),W^u_{\mathrm{loc},r}(0)\cap \partial
Q(r) \right) \rightarrow 0.
\]
Furthermore, if $f$ is a non-degenerate local Lyapunov function
for $F$ at $0$, there holds
\[
\limsup_{n\rightarrow \infty} f(\phi(t_n,\xi_n)) < f(0).
\]
Finally, there exists $r^{\prime}<r$ such that
\begin{eqnarray*}
\sup \set{f(\xi)}{\xi \in \partial Q(r) \mbox{ and } \exists t<0
  \mbox{ such that } \phi(-t,\xi)\in Q(r^{\prime}), \; \phi([-t,0]
  \times \{\xi\}) \subset Q(r)} \\ < \inf \set{f(\xi)}{\xi \in Q(r^{\prime})}.
\end{eqnarray*}
\end{lem}

\proof If the vector field is linear, $F(\xi)=L\xi$, the first
conclusion is immediate: actually for any $(\xi_n)\subset H$
converging to $0$ and any $(t_n)\subset [0,+\infty[$, there holds
\begin{equation}
\label{llin} \lim_{n\rightarrow \infty} \dist \left(e^{t_n L}
\xi_n, H^u\right)=0.
\end{equation}
By the Grobman-Hartman theorem, if $r_1>0$ is small enough the
local flow $\phi$ restricted to $Q(r_1)$ is conjugated to its
linearization $(t,\xi) \mapsto e^{tL} \xi$, by a bi-uniformly
continuous homeomorphism\footnote{We recall that this conjugacy
is found as a
  fixed point of a contraction $T$ on a suitable space of continuous
  maps (see \cite{shu87}, chapter 7). Since for $\alpha\in ]0,1[$ small enough,
  the space of $\alpha$-H\"older continuous maps is $T$-invariant,
  such a conjugacy turns out to be  H\"older
  continuous together with its inverse. In general, it needs not be
  even Lipschitz continuous.}.
By Theorem \ref{stabile}, we may also assume that $r_1$ is so
small that $W^u_{\mathrm{loc},r_1}(0)$ is the graph of a uniformly
continuous map $\sigma^u :H^u(r_1) \rightarrow H^s(r_1)$.

Let $r<r_1$ and set $\eta_n:=\phi(t_n,\xi_n) \in \partial Q(r)$,
with $\xi_n \rightarrow 0$ and $t_n\geq 0$. By the linear case
(\ref{llin}) and by the uniform continuity of the conjugacy,
there exists $(\eta_n^{\prime}) \subset W^u_{\mathrm{loc},r_1}
(0)$ such that $\|\eta_n - \eta_n^{\prime}\|$ is infinitesimal.
Setting $\eta_n^{\prime\prime}=(P^u \eta_n, \sigma^u(P^u
\eta_n))\in W^u_{\mathrm{loc},r}(0)\cap \partial Q(r)$, by the
uniform continuity of $\sigma^u$ we have
\begin{eqnarray*}
\dist\left(\eta_n, W^u_{\mathrm{loc},r}(0)\cap \partial
Q(r)\right) \leq \|\eta_n - \eta_n^{\prime\prime}\| \leq \|\eta_n
- \eta_n^{\prime}\| + \|P^u \eta_n^{\prime} - P^u
\eta_n^{\prime\prime}\| + \|P^s \eta_n^{\prime} - P^s
\eta_n^{\prime\prime}\| \\ =  \|\eta_n - \eta_n^{\prime}\|  +
\|P^u \eta_n^{\prime} - P^u \eta_n\| + \|\sigma^u(P^u
\eta_n^{\prime}) - \sigma^u(P^u \eta_n)\| \rightarrow 0,
\end{eqnarray*}
proving the first claim. Since the local unstable manifold is
tangent to $H^u$ at $0$, since $Df(0)=0$ and $-D^2 f(0)$ is
coercive on $H^u$, by $o(r)$ considerations we have
\[
\sup \set{f(\xi)}{\xi \in W^u_{\mathrm{loc},r}(0)\cap \partial
Q(r)} < f(0),
\]
if $r>0$ is small enough.
Since $f$ is uniformly continuous on $Q(r)$ for $r$ small enough,
the second claim follows from the first one. The last claim is an
immediate consequence of the second one, arguing by contradiction. 
\qed

Given two metric spaces $X$ and $Y$ and a positive number
$\theta$, $\mathrm{Lip}_{\theta}(X,Y)$ will denote the space of
$\theta$-Lipschitz maps from $X$ to $Y$, endowed with the $C^0$
topology. The following version of the graph transform theorem is
proved in \cite{ama01} Proposition A.3 and Addendum A.5 (see also
\cite{shu87}, chapter 5).

\begin{prop}
\label{grafici} For any $r>0$ small enough there is a continuous
(nonlinear) semigroup
\[
\Gamma: [0,+\infty] \times \mathrm{Lip}_1 (H^u(r), H^s(r))
\rightarrow \mathrm{Lip}_1 (H^u(r), H^s(r))
\]
such that for every $\sigma\in \mathrm{Lip}_1 (H^u(r), H^s(r))$
there holds:
\begin{enumerate}
\item $\Gamma(0,\sigma) = \sigma$, and $\Gamma(t+s,\sigma) =
  \Gamma(t,\Gamma(s,\sigma))$, for every $t,s\in [0,+\infty]$;
\item for every $t\in [0,+\infty[$, the restriction of $\phi_t$
to
  $Q(r)$ maps the graph of $\sigma$ onto the graph of
  $\Gamma(t,\sigma)$, that is
\[
\graf \Gamma(t,\sigma) = \set{\phi(t,\xi)}{\xi \in \graf \sigma
  \mbox{ and } \phi([0,t]\times \{\xi\}) \subset Q(r)};
\]
\item $\graf \Gamma(+\infty,\sigma) = W^u_{\mathrm{loc},r}(0)$;
\item for any $\theta>0$ there exists $r_0\in ]0,r]$ and $t_0\geq
0$ such that the restriction of $\Gamma(t,\sigma)$ to $H^u(r_0)$
is in $\mathrm{Lip}_{\theta}(H^u(r_0),H^s(r_0))$, for any $t\in
[t_0,+\infty]$ and any $\sigma\in \mathrm{Lip}_1(H^u(r),H^s(r))$.
\end{enumerate}
Furthermore:
\begin{enumerate}
\setcounter{enumi}{4} \item if $V \subset H^u(r)$ is open and
$\sigma\in \mathrm{Lip}_1(V,H^s(r))$ is such that $\graf \sigma
\cap W^s_{\mathrm{loc},r}(0) \neq \emptyset$, then there exists
$t\geq 0$ and $\sigma^{\prime}\in \mathrm{Lip}_1(H^u(r),H^s(r))$
such that the restriction of $\phi_t$ to $Q(r)$ maps the graph of
$\sigma$ onto the graph of $\sigma^{\prime}$, that is
\[
\graf \sigma^{\prime} = \set{\phi(t,\xi)}{\xi \in \graf \sigma
  \mbox{ and } \phi([0,t]\times \{\xi\}) \subset Q(r)}.
\]
\end{enumerate}
\end{prop}

\paragraph{Global statements.}
Now let $F$ be a $C^1$ vector field on the real Hilbert manifold
$M$, and let $\phi:\Omega(F) \rightarrow M$, $\Omega(F)\subset \R
\times M$, denote its local flow. Let $x$ be a hyperbolic rest
point of $F$. We can identify a neighborhood of $x$ in $M$ with a
neighborhood of $0$ in the Hilbert space $H$, identifying $x$ with
$0$. We denote by
$H=H^u \oplus H^s$ the splitting of $H$ associated to the hyperbolic
operator $\nabla F(x) = DF(0)$, and we
endow $H$ with an equivalent inner product adapted to $\nabla F(x)$,
as in the previous section. For $r>0$ small enough, we set
\[
Q(r) = H^u(r) \times H^s(r), \quad
Q^+(r) = \partial H^u(r)\times
H^s(r), \quad Q^-(r) = H^u(r)\times \partial H^s(r).
\]
Lemma \ref{llmm} and the last statement of Lemma \ref{hg} have
the following consequence.

\begin{prop}
\label{conley} For any $r>0$ small enough there holds:
\begin{enumerate}
\item if $p\in Q(r)$ and $\phi(t,p)\notin Q(r)$ for some $t>0$, then
  there exists $s\in [0,t[$ such that $\phi(s,p)\in Q^+(r)$;
\item if $p\in Q(r)$ and $\phi(t,p)\notin Q(r)$ for some $t<0$, then
  there exists $s\in ]t,0]$ such that $\phi(s,p)\in Q^-(r)$.
\end{enumerate}
Moreover, if $F$ admits a global $C^1$ Lyapunov function which is
twice differentiable and non-degenerate at $x$:
\begin{enumerate}
\setcounter{enumi}{2} \item if $p\in Q^+(r)$, then $\phi(t,p)\notin
Q(r)$ for any $t>0$; \item if $p\in Q^-(r)$, then $\phi(t,p)\notin
Q()$ for any $t<0$.
\end{enumerate}
\end{prop}

\noindent The unstable and stable manifolds of $x$ are the
$\phi$-invariant subsets of $M$
\begin{eqnarray*}
W^u(x) = \set{p \in M}{]-\infty,0]\times \{p\} \subset \Omega(F)
  \mbox{ and } \lim_{t\rightarrow -\infty} \phi(t,p) = x}, \\
W^s(x) = \set{p \in M}{[0,+\infty[ \times \{p\} \subset \Omega(F)
  \mbox{ and } \lim_{t\rightarrow +\infty} \phi(t,p) = x}.
\end{eqnarray*}
The local stable manifold theorem (Theorem \ref{stabile}) and
Proposition \ref{conley} imply:
\begin{thm}
\label{thesta}
The sets $W^u(x)$ and $W^s(x)$ are images of $C^1$
injective immersions
\[
e^u : H^u \rightarrow M, \quad e^s : H^s \rightarrow M,
\]
such that $e^u(0)=e^s(0)=x$, and $D e^u(0)$ and $D e^s(0)$ are the
identity mappings. If moreover $F$ admits a global $C^1$ Lyapunov
function which is twice differentiable and non-degenerate at $x$,
then for any $r>0$ small enough
\[
W^u(x) \cap Q(r) = W^u_{\mathrm{loc},r}(0), \quad W^s(x)
\cap Q(r) = W^s_{\mathrm{loc},r}(0),
\]
and the maps $e^u$, $e^s$ are embeddings, so that $W^u(x)$ and
$W^s(x)$ are $C^1$ submanifolds of $M$.
\end{thm}

\providecommand{\bysame}{\leavevmode\hbox to3em{\hrulefill}\thinspace}
\providecommand{\MR}{\relax\ifhmode\unskip\space\fi MR }
\providecommand{\MRhref}[2]{%
  \href{http://www.ams.org/mathscinet-getitem?mr=#1}{#2}
}
\providecommand{\href}[2]{#2}

\end{document}